\documentstyle{article}

\def\A{{\cal A}}
\def\C{{\cal C}}
\def\G{{\cal G}}
\def\X{{\cal X}}
\def\H{{\cal H}}
\def\P{{\cal P}}

\def\pl{\!+\!}
\def\mn{\!-\!}

\def\prop#1#2{\vspace{2ex} \noindent{\sc #1.} {\it #2} \par \vspace{2ex}}
\def\propo#1#2{\vspace{2ex} \noindent{\sc #1.} {\it #2}}
\def\oprop#1#2{\noindent{\sc #1.} {\it #2} \par \vspace{2ex}}
\def\dkz{\noindent{\sc Proof. }}
\def\qed{\hfill $\dashv$\vspace{2ex}}

\def\uhy{_{\raisebox{-.5pt}{$\cup$} H-Y}}
\def\ccup{\raisebox{-.5pt}{$\cup$}}
\def\R{{\mathbf R}}

\begin{document}

\title{\textbf{Hypergraph Polytopes}}
\author{{\sc Kosta Do\v sen} and {\sc Zoran Petri\' c}
\\[1ex]
{\small Mathematical Institute, SANU}\\[-.5ex]
{\small Knez Mihailova 36, p.f.\ 367, 11001 Belgrade,
Serbia}\\[-.5ex]
{\small email: \{kosta, zpetric\}@mi.sanu.ac.rs}}
\date{}
\maketitle

\begin{abstract}
\noindent We investigate a family of polytopes introduced by E.M.\
Feichtner, A.\ Postnikov and B.\ Sturmfels, which were named
nestohedra. The vertices of these polytopes may intuitively be
understood as constructions of hypergraphs. Limit cases in this
family of polytopes are, on the one end, simplices, and, on the
other end, permutohedra. In between, as notable members one finds
associahedra and cyclohedra. The polytopes in this family are
investigated here both as abstract polytopes and as realized in
Euclidean spaces of all finite dimensions. The later realizations
are inspired by J.D.\ Stasheff's and S.\ Shnider's realizations of
associahedra. In these realizations, passing from simplices to
permutohedra, via associahedra, cyclohedra and other interesting
polytopes, involves truncating vertices, edges and other faces.
The results presented here reformulate, systematize and extend
previously obtained results, and in particular those concerning
polytopes based on constructions of graphs, which were introduced
by M.\ Carr and S.L.\ Devadoss.
\end{abstract}

\vspace{.5ex}


\vspace{2ex}

\begin{tabbing}
{\bf Contents}
\\*[1.5ex]
1. \hspace{.3em} \= {\it Introduction}
\\*
2. \> {\it Connected hypergraphs}
\\
3. \> {\it Constructions}
\\
4. \> {\it Saturation and cognate hypergraphs}
\\
5. \> {\it Constructs and abstract polytopes of hypergraphs}
\\
6. \> {\it Constructions of ASC-hypergraphs}
\\
7. \> {\it Continuations of constructions}
\\
8. \> {\it Abstract polytopes of hypergraphs are abstract
polytopes}
\\
9. \> {\it Realizations}
\\[1.5ex]
Appendix A. \hspace{.3em} \= {\it Constructs and tubings}
\\
Appendix B. \> {\it Hypergraph polytopes of dimension 3 and lower}
\\[1.5ex]
{\it References}
\end{tabbing}

\section{\large \textbf{Introduction}}
One key to understanding the permutohedron is that it is a
truncated simplex. Our results here are a development of that
idea. They present the abstract underpinnings of these
truncations.

We investigate a family of polytopes that like permutohedra may be
obtained by truncating the vertices, edges and other faces of
simplices, in any finite dimension. The permutohedra are limit
cases in that family, where all possible truncations have been
made. The limit cases at the other end, where no truncation has
been made, are simplices, like the tetrahedron in three
dimensions.

As notable intermediate cases, with some truncations, whose
principles we make manifest, we have in this family associahedra
and cyclohedra (see Appendix~B; for historical references
concerning associahedra, cyclohedra and permutohedra see
\cite{S97a}, \cite{S97b}, \cite{Z95}, Lecture~0, Example 0.10, and
\cite{GR63}). There are also other interesting polytopes in the
vicinity of these, which are not very well known, or are quite
unknown. These other polytopes have an application in category
theory and the theory of operads similar to that of associahedra
(see \cite{DP10b}, and the end of this introduction).

This family of polytopes was introduced as a family in \cite{FS05}
(Section~3) and \cite{P09} (Section~7). The polytopes in it were
named nestohedra in \cite{PRW08} (Section~6).

The polytopes in this family are defined with respect to
hypergraphs (in the sense of \cite{B89}; see the beginning of
Section~1). These hypergraphs are essentially a special kind of
building sets, which are defined in \cite{FK04} with respect to
arbitrary finite meet semilattices. For hypergraphs we have
instead finite set lattices, with meet and join being respectively
intersection and union (see \cite{FS05} and \cite{P09}).

A vertex of one of these polytopes may be identified with another
hypergraph that may intuitively be understood as a construction of
the original hypergraph. The faces of greater dimension of the
polytope correspond to partial constructions. These constructions,
partial and not partial, which we call constructs, are called
nested sets in \cite{FS05} and \cite{P09} (a term that was
introduced with respect to the more general notion of building set
in \cite{FK04}; references concerning antecedents of this notion
are in \cite{FK04}, beginning of the Introduction, and
\cite{PRW08}, beginning of Section~6).

In most of our text, we deal with the hypergraph polytopes in an
abstract manner, based on the definition of abstract polytope of
\cite{MS02} (Section 2A). We believe that this point of view is
novel.

We devote however one part of our work (Section~9) to the
Euclidean realizations of hypergraph polytopes. This approach to
realizing these polytopes, which is inspired by \cite{S97b}
(Appendix~B), and is based on truncating simplices, may be found
in the literature in cases where the hypergraphs can be identified
with graphs (see \cite{CD06}, \cite{D09}, \cite{DF08} and further
references that may be found in these papers; the polytopes in
question are called there graph-associahedra). Here this approach
is extended to all hypergraph polytopes. The approach to realizing
hypergraph polytopes of \cite{FS05} and \cite{P09}, which is based
on the Minkowsky sum of simplices, is different (see, however,
also Remark 6.6 of \cite{PRW08}, \cite{V10}, \cite{BV10} and
references in there).

Another difference of our approach is that for us inductive
definitions play a more important role than it is the case in the
other approaches. We find that these definitions enable us to
clarify and simplify matters; for proving some results they have a
clear advantage. We present several alternative views on the same
subject matter---in particular, three equivalent notions of
constructions. (These notions are closely related to notions in
\cite{FS05}, end of Section~3, and \cite{P09}, Definition 7.7; we
studied them first in \cite{DP10a}.) When we restrict ourselves to
graphs, then the notion of construct, which is for us a secondary
notion, derived from the primary, more basic, notion of
construction, amounts to the notion of tubing of \cite{DF08} (see
Appendix~A, where these matters are treated in detail; this
appendix provides a bridge between the approach through nested
sets and the approach through tubings).

Another novelty that we give may be an inductive definition of an
abstract polytope (see Section~8), equivalent to the definition of
\cite{MS02}. We find this inductive definition useful for showing
that abstract hypergraph polytopes are indeed abstract polytopes.
For that, we rely on the results of Section~7, which are closely
related to the results of \cite{Z06}, though the presentation is
different.

We survey all the hypergraph polytopes up to and including
dimension 3 in Appendix~B. To obtain intuitive pictures, the
reader may consult this appendix while going through the previous
exposition.

Our investigation of the matters covered here started in
\cite{DP10a}. In \cite{DP10b}, which is about a problem in
category theory and the theory of operads, one finds an
application of the ideas of \cite{DP10a}. In general, in these two
preceding papers we were less concerned with the theory of
polytopes, abstract or realized.

In \cite{DP10a} we worked in the direction from the permutohedra
towards other hypergraph polytopes, which is the direction of
\cite{T97} and \cite{P09}. We could not reach simplices, because
we stuck to graphs only, and did not envisage hypergraphs. We were
collapsing different vertices of a permutohedron into a single
vertex (which is akin to what is done in \cite{T97}). Now we work
in the opposite direction, by truncating, starting from the
simplices, as in \cite{S97b} (Appendix~B). The two approaches,
with two opposite directions, collapsing and truncating, cover
however essentially the same ground (provided that by introducing
hypergraphs we allow collapsing to go all the way up to
simplices). They have an identical basic core, and our goal here
is to present clearly this core.

As one can base an alternative proof of Mac Lane's coherence
theorem for monoidal categories of \cite{ML63} (see also
\cite{ML98}, Section VII.2) on Stasheff's results of \cite{S63}
concerning associahedra (see also \cite{S97a}, \cite{S97b} and
references therein), so one can base an alternative proof of a
categorial coherence result of \cite{DP10b}, concerning operads,
on the results presented here. We will however deal with these
matters of category theory on another occasion.

The first version of this paper posted in the arXiv, which differs
unessentially from the present one, was written without our being
aware of \cite{FK04}, \cite{FS05}, \cite{P09}, \cite{PRW08} and
\cite{Z06}. We were also not aware of the papers \cite{BV10},
\cite{V10} and \cite{C07}, dealing with matters related to our
truncations. We would probably have presented matters differently
if we knew about these references from the outset, but perhaps our
independent approach, for which we believe that it is sometimes
simpler, sheds a new light on the matter.

\section{\large \textbf{Connected hypergraphs}}
In this section we define the basic notions that we need
concerning hypergraphs.

For $C$ a finite (possibly empty) set, consider families of sets
$H$ such that $H\subseteq \P C$, i.e.\ families of subsets of $C$.
When $\emptyset\not\in H$ and $C$ is the union $\bigcup H$ of all
the members of $H$, the family $H$ is a \emph{hypergraph} on $C$
(see \cite{B89}, Section 1.1). The members of $C$ correspond to
the vertices of a graph, and the members of $H$ that are pairs,
i.e.\ two-element sets, correspond to the edges of a graph. It is
not necessary to mention always the \emph{carrier} $C$ of a
hypergraph, since $C=\bigcup H$, and every hypergraph $H$ is a
hypergraph on $\bigcup H$; sometimes however mentioning $C$ is
useful, and clarifies matters. A hypergraph on $\{x,y,z,u,v\}$ is,
for example, the family
\[E=\{\{x,y\},\{x,y,z\},\{y,z\},\{u\},\{v\}\}.\]

The \emph{empty hypergraph} is the hypergraph $\emptyset$ on
$\emptyset$. There is no other hypergraph on $\emptyset$, since
$\{\emptyset\}$ is not a hypergraph on $\emptyset$; we have that
$\emptyset\in\{\emptyset\}$. So $H=\emptyset$ iff $\bigcup
H=\emptyset$.

We allow the empty hypergraph, and spend some time in explaining
limit matters pertaining to it, but the reader should not imagine
that this is extremely important. In much of our text the empty
hypergraph fits nicely into the picture, but in some parts (see
the end of Section~7 and the beginning of Section~9) we treat it
separately. Our main interest is in nonempty hypergraphs, and the
limit case of the empty hypergraph could as well have been
omitted.

A \emph{hypergraph partition} of a hypergraph $H$ is a partition
$\{H_1,\ldots,H_n\}$, with $n\geq 1$, of $H$ such that $\{\bigcup
H_1,\ldots,\bigcup H_n\}$ is a partition of $\bigcup H$. For
example, the sets
\begin{tabbing}
\hspace{7em}\=$E'$\hspace{.5em}\=$=\{\{\{x,y\},\{x,y,z\},
\{y,z\}\},\{\{u\}\},\{\{v\}\}\}$,\\*
\>$E''\;$\>$=\{\{\{x,y\},\{x,y,z\},\{y,z\},\{u\}\},\{\{v\}\}\}$,\\*
\>$E'''\;$\>$=\{E\}$
\end{tabbing}
are hypergraph partitions of $E$. The partition
\[\{\{\{x,y\},\{x,y,z\}\},\{\{y,z\},\{u\},\{v\}\}\}\]
of $E$ is not a hypergraph partition of $E$, because we have that
$\bigcup \{\{x,y\},\{x,y,z\}\}$ $=\{x,y,z\}$ and
$\bigcup\{\{y,z\},\{u\},\{v\}\}=\{y,z,u,v\}$, and
$\{\{x,y,z\},\{y,z,u,v\}\}$ is not a partition of $\{x,y,z,u,v\}$.

The \emph{trivial} partition $\{H\}$ of $H$ exists if $H$ is
nonempty, and it is a hypergraph partition. If $H=\emptyset$, then
$\{\emptyset\}$ is not a partition of $\emptyset$, because all the
members of a partition must be nonempty. The empty hypergraph has
however one, and only one, partition; this is the empty partition
$\emptyset$, which is a hypergraph partition of~$\emptyset$.

A hypergraph $H$ is \emph{connected} when it has only one
hypergraph partition; if $H$ is nonempty, then this unique
hypergraph partition is the trivial partition $\{H\}$, and if
$H=\emptyset$, then this hypergraph partition is $\emptyset$. The
hypergraph $E$ above is not connected; the family
\[\{\{x,y\},\{x,y,z\},\{y,z\},\{z,u\}\}\]
is a connected hypergraph on $\{x,y,z,u\}$.

For a hypergraph $H$, let the \emph{intersection graph} of $H$ be
the graph $\Omega(H)$ whose vertices are the elements of $H$,
which are connected by an edge when they have a nonempty
intersection (see \cite{H69}, Chapter~2). For example, $\Omega(E)$
is
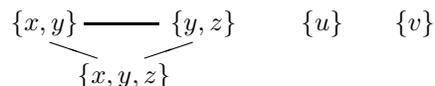
\begin{center}
\begin{picture}(140,30)
\put(20,12){\line(-3,1){18}}

\put(38,12){\line(3,1){18}}

\put(15,25){\line(1,0){28}}

\put(30,5){\makebox(0,0){$\{x,y,z\}$}}
\put(0,25){\makebox(0,0){$\{x,y\}$}}

\put(60,25){\makebox(0,0){$\{y,z\}$}}

\put(105,25){\makebox(0,0){$\{u\}$}}

\put(140,25){\makebox(0,0){$\{v\}$}}
\end{picture}
\end{center}

A \emph{path} of a hypergraph $H$ is a sequence $X_1,\ldots,X_n$,
with $n\geq 1$, of distinct elements of $H$ that make a path in
$\Omega(H)$ (see \cite{H69}, Chapter~2, for the notion of path in
a graph; this is a sequence of distinct vertices such that
consecutive vertices are joined by edges). Three examples of paths
in $E$ are the sequences
\[
\begin{array}{l}
\{x,y\},\{y,z\},\\
\{y,z\},\{x,y\},\{x,y,z\},\\
\{x,y\}
\end{array}
\]
If $n=1$, then the path $X_1$ is just the element $X_1$ of $H$. In
the third example, the element $\{x,y\}$ of $E$ is a path of~$E$.

For $x,y\in\bigcup H$, we say that a path $X_1,\ldots,X_n$ of $H$
\emph{joins} $x$ with $y$ when $x\in X_1$ and $y\in X_n$. So the
path $\{x,y\},\{y,z\}$ of $E$ joins $x$ with $z$, but it also
joins $x$ with $y$, $y$ with $z$ and $y$ with $y$. The path
$\{u\}$ of $E$ joins $u$ with~$u$.

One can verify the following.

\propo{Remark 2.1}{A nonempty hypergraph $H$ is connected}

{\it iff $\Omega(H)$ is connected;}

{\it iff for every $x,y\in\bigcup H$ there is a path of $H$ that
joins $x$ with~$y$.}

\vspace{2ex}

\noindent This shows that our notion of connected hypergraph is
the same as the notion in \cite{B89} (Section 1.2). We have in
this remark the assumption that $H$ is nonempty because otherwise
the graph $\Omega(H)$ would be without vertices (and edges), and
this presumably goes counter to common usage in graph theory (see
\cite{H69}, Chapter~2, and \cite{HR74}; cf.\ Appendix A).
Otherwise, if $\Omega(H)$ is allowed to be without vertices, then
we may lift the assumption of nonemptiness for~$H$.

For every nonempty hypergraph $H$ there is a unique hypergraph
partition $\{H_1,\ldots,H_n\}$, with $n\geq 1$, of $H$ such that
for every $i\in\{1,\ldots,n\}$ we have that $H_i$ is a connected
hypergraph on $\bigcup H_i$. We call this hypergraph partition the
\emph{finest} hypergraph partition of $H$. For example, $E'$ above
is the finest hypergraph partition of~$E$.

The trivial partition $\{H\}$ is the \emph{coarsest} hypergraph
partition of a nonempty hypergraph $H$. If $H$ is nonempty and
connected, then the finest and coarsest hypergraph partitions of
$H$ coincide. (The empty hypergraph, which also happens to be
connected by our definition, has only one partition, namely
$\emptyset$, which may be taken as the finest and coarsest
hypergraph partition of this hypergraph.)

\section{\large \textbf{Constructions}}
In this section we introduce three equivalent notions that may
intuitively be understood as constructions of hypergraphs. They
are called construction, f-construction and s-construction. The
first two notions are based on sets, while the last is based on
words, i.e.\ finite sequences. Of the first two notions, the
notion of f-construction (``f'' comes from ``forest'') is perhaps
more intuitive---it involves a more direct record of constructing.
But the equivalent notion that we call simply construction is the
notion on which we rely in the remainder of the text, and to which
for this reason we give prominence. (Since our constructions are
hypergraphs, we could call them h-constructions, but it would be
onerous to write ``h'' all the time.) The third equivalent notion,
the notion of s-construction (``s'' may be associated with
``syntax''), is based on a notion investigated in~\cite{DP10a}. It
provides the most economical notation.

For $F\subseteq \P C$ and $Y\subseteq C$ let
\[F_Y=_{df}\{X\in F\mid X\subseteq Y\}.\]
We are interested in this definition in particular when $F$ is a
hypergraph $H$ and $Y\subseteq\bigcup H$.

We call a hypergraph $H$ \emph{atomic} when for every $x$ in
$\bigcup H$ we have that $\{x\}\in H$ (cf.\ Lemma 3.9 of
\cite{FS05} and Definition 7.1 of \cite{P09}). Note that the empty
hypergraph is atomic, for trivial reasons. One can verify the
following.

\prop{Remark 3.1}{The hypergraph $H$ is atomic iff for every
subset $Y$ of $\bigcup H$ we have that $H_Y$ is a hypergraph
on~$Y$.}

\noindent With the help of this remark we establish easily the
following.

\prop{Remark 3.2}{If $H$ is an atomic hypergraph and
$Y\subseteq\bigcup H$, then $H_Y$ is an atomic hypergraph on~$Y$.}

For $H$ an atomic hypergraph, we define families of subsets of
$\bigcup H$ that we call \emph{constructions} of $H$. This
definition is by induction on the cardinality $|\bigcup H|$
of~$\bigcup H$:
\begin{itemize}
\item[(0)] if $|\bigcup H|=0$, then $H$ is the empty hypergraph
$\emptyset$, and $\emptyset$ is the only construction of
$\emptyset$;

\item[(1)] if $|\bigcup H|\geq 1$, and $H$ is connected, and $K$
is a construction of $H_{\ccup H-\{x\}}$ for $x\in\bigcup H$, then
$K\cup\{\bigcup H\}$ is a construction of~$H$;

\item[(2)] if $|\bigcup H|\geq 2$, and $H$ is not connected, and
$\{H_1,\ldots,H_n\}$, where $n\geq 2$, is the finest hypergraph
partition of $H$, and for every $i\in\{1,\ldots,n\}$ we have that
$K_i$ is a construction of $H_i$, then $K_1\cup\ldots\cup K_n$ is
a construction of~$H$.
\end{itemize}
This concludes our inductive definition of a construction of $H$.
Note that $\bigcup H$ in clause (1) is $\bigcup K\cup\{x\}$, and
$x\not\in\bigcup K$.

For this definition to be correct, in clause (1) we must verify
that $H_{\ccup H-\{x\}}$ is an atomic hypergraph for every $x$ in
$\bigcup H$, and in clause (2) we must verify that $H_i$ is an
atomic hypergraph for every $i\in\{1,\ldots,n\}$. For both of
these verifications we use Remark 3.2. The present remark about
correctness of definition applies also to the inductive
definitions of f-constructions and w-constructions, to be given
later in this section.

It is easy to verify that a construction of an atomic hypergraph
$H$ on $\bigcup H$ is itself a hypergraph on $\bigcup H$ (though
not necessarily atomic). In particular cases, it will be a
subfamily of $H$ (see Section~4).

To give examples of constructions, consider the atomic hypergraph
\[A=\{\{x\},\{y\},\{z\},\{u\},\{x,y\},\{y,z\},\{z,u\},\{x,y,z\}\}.\]
We could draw this hypergraph in the following manner:
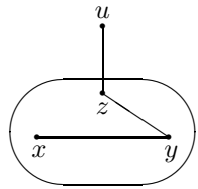
\begin{center}
\begin{picture}(50,65)(0,-6)
\put(0,10){\line(1,0){50}} \put(50,10){\line(-3,2){25}}
\put(25,26.5){\line(0,1){25}}

\put(1,10){\makebox(0,0){\circle*{2}}}
\put(51,10){\makebox(0,0){\circle*{2}}}
\put(26,27){\makebox(0,0){\circle*{2}}}
\put(26,52){\makebox(0,0){\circle*{2}}}

\put(1,7){\makebox(0,0)[t]{$x$}} \put(51,7){\makebox(0,0)[t]{$y$}}
\put(25,23){\makebox(0,0)[t]{$z$}}
\put(25,55){\makebox(0,0)[b]{$u$}}

\put(25,12){\oval(70,40)}

\end{picture}
\end{center}
In such drawings, the circles corresponding to singletons are
taken for granted, and instead of the circles corresponding to
two-element sets, we draw edges between the two elements of such
sets, as for graphs.

A construction of $A$ is the hypergraph
\[L=\{\{u\},\{z,u\},\{y,z,u\},\{x,y,z,u\}\}\]
on $\{x,y,z,u\}$. Here is how $L$ was obtained by our definition.
We had first $L_0=\emptyset$ as a construction of the atomic
hypergraph $\emptyset$. Then $L_1=L_0\cup\{\{u\}\}=\{\{u\}\}$ was
a construction of the atomic hypergraph $\{\{u\}\}$. Next we had
$L_2=L_1\cup \{\{z,u\}\}=\{\{u\},\{z,u\}\}$ as a construction of
the atomic hypergraph $\{\{z\},\{u\},\{z,u\}\}$. Then we had
$L_3=L_2\cup\{\{y,z,u\}\}=\{\{u\},\{z,u\},\{y,z,u\}\}$ as a
construction of the atomic hypergraph
\[A'=\{\{y\},\{z\},\{u\},\{y,z\},\{z,u\}\}.\]
Finally, we have $L=L_3\cup\{\{x,y,z,u\}\}$ as a construction of
$A$. In all that, we applied only clauses (0) and (1) of our
definition of a construction.

Here is another construction of $A$:
\[M=\{\{y\},\{u\},\{y,z,u\},\{x,y,z,u\}\}.\]
It was obtained from the construction $M_1=\{\{u\}\}$ of
$\{\{u\}\}$ and the construction $M_2=\{\{y\}\}$ of $\{\{y\}\}$ by
applying clause (2) of the definition, which yields $M_1\cup
M_2=\{\{y\},\{u\}\}$ as a construction of the atomic hypergraph
$\{\{y\},\{u\}\}$, which is not connected. Then we had by clause
(1) that $M_3=M_1\cup M_2\cup\{\{y,z,u\}\}$ is a construction of
$A'$, mentioned above, and, finally, $M=M_3\cup\{\{x,y,z,u\}\}$ is
a construction of~$A$.

Note that $L$ and $M$ would be constructions also of the atomic
hypergraph $A''$ obtained from $A$ by rejecting $\{x,y,z\}$ (we
deal with that matter in Section~4). They would also be
constructions of atomic hypergraphs more different from $A$ than
$A''$ . Such are, for example, the hypergraphs
\begin{tabbing}
\hspace{1.5em}\=(0)\hspace{1em}\=$s(\emptyset)=e$,\kill

\>\>$A^\circ=A''\cup \{\{u,x\}\}$,\\
\>\>$A^*=(A''\mn\{\{x,y\}\})\cup\{\{x,z\}\}$,
\end{tabbing}
which, together with $A''$, may be drawn as follows:
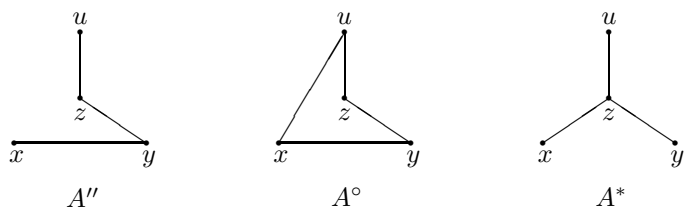
\begin{center}
\begin{picture}(50,80)(100,-15)
\put(0,10){\line(1,0){50}}

\put(50,10){\line(-3,2){25}}

\put(25,26.5){\line(0,1){25}}

\put(1,10){\makebox(0,0){\circle*{2}}}
\put(51,10){\makebox(0,0){\circle*{2}}}
\put(26,27){\makebox(0,0){\circle*{2}}}
\put(26,52){\makebox(0,0){\circle*{2}}}

\put(1,7){\makebox(0,0)[t]{$x$}}

\put(51,7){\makebox(0,0)[t]{$y$}}

\put(25,23){\makebox(0,0)[t]{$z$}}

\put(25,55){\makebox(0,0)[b]{$u$}}


\put(100,10){\line(1,0){50}} \put(100,10){\line(3,5){25}}
\put(150,10){\line(-3,2){25}} \put(125,26.5){\line(0,1){25}}

\put(101,10){\makebox(0,0){\circle*{2}}}
\put(151,10){\makebox(0,0){\circle*{2}}}
\put(126,27){\makebox(0,0){\circle*{2}}}
\put(126,52){\makebox(0,0){\circle*{2}}}

\put(101,7){\makebox(0,0)[t]{$x$}}
\put(151,7){\makebox(0,0)[t]{$y$}}
\put(125,23){\makebox(0,0)[t]{$z$}}
\put(125,55){\makebox(0,0)[b]{$u$}}


\put(200,10){\line(3,2){25}} \put(250,10){\line(-3,2){25}}
\put(225,26.5){\line(0,1){25}}

\put(201,10){\makebox(0,0){\circle*{2}}}
\put(251,10){\makebox(0,0){\circle*{2}}}
\put(226,27){\makebox(0,0){\circle*{2}}}
\put(226,52){\makebox(0,0){\circle*{2}}}

\put(201,7){\makebox(0,0)[t]{$x$}}
\put(251,7){\makebox(0,0)[t]{$y$}}
\put(225,23){\makebox(0,0)[t]{$z$}}
\put(225,55){\makebox(0,0)[b]{$u$}}


\put(26,-7){\makebox(0,0)[t]{$A''$}}
\put(126,-7){\makebox(0,0)[t]{$A^\circ$}}
\put(226,-7){\makebox(0,0)[t]{$A^*$}}
\end{picture}
\end{center}
By removing $x$ from these three hypergraphs we obtain the
hypergraph $A'$. (These three hypergraphs should be compared with
Examples 5.15, 5.13 and 5.16 of \cite{DP10a}, and with
$H'_{4321}$, $H^\circ_{4441}$ and $H^*_{4331}$ in Appendix~B.)

We will make a comment on the intuitive meaning of our
constructions after introducing the equivalent notion of
f-construction, and after giving analogous examples of
f-constructions.

The definition of an \emph{f-construction} of an atomic hypergraph
$H$ is again by induction on $|\bigcup H|$, and its clauses (0)
and (2) are exactly as in the definition of a construction above,
with ``construction'' replaced by ``f-construction''. For clause
(1) of the new definition we make that replacement, and moreover
$K\cup\{\bigcup H\}$ is replaced by $\{K\cup\{x\}\}$. This
concludes the definition of an f-construction.

To give examples of f-constructions, consider again the atomic
hypergraph $A$ above. An f-construction of $A$ is
$L^f=\{\{x,\{y,\{z,\{u\}\}\}\}\}$. (Note that $L^f$ is a
singleton.) We will show exactly later how this f-construction
corresponds to the construction $L$ above. Here is how $L^f$ was
obtained by our definition. We had first $L^f_0=\emptyset=L_0$ as
a construction of $\emptyset$. Then
$L^f_1=\{L^f_0\cup\{u\}\}=\{\{u\}\}=L_1$ was an f-construction of
$\{\{u\}\}$. Next we had
$L^f_2=\{L^f_1\cup\{z\}\}=\{\{z,\{u\}\}\}$ as a construction of
$\{\{z\},\{u\},\{z,u\}\}$. Then we had
$L^f_3=\{L^f_2\cup\{y\}\}=\{\{y,\{z\{u\}\}\}\}$ as a construction
of $A'$ above. Finally, we have $L^f=\{L_3\cup\{x\}\}$ as a
construction of $A$. In this example, we applied only the new
clauses (0) and (1).

Another f-construction of $A$ is
$M^f=\{\{x,\{z,\{y\},\{u\}\}\}\}$, which corresponds to $M$. In
obtaining $M^f$ by the definition of an f-construction we apply
also clause~(2).

These examples should explain the denomination ``construction'' in
our constructions and f-constructions. The hypergraphs $L$ and
$M$, as well as the sets $L^f$ and $M^f$, may be understood as
constructions of $A$ in time. Within a connected part of $A$ the
construction proceeds by adding in $L^f$ and $M^f$ a chosen
vertex, and this choice induces a temporal order. Clause (1)
serves for that. Connected parts of $A$ that are mutually
disconnected are added simultaneously, without order, and clause
(2) serves for that.

We define next another notion equivalent to the notion of
construction. For an atomic hypergraph $H$ we define first words
in the alphabet $\bigcup H$, which we call \emph{w-constructions}.
This definition is again by induction on $\bigcup H$:
\begin{itemize}
\item[(0)] if $|\bigcup H|=0$, then the empty word $e$ is the only
w-construction of the hypergraph $\emptyset$;

\item[(1)] if $|\bigcup H|\geq 1$, and $H$ is connected, and $t$
is a w-construction of $H_{\ccup H-\{x\}}$ for $x\in\bigcup H$,
then $xt$ is a w-construction of~$H$;

\item[(2)] if $|\bigcup H|\geq 2$, and $H$ is not connected, and
$\{H_1,\ldots,H_n\}$, where $n\geq 2$, is the finest hypergraph
partition of $H$, and for every $i\in\{1,\ldots,n\}$ we have that
$t_i$ is a w-construction of $H_i$, then $(t_1\pl\ldots\pl t_n)$
is a w-construction of~$H$.
\end{itemize}
This concludes our inductive definition of a w-construction of
$H$.

Consider equivalence classes of w-constructions of $H$ obtained by
factoring through the commutativity of +. We call these
equivalence classes \emph{s-constructions}. We refer to an
s-construction by any w-construction that belongs to it. Our
s-constructions are analogous to the \textbf{S}-forests of
\cite{DP10a} (Section~5).

Here are two examples of s-constructions of the hypergraph $A$
above. These are $L^s$, which is $xyzu$, and $M^s$, which is
$xz(y\pl u)$; they correspond to $L$ and $M$ respectively. The
s-construction $xz(u\pl y)$ is the same as~$M^s$.

Our task now is to show that the notions of construction,
f-construction and s-construction are all equivalent. By this we
mean that there are structure-preserving bijections between the
sets $\C(H)$, $\C^f(H)$ and $\C^s(H)$ of, respectively,
constructions, f-constructions and s-constructions of an atomic
hypergraph~$H$.

We define first a map $f\!:\C(H)\rightarrow\C^f(H)$ by induction
on $\bigcup H$, in parallel with the clauses of the inductive
definitions of construction and f-construction:
\begin{tabbing}
\hspace{1.5em}\=(0)\hspace{1em}\=$f(\emptyset)=\emptyset$,\\
\>(1)\>$f(K\cup\{\bigcup H\})=\{f(K)\cup\{x\}\}$,\\
\>(2)\>$f(K_1\cup\ldots\cup K_n)=f(K_1)\cup\ldots\cup f(K_n)$.
\end{tabbing}
The conditions concerning $H$, $K$, $K_1,\ldots,K_n$ are taken
from the clauses (1) and (2) of the definition of a construction.
For (1) we have that $H$ is connected and $K$ is a construction of
$H_{\ccup H-\{x\}}$ for $x\in\bigcup H$. For (2) we have that $H$
is not connected and $K_1,\ldots,K_n$, with $n\geq 2$, are
constructions of respectively $H_1,\ldots,H_n$ for
$\{H_1,\ldots,H_n\}$ being the finest hypergraph partition of $H$.
We proceed analogously for the two other maps below.

Next we define analogously a map $s\!:\C^f(H)\rightarrow\C^s(H)$:
\begin{tabbing}
\hspace{1.5em}\=(0)\hspace{1em}\=$s(\emptyset)=e$,\\
\>(1)\>$s(\{K\cup\{x\}\})=x\:s(K)$,\\
\>(2)\>$s(K_1\cup\ldots\cup K_n)=s(K_1)\pl\ldots\pl s(K_n)$.
\end{tabbing}
Finally, we define analogously a map
$c\!:\C^s(H)\rightarrow\C(H)$:
\begin{tabbing}
\hspace{1.5em}\=(0)\hspace{1em}\=$c(e)=\emptyset$,\\
\>(1)\>$c(xt)=c(t)\cup\{\bigcup H\}$,\\
\>(2)\>$c(t_1\pl\ldots\pl t_n)=c(t_1)\cup\ldots\cup c(t_n)$.
\end{tabbing}

Then to verify that $f$, $s$ and $c$ are bijections it is enough
to verify the following three equations:
\[c(s(f(K)))=K,\hspace{3em}f(c(s(K)))=K,\hspace{3em}s(f(c(t)))=t,\]
which is quite straightforward.

Constructions, f-constructions and s-constructions bear a forest
structure (a forest is a disjoint union of trees, with a tree
being a limit case). This structure is clearer in f-constructions
and s-constructions. In f-constructions, the $x$ added in clause
(1) is the \emph{root} of the tree $\{K\cup\{x\}\}$. Here $K$, if
it is nonempty, is equal to a forest $K_1\cup\ldots\cup K_n$, with
$n\geq 1$, and the roots of the trees $K_1,\ldots,K_n$ are the
immediate successors into which $x$ branches.

\section{\large \textbf{Saturation and cognate hypergraphs}}
Two different atomic hypergraphs may have the same constructions.
Such are, for example, $A$ and $A''$ of Section~3. In this section
we concentrate on atomic hypergraphs that have the same
constructions in order to find among them a representative that is
easiest to work with.

In the set of all atomic hypergraphs on the same carrier that have
the same constructions there is a greatest one, which has a
property we will call saturation. We will characterize the
equivalence relation that the hypergraphs in this set bear to each
other in terms of a relation where their difference is reduced to
atomic differences. Such an atomic difference consists in one
hypergraph having a member that is not in the other hypergraph,
but this one member---a dispensable member---does not, roughly
speaking, increase connectedness. Two hypergraphs are called
cognate when they differ only with respect to dispensable members,
and a hypergraph is saturated when all possible dispensable
members are in it.

It will simplify the exposition later if we concentrate on
saturated hypergraphs (see Section~6). The results of the present
section justify this simplification; they show that it makes no
difference with respect to constructions.

We call a hypergraph $H$ \emph{saturated} when for every
$X_1,X_2\in H$ if $X_1\cap X_2\neq\emptyset$, then $X_1\cup X_2\in
H$ (cf.\ Lemma 3.9 of \cite{FS05} and Definition 7.1 of
\cite{P09}). One can verify the following.

\propo{Remark 4.1}{The hypergraph $H$ is saturated}

{\it iff for every $Y\subseteq\bigcup H$, if $H_Y\mn\{Y\}$ is a
connected hypergraph on $Y$, then $Y\in H$;}

{\it iff for every $Y\subseteq\bigcup H$, if $H_Y$ is a connected
hypergraph on $Y$, then $Y\in H$;}

{\it iff for every $Y\subseteq\bigcup H$ we have that $H_Y$ is a
connected hypergraph on $Y$ iff $Y\in H$.}

\vspace{2ex}

The hypergraph $E$ of Section~2 is saturated, and
$E\mn\{\{x,y,z\}\}$, which is $\{\{x,y\},\{y,z\},\{u\},\{v\}\}$,
is not saturated. The hypergraph $A$ of Section~3 is not
saturated. The empty hypergraph is saturated, for trivial reasons.

For a hypergraph $H$, we say that a subset $Y$ of $\bigcup H$ is
\emph{dispensable} in $H$ when $H_Y\mn\{Y\}$ is a connected
hypergraph on $Y$. For example, $\{x,y,z\}$ is dispensable in the
hypergraph $E$. It is also dispensable in $E\mn\{\{x,y,z\}\}$.
Note that singleton members of a hypergraph are never dispensable.

By Remark 4.1, we have that the hypergraph $H$ is saturated iff
every subset of its carrier dispensable in $H$ is an element of
$H$. In terms of dispensability we can also formulate a notion
dual to saturation. We will say that the hypergraph $H$ is
\emph{bare} when no subset of its carrier dispensable in $H$ is an
element of~$H$.

We can prove the following.

\prop{Proposition 4.2}{Suppose $Y$ is dispensable in $H$. Then $Z$
is dispensable in $H$ iff $Z$ is dispensable in $H\cup\{Y\}$.}

\dkz The equivalence of the proposition from left to right is
trivial. For the other direction, suppose $Z$ is dispensable in
$H\cup\{Y\}$.

So $(H\cup\{Y\})_Z\mn\{Z\}$ is a connected hypergraph on $Z$. If
$Y\not\subseteq Z$, then $(H\cup\{Y\})_Z=H_Z$, and we are done.

Suppose $Y\subseteq Z$. Then, by Remark 2.1, for every $x,y\in Z$
there is a path $X_1,\ldots,X_n$ of $(H\cup\{Y\})_Z\mn\{Z\}$ that
joins $x$ with $y$; so $x\in X_1$ and $y\in X_n$. If for every
$i\in\{1,\ldots,n\}$ we have that $X_i\neq Y$, then
$X_1,\ldots,X_n$ is a path of $H_Z\mn Z$, and we are done.

Suppose for some $i$ we have that $X_i=Y$. Then we have the
following four cases: (1) $1<i<n$, (2) $1=i<n$, (3) $1<i=n$ and
(4) $1=i=n$.

In case (1) our path is of the form
\[X_1,\ldots,X_{i-1},Y,X_{i+1},\ldots,X_n\]
with $x'\in X_{i-1}\cap Y$ and $y'\in Y\cap X_{i+1}$, where
$X_1,\ldots,X_{i-1}$ is a path of $H_Z\mn\{Z\}$ that joins $x$
with $x'$ and $X_{i+1},\ldots,X_n$ is a path of $H_Z\mn\{Z\}$ that
joins $y'$ with $y$. Since $Y$ is dispensable in $H$, we have that
$H_Y\mn\{Y\}$ is a connected hypergraph on $Y$. By Remark 2.1,
this means that there is a path $Y_1,\ldots,Y_m$ of $H_Y\mn\{Y\}$
that joins $x'$ with $y'$. For every $j\in\{1,\ldots,m\}$ we have
that $Y_j\subset Y$. Since $Y\subseteq Z$, we obtain that
$Y_j\subset Z$, and hence $Y_1,\ldots,Y_m$ is a path of
$H_Z\mn\{Z\}$. So either
\[X_1,\ldots,X_{i-1},Y_1,\ldots,Y_m,X_{i+1},\ldots,X_n\]
is a path of $H_Z\mn\{Z\}$ that joins $x$ with~$y$, or it can
easily be transformed into such a path by contracting, if for some
$l\in\{1,\ldots,n\}\mn\{i\}$ and some $j\in\{1,\ldots,m\}$ we have
that $X_l$ is $Y_j$.

We take that $x'$ is $x$ in cases (2) and (4), and we take that
$y'$ is $y$ in cases (3) and (4). In all these three remaining
cases we proceed analogously to case (1). This is enough to
conclude that $Z$ is dispensable in $H$.\qed

We will say that a hypergraph $H\cup\{Y\}$ \emph{enhances} the
hypergraph $H$ when $Y$ is dispensable in $H$ and $Y\not\in H$.

As a corollary of Proposition 4.2 we have that $Y$ is dispensable
in $H$ iff $Y$ is dispensable in $H\cup\{Y\}$. (It is easy to
prove this corollary directly.) This corollary shows that in the
definition of enhancement the dispensability of $Y$ in $H$ amounts
to the dispensability of $Y$ in $H\cup\{Y\}$, and the later
dispensability could serve for the definition.

Consider the equivalence relation on the set of hypergraphs on the
same carrier $C$ obtained as the reflexive, symmetric and
transitive closure of the relation of enhancement. When two
hypergraphs on $C$ are in this relation we say that they are
\emph{cognate}. As a corollary of Proposition 4.2 we obtain the
following.

\prop{Proposition 4.3}{Suppose $H$ and $J$ are cognate hypergraphs
on $C$. Then for every $Z\subseteq C$ we have that $Z$ is
dispensable in $H$ iff $Z$ is dispensable in~$J$.}

A \emph{cognate set} of hypergraphs is an equivalence class of
hypergraphs with respect to the cognation equivalence relation.
With the help of Proposition 4.3, we establish that a cognate set
is a lattice with respect to intersection and union. It has a
greatest element, the union of all its members, which is a
saturated hypergraph, and it has a least element, the intersection
of all its members, which is a bare hypergraph.

We need the following remark.

\prop{Remark 4.4}{For $Y\subseteq Z\subseteq\bigcup H$ we have
that $Y$ is dispensable in $H$ iff $Y$ is dispensable in $H_Z$.}

\noindent This is because $(H_Z)_Y=H_Y$. We also need the
following lemma.

\prop{Lemma 4.5}{Suppose $Y$ is dispensable in $H$. Then $H$ is
connected iff $H\cup\{Y\}$ is connected.}

\dkz This proof will be quite similar to the proof of Proposition
4.2. The equivalence of the lemma from left to right is trivial.
For the other direction suppose there is a path $X_1,\ldots,X_n$
of $H\cup\{Y\}$ that joins $x$ with $y$, for $x,y\in\bigcup H$. If
for every $i\in\{1,\ldots,n\}$ we have that $X_i\neq Y$, then
$X_1,\ldots,X_n$ is a path of $H$, and we are done.

Suppose for some $i$ we have that $X_i=Y$. Then we have the
following four cases: (1) $1<i<n$, (2) $1=i<n$, (3) $1<i=n$ and
(4) $1=i=n$.

In case (1) we proceed as in the proof of Proposition 4.2 until we
reach the path $Y_1,\ldots,Y_m$ of $H_Y\mn\{Y\}$ that should
replace $Y$. This path is made of elements of $H$, and after the
replacement we have a path of $H$ that joins $x$ with $y$. In the
remaining cases we proceed analogously.\qed

Then we can prove the following.

\prop{Proposition 4.6}{Suppose $H$ and $J$ are cognate hypergraphs
on $C$. Then for every $Z\subseteq C$ we have that $H_Z$ is
connected iff $J_Z$ is connected.}

\dkz It is enough to prove this proposition when $J$ is
$H\cup\{Y\}$ for $Y$ dispensable in $H$. It is clear that if $H_Z$
is connected, then $(H\cup\{Y\})_Z$ is connected.

Suppose $(H\cup\{Y\})_Z$ is connected. We have that
$(H\cup\{Y\})_Z$ is different from $H_Z$ only when $Y\subseteq Z$.
Then by Remark 4.4 we have that $Y$ is dispensable in $H_Z$. Since
$Y\subseteq Z$, we also have that $(H\cup\{Y\})_Z=H_Z\cup\{Y\}$.
It suffices to apply Lemma 4.5 to obtain that $H_Z$ is
connected.\qed

Remember that for an atomic hypergraph $H$ the set $\C(H)$ is the
set of all constructions of $H$. We will prove the following.

\prop{Proposition 4.7}{If $H$ and $J$ are cognate atomic
hypergraphs, then $\C(H)=\C(J)$.}

\dkz It is enough to prove this proposition when $J$ is
$H\cup\{Y\}$ for $Y$ dispensable in $H$. We establish first that
every construction of $H$ is a construction of $H\cup\{Y\}$.

We proceed by induction on the cardinality of $\bigcup H$, as in
the inductive definition of a construction. If $\bigcup
H=\emptyset$, then $H=H\cup\{Y\}=\emptyset$.

Suppose $|\bigcup H|\geq 1$ and $H$ is connected. By Lemma 4.5, we
have that $H\cup\{Y\}$ is connected. Then a construction of $H$ is
of the form $K\cup\{\bigcup H\}$ for $K$ a construction of
$H_{\ccup H-\{x\}}$, where $x\in\bigcup H$. If $x\in Y$, then
$(H\cup\{Y\})_{\ccup H-\{x\}}=H_{\ccup H-\{x\}}$. So $K$ is a
construction of $(H\cup\{Y\})_{\ccup H-\{x\}}$. If $x\not\in Y$,
then we have that $Y\subseteq\bigcup H\mn\{x\}\subseteq\bigcup H$,
and by Remark 4.4 and the induction hypothesis, $K$ is a
construction of $(H\cup\{Y\})_{\ccup H-\{x\}}$. Hence
$K\cup\{\bigcup H\}$ is a construction of $H\cup\{Y\}$.

Suppose $|\bigcup H|\geq 2$ and $H$ is not connected. By Lemma
4.5, we have that $H\cup\{Y\}$ is not connected. Suppose
$\{H_1,\ldots,H_n\}$, where $n\geq 2$, is the finest hypergraph
partition of $H$. Then a construction of $H$ is of the form
$K_1\cup\ldots\cup K_n$ for $K_1,\ldots,K_n$ constructions of
$H_1,\ldots,H_n$ respectively. For some $i\in\{1,\ldots,n\}$ we
must have that $Y\subseteq\bigcup H_i$. Hence $Y$ is dispensable
in $H_i$ by Remark 4.4 (we have that $Y\subseteq\bigcup
H_i=Z\subseteq\bigcup H$). Then, by the induction hypothesis,
$K_i$ is a construction of $H_i\cup\{Y\}$, and $K_1\cup\ldots\cup
K_n$ is a construction of $H\cup\{Y\}$.

We proceed analogously to establish in the converse direction that
every construction of $H\cup\{Y\}$ is a construction of~$H$. \qed

The \emph{saturated closure} $\bar{H}$ of a hypergraph $H$ is the
saturated hypergraph in the cognate set of hypergraphs to which
$H$ belongs. We can prove the following.

\prop{Proposition 4.8}{For an atomic hypergraph $H$ we have that
$\bigcup\C(H)=\bar{H}$.}

\dkz We prove this proposition first for $H=\bar{H}$. From left to
right, suppose that for some construction $K$ of $H$ we have that
$Y\in K$. Then $Y\subseteq\bigcup H$ and $H_Y$ is a connected
hypergraph on $Y$. So, by Remark 4.1, we have that $Y\in H$.

From right to left, suppose that $Y\in H$. We show by induction on
$k=|\bigcup H\mn Y|$ that there is a construction $K$ of $H$ such
that $Y\in K$. If $k=0$, then $Y=\bigcup H$, and $H$ is a
connected hypergraph on $Y$. For an arbitrary construction $K$ of
$H$ we must have that $Y\in K$. If $k>0$ and $x\in\bigcup H\mn Y$,
then by the induction hypothesis we have a construction $K'$ of
$H_{\ccup H-\{x\}}$ such that $Y\in K'$. For $\{H_1,\ldots,H_n\}$
being the finest hypergraph partition of $H$, and $x\in\bigcup
H_i$, we define the construction $K$ of $H$ by $K=K'\cup\{\bigcup
H_i\}$, and we have that $Y\in K$.

It remains to remark that $\C(H)=\C(\bar{H})$ for an arbitrary
atomic hypergraph $H$, which we have by Proposition~4.6.\qed

The following proposition, which completes Proposition 4.7,
characterizes cognate classes in terms of constructions. Atomic
hypergraphs are cognate iff they have the same constructions.

\prop{Proposition 4.9}{Suppose $H$ and $J$ are atomic hypergraphs.
Then $H$ and $J$ are cognate iff $\C(H)=\C(J)$.}

\dkz From left to right we have Proposition 4.7. For the other
direction, suppose $\C(H)=\C(J)$. Then, with the help of
Proposition 4.8, we obtain that $\bar{H}=\bar{J}$. So $H$ and $J$
are in the same cognate set.\qed

\vspace{-2ex}

\section{\large \textbf{Constructs and abstract polytopes of
hypergraphs}}
In this section we define the abstract polytopes of atomic
hypergraphs with the help of the notion of construct, which is a
notion derived from our notion of construction. The notion of
construction is the notion on which all the burden rests. The
proof that the polytopes so defined are indeed abstract polytopes
will be given in Section~8.

We start with the following remark.

\prop{Remark 5.1}{For every construction $K$ of an atomic
hypergraph $H$ we have that $|K|=|\bigcup H|$.}

\noindent This is established in a straightforward manner by
induction on the size of~$K$.

Note also that if $\{H_1,\ldots,H_n\}$, with $n\geq 1$, is the
finest hypergraph partition of the atomic hypergraph $H$, then for
every $i\in\{1,\ldots,n\}$ we have that $\bigcup H_i\in K$. We say
that $H_i$ is a \emph{connected component} of $H$, and $\bigcup
H_i$ is a \emph{connected component} of the carrier $\bigcup H$ of
$H$. The number $n$ is the \emph{connectedness number} of $H$; it
is the number of connected components of $H$, or of~$\bigcup H$.
The atomic hypergraph $\emptyset$ and its carrier $\emptyset$ have
just one partition $\emptyset$, with 0 connected components; so
the connectedness number of this hypergraph is 0.

A \emph{construct} of an atomic hypergraph $H$ is a subfamily (not
necessarily proper) of a construction of $H$ that contains every
connected component of the carrier $\bigcup H$ of $H$. For
example,
\[\{\{u\},\{y,z,u\},\{x,y,z,u\}\}\]
is a construct of the hypergraph $A$ of Section~3. It is a
subfamily of both of the constructions $L$ and $M$ of~$A$.

It is clear that every construct of $H$ is a hypergraph on
$\bigcup H$, as $H$ is. It is also clear, in accordance with
Remark 5.1, that for $C$ a construct of $H$ we have that
$|C|\leq|\bigcup H|$.

The constructs of the atomic hypergraph $H$ serve to define as
follows the \emph{abstract polytope} of $H$, which we designate by
$\A(H)$ (for the definition of an abstract polytope in general,
and related notions used below, see \cite{MS02}, Section 2A, and
our Section~8 below).

The elements of $\A(H)$, i.e.\ the faces of $\A(H)$, are all the
constructs of $H$ plus the set $\bar{H}^*$, which is
$\bar{H}\cup\{*\}$, for $\bar{H}$ being the saturated closure of
$H$ (see Section~4) and $*$ a new element that is not in~$\bigcup
H$.

We take $\A(H)$ as a partial order with the inverse of the subset
relation; i.e.\ for the faces $C_1$ and $C_2$ of $\A(H)$ we have
that $C_1\leq C_2$ when $C_2\subseteq C_1$. So the incidence
relation of $\A(H)$ is the symmetric closure of the subset
relation. In the partial order $\A(H)$ the element $\bar{H}^*$ is
the least element.

Proposition~4.8 states that the union of all the constructions of
$H$ is $\bar{H}$. Hence the union of all the constructs of $H$ is
$\bar{H}$. So it seems we could take simply $\bar{H}$ instead of
$\bar{H}^*$ as the least element. We did not do that for the
following reason.

If all the elements of $H$ are singletons (which means that $H$
can also be empty), then $\bar{H}$ coincides with $H$, which is
the only construct of $H$, and the only construction of $H$. We
want however to distinguish even in that case the construct
$\bar{H}$ from $\bar{H}^*$.

The choice of $\bar{H}^*$ is also dictated by our wish to base the
incidence relation in $\A(H)$ on the subset relation at every
level. We could however obtain the same effect by having instead
of $\bar{H}^*$ any set in which $\bar{H}^*$ is included; for
example, the power set of $\bigcup H$, with $\emptyset$ being $*$,
or even a universal set in which all the sets $\bar{H}^*$ are
included. We could also replace $\bar{H}^*$ by anything different
from the other elements of $\A(H)$ if we do not insist that
incidence with respect to it must be based on the subset relation.

If $n$ is the connectedness number of $H$, then the rank $r$ of
$\A(H)$ is $|\bigcup H|\mn n$. In general, we have that $r\geq 0$.
The rank is the dimension of the realization of the abstract
polytope in space (see Section~9). In our example with the
hypergraph $A$ of Section~2, we have that $|\bigcup A|=4$ and
$n=1$; so the rank of $\A(A)$ is $3$. The abstract polytope
$\A(A)$ corresponds to the three-dimensional associahedron $K_5$
(see ($H'_{4321}$) in Appendix~B, and references therein).

The least face $F_{-1}$ of $\A(H)$ is $\bar{H}^*$, and the
greatest face $F_r$ is the set $\{\bigcup H_1,\ldots,\bigcup
H_n\}$ of the connected components of the carrier $\bigcup H$ of
$H$. (If all the elements of $H$ are singletons, then $F_r\cup
\{*\}=F_{-1}$.) With the hypergraph $A$, we have that $F_{-1}$ of
$\A(A)$ is $A\cup \{\{y,z,u\},\{x,y,z,u\},*\}$, while $F_3$ of
$\A(A)$ is $\{\{x,y,z,u\}\}$.

The vertices, i.e.\ the faces of rank 0, of $\A(H)$ are the
constructions of $H$. By Remark 5.1, the cardinality of every
vertex is $|\bigcup H|$. With the hypergraph $A$, we have 14
vertices in $\A(A)$, among which we find $L$ and $M$ of Section~3.

Besides $L$, there are seven more vertices of the same type. These
eight vertices correspond to the s-constructions
\[
\begin{array}{llll}
xyzu, & xyuz, & xuyz, & xuzy,\\
uzyx, & uzxy, & uxzy, & uxyz,
\end{array}
\]
the first s-construction $xyzu$ corresponding exactly to $L$. How
these eight vertices are distributed in $K_5$ may be seen in the
second picture of $K_5$ in Appendix~B (see ($H'_{4321}$)), which
is based on a picture in \cite{DP10a} (Example 5.15, where $xyzu$
is written $x\!\cdot\! y\!\cdot\! z\!\cdot\! u$; we omit $\cdot$
now). There are two vertices of type $M$, which correspond to the
s-constructions
\[
\begin{array}{ll}
xz(y\pl u), & uy(x\pl z),
\end{array}
\]
with the first s-construction corresponding exactly to $M$.

There are four vertices of another type, corresponding to the
s-constructions
\[
\begin{array}{ll}
y(x\pl(zu)), & y(x\pl(uz)),\\
z((xy)\pl u), & z((yx)\pl u),
\end{array}
\]
where, for example, the first s-construction $y(x\pl(zu))$
corresponds to the construction
\[N=\{\{x\},\{u\},\{z,u\},\{x,y,z,u\}\}.\]
(Here it is clear how much the notation of s-constructions is more
economical.) With that we have obtained all the 14 vertices of
$\A(A)$, whose distribution may be seen in the picture of $K_5$ of
Appendix~B, mentioned above.

The edges, i.e\ the faces of rank 1, of $\A(H)$ are all the
constructs of $H$ of cardinality $|\bigcup H|\mn 1$. For example,
the edge joining $L$ and $M$ in $\A(A)$ is $L\cap
M=\{\{u\},\{y,z,u\},\{x,y,z,u\}\}$, while the edge joining $L$ and
$N$ is $L\cap N=\{\{u\},\{z,u\},\{x,y,z,u\}\}$. We will ascertain
later that for every edge of $\A(H)$ there are exactly two
different vertices such that our edge is their intersection. (This
follows from property (P4) when $i=0$; see Section~8.)

In general, for $k\geq 0$, the faces of rank $k$ of $\A(H)$ are
all the constructs of $H$ of cardinality $|\bigcup H|\mn k$, and
if $k=-1$, then the unique face of rank $-1$ is $\bar{H}^*$, whose
cardinality is $|\bar{H}|\mn k=|\bar{H}|\pl 1$.

So the facets, i.e.\ the faces of rank $r\mn 1$, of $\A(H)$, where
$r$ is the rank of $\A(H)$, are all the constructs of $H$ of
cardinality $|\bigcup H|\mn (r\mn 1)=|\bigcup H|\mn (|\bigcup
H|\mn n\mn 1)=n\pl 1$. Besides the $n$ connected components of
$\bigcup H$ we find in each facet a single additional member. This
member is from $\bar{H}$ if $r>0$, and it is $*$ if $r=0$.

We have that $r=0$ for $\A(H)$ in the following two cases. The
first case is when $H=\{\{x_1\},\ldots,\{x_n\}\}$, for $n\geq 1$.
Then\[\A(H)=\{\{\{x_1\},\ldots,\{x_n\}\},\{\{x_1\},\ldots,\{x_n\},
*\}\},\] with $F_0$ being the vertex $\{\{x_1\},\ldots,\{x_n\}\}$,
and $F_{-1}$ being the facet $\{\{x_1\},\ldots,$ $\{x_n\},*\}$.
The second case is when $H=\bigcup H=\emptyset$. Then
$\A(H)=\{\emptyset,\{*\}\}$, with $F_0$ being the vertex
$\emptyset$, and $F_{-1}$ being the facet $\{*\}$. Both situations
are anomalous for having a vertex strictly above a facet. When
$r=1$, the vertices and facets coincide, and when $r>1$, the
vertices are strictly below the facets.

We have the following.

\prop{Proposition 5.2}{If $r>0$ is the rank of $\A(H)$, then every
vertex of $\A(H)$ is incident with $r$ facets.}

\dkz Every vertex of $\A(H)$, i.e.\ every construction $K$ of $H$,
is of cardinality $|\bigcup H|$, by Remark 5.1. The facets of
$\A(H)$ have each besides the $n$ connected components of $\bigcup
H$ a single additional member. Those facets with which $K$ is
incident have as this additional member a member of $K$ different
from the $n$ connected components of $\bigcup H$. There are
$|\bigcup H|\mn n$ such members in $K$.\qed

If $r=0$, then the additional member of a facet mentioned in this
proof is not from $K$, but it is $*$. If the rank of $\A(H)$ is 0,
then every vertex of $\A(H)$ is incident with 1 facet (in this
case there is a single vertex strictly above the single facet).

In our example with the hypergraph $A$, we have as the facets of
$\A(A)$ the two-element sets that besides $\{x,y,z,u\}$ have as an
additional member one of
\[
\{x\},\{y\},\{z\},\{u\},\{x,y,z\},\{y,z,u\},
\{x,y\},\{y,z\},\{z,u\}.
\]
The first six facets in the ensuing list correspond to pentagons,
while the last three correspond to squares. The facet incident
with the vertices $L$, $M$ and $N$ is $\{\{u\},\{x,y,z,u\}\}$.

If the atomic hypergraph $H$ is not connected, and its connected
components are $H_1,\ldots,H_n$, for $n\geq 2$, then $\A(H)$ may
be obtained out of $\A(H_1),\ldots,\A(H_n)$ in the following
manner. Let $C_1,\ldots,C_n$ be constructs, i.e.\ faces of rank at
least 0, of $\A(H_1),\ldots,\A(H_n)$ respectively. Let the
constructs $C_1,\ldots,C_n$ be respectively of cardinalities
$k_1,\ldots,k_n$. Then ${C_1\cup\ldots\cup C_n}$ of cardinality
${k_1\pl\ldots\pl k_n}$ is a face of $\A(H)$; it is a face of rank
$k$ when $k_1\pl\ldots\pl k_n=|\bigcup H|\mn k$. This is how we
obtain all the faces of $\A(H)$ of rank at least 0. The face
$F_{-1}$ of $\A(H)$ is, as always, $\bar{H}^*$.

We have that $\A(H)\mn\{\bar{H}^*\}$, i.e.\ $\A(H)\mn\{F_{-1}\}$,
is isomorphic to the cartesian product
$(\A(H_1)\mn\{\bar{H}_1^*\})\times\ldots\times(\A(H_n)\mn\{\bar{H}_n^*\})$
(cf.\ the product $\cdot$ of Section~8). Isomorphism means here
the existence of an order-preserving bijection. We may conceive of
$\A(H)$ as obtained from $\A(H_1),\ldots,\A(H_n)$ by an operation
$\otimes$ related to $\times$. Binary $\otimes$ differs from
binary $\times$ by having instead of an ordered pair the union of
the two members of the ordered pair; these members are disjoint,
and their disjoint union corresponds bijectively to the ordered
pair. We conflate moreover the least faces $F_{-1}$ into a single
one. Hence we have $\A(H)=\A(H_1)\otimes\ldots\otimes\A(H_n)$.
(For an example, see $\A(H_{4200})$ in Appendix~B.)

It remains to verify that $\A(H)$ for an arbitrary atomic
hypergraph $H$ is indeed an abstract polytope of rank $|\bigcup
H|\mn n$, in the sense of \cite{MS02} (Section 2A), and this will
be done in Section~8. In the remainder of this section, and in the
next two sections, we consider various properties of abstract
polytopes. This will help us for the results of Section~8, and
also for those of Section~9, where we deal with the realizations
of our abstract polytopes in Euclidean spaces.

The partial order $\A(H)$ is a lattice, with join being
intersection. The meet of two faces of $\A(H)$ is their union if
this union is a construct of $H$, and otherwise it is $\bar{H}^*$,
i.e.\ $F_{-1}$. For example, for $\alpha$ being $\{x,y,z,u\}$ the
union of the facets $\{\{x\},\alpha\}$ and $\{\{y\},\alpha\}$ of
$\A(A)$ is not a construct of $A$; so the meet of these two facets
is $\bar{A}^*$.

In general, the lattice $\A(H)$ is not distributive. For example,
$\A(A)$ is not distributive because
\[
\begin{array}{l}
\{\{y\},\alpha\}\wedge(\{\{x\},\alpha\}\vee\{\{z\},\alpha\})=\{\{y\},\alpha\},\\
(\{\{y\},\alpha\}\wedge\{\{x\},\alpha\})\vee
(\{\{y\},\alpha\}\wedge\{\{z\},\alpha\})=\bar{A}^*
\end{array}
\]
(here the join $\vee$ is intersection, while the meet $\wedge$ is
either union or its result is $\bar{A}^*$, as explained above).

A more natural lattice than $\A(H)$ is the dual lattice (namely,
$\A(H)$ upside down). In the dual lattice the meet would be
intersection, and the join would be the other operation involving
union. As an abstract polytope, we need however $\A(H)$ as it is,
and not the dual lattice, which would give another polytope.

For a face $F$ of $\A(H)$ different from $F_{-1}$, consider the
section $F_r/F$ of $\A(H)$, i.e.\ the set of all the constructs of
$H$ of which the greatest face $F_r$ is a subset, and which are
subsets of the construct $F$. This lattice is isomorphic to the
lattice $\langle\P(F\mn F_r),\cup,\cap\rangle$, with $\P(F\mn
F_r)$ being the power set of $F\mn F_r$, meet being $\cup$, join
being $\cap$; the greatest element of this lattice is $\emptyset$,
and the least element is $F\mn F_r$. Hence, by Proposition 2.16 of
\cite{Z95} (Section 2.5), we may conclude that a geometric
realization of $\A(H)$ whose face lattice is isomorphic to $\A(H)$
must be a simple polytope. (This means that each of its vertex
figures, which are figures obtained by truncating a vertex, is a
simplex; we deal with these matters in Section~9.) Another way to
reach the same conclusion is to rely on Proposition 5.2, and
appeal again to Proposition 2.16 of \cite{Z95}.

We can prove the following.

\prop{Proposition 5.3}{If for some $k\in\{-1,\ldots,r-1\}$ all the
faces of rank $k$ of $\A(H_1)$ and $\A(H_2)$ are the same, then
$\A(H_1)=\A(H_2)$.}

\dkz Let $k\in\{0,\ldots,r\mn 1\}$. Then the union of all the
faces of rank $k$ of $\A(H)$ is $\bar{H}$; this follows from
Proposition 4.8. Hence if the faces of rank $k$ of $\A(H_1)$ and
$\A(H_2)$ are the same, then $\bar{H}_1=\bar{H}_2$, and hence
$\A(H_1)=\A(H_2)$. When $k=-1$, we reason similarly, with
$\bar{H}^*$ instead of $\bar{H}$.\qed

\noindent So, in particular, $\A(H)$ is completely determined by
its vertices, or by its facets.

\section{\large \textbf{Constructions of ASC-hypergraphs}}
Let an \emph{ASC-hypergraph} be a hypergraph that is atomic (see
the beginning of Section~3), saturated (see the beginning of
Section~4) and connected (see Section~2). Concentrating on these
hypergraphs will simplify our exposition. With their help, we will
give in this section noninductive characterizations of
constructions (Proposition 6.11) and of constructs (Proposition
6.13), which presents an alternative to our inductive definition
of a construction in Section~3. These alternative
characterizations will serve for the results of Section~9, which
are about the realizations of polytopes of hypergraphs. They will
also serve in Appendix~A to explain the relationship between
constructs and tubings.

Putting atomicity in the definition of an ASC-hypergraph is
essential, as it was essential up to now when we dealt with
constructions, and the derived notions of construct and abstract
polytope of a hypergraph. Putting in this definition the other two
properties---saturation and connectedness---is however only a
matter of convenience.

Saturation could be omitted at the cost of having a little bit
more complicated formulations, in which ``$Y$ belongs to the
hypergraph $H$'' is replaced by ``$Y$ is a connected subset of the
hypergraph $H$'', which should mean that $Y$ is a subset of
$\bigcup H$ such that $H_Y$ is connected. By Proposition 4.7, an
arbitrary atomic hypergraph and its saturated closure have the
same constructions. So they do not differ essentially when we deal
with constructions; they yield the same result.

Connectedness too is assumed to organize reasonably the
exposition. It covers what we need most when we deal with abstract
polytopes of hypergraphs. The abstract polytopes of arbitrary
atomic hypergraphs may be derived from the abstract polytopes of
ASC-hypergraphs. The connected components $H_1,\ldots,H_n$ of the
saturated closure of an atomic hypergraph $H$ are ASC-hypergraphs,
and $\A(H)$ is equal to $\A(H_1)\otimes\ldots\otimes\A(H_n)$ (see
Section~5 for the operation $\otimes$). We will return to this
matter at the end of Section~7.

First we have the following, which is proved in a straightforward
manner.

\prop{Remark 6.1}{If $H$ is an ASC-hypergraph and $Y\in H$, then
$H_Y$ is an ASC-hypergraph.}

\noindent For that we rely on Remarks 3.1 and~3.2.

Our goal in this section is to characterize noninductively
constructions of ASC-hypergraphs. Before we start going towards
this goal, we will see in the next proposition how we could define
inductively constructions of ASC-hypergraphs. This is an
adaptation of our inductive definition of a construction of
Section~3 to ASC-hypergraphs specifically. Constructions of
hypergraphs that are not connected do not occur separately any
more---they are incorporated into constructions of connected
hypergraphs. In the proof of this proposition we rely on the fact
that the empty hypergraph is an ASC-hypergraph (see the remarks
concerning the empty hypergraph before Remark 3.1, after Remark
4.1, and in Section~2).

\prop{Proposition 6.2}{If $K$ is a construction of an
ASC-hypergraph $H$, then either $K=\emptyset$, or $K=\{\{x\}\}$,
or} \vspace{-3ex}
\begin{itemize}
\item[(K)] {\it $K=K_1\cup\ldots\cup K_n\cup\{\bigcup H\}$, where
$n\geq 1$, and $\bigcup H=\bigcup H_1\cup\ldots\cup\:\bigcup
H_n\cup\{x\}$ for $\{H_1,\ldots,H_n\}$ being the finest hypergraph
partition of $H_{\ccup H-\{x\}}$, and for every
$i\in\{1,\ldots,n\}$ we have that $K_i$ is a construction of the
ASC-hypergraph $H_i$.}
\end{itemize}

\dkz If $|\bigcup H|=0$, then $K=\emptyset$ by clause (0) of the
definition of a construction. If $|\bigcup H|>0$, then $K$ is
obtained by applying clause (1) of our definition as the last
clause.

If clause (0) preceded immediately this application of clause (1),
then $\bigcup H=\{x\}$, and $K=\{\{x\}\}$. If another application
of clause (1) or an application of clause (2) preceded immediately
this application of clause (1), then we have (K); if that was
another application of (1), then $n=1$, and if that was an
application of (2), then $n\geq 2$.

That $H_i$, which is equal to $H_{\ccup H_i}$, is atomic follows
from the atomicity of $H$, Remark 3.1 and $\bigcup H_i\in H$,
which holds because $H$ is saturated. We can then conclude easily
that $H_i$ is an ASC-hypergraph.\qed

\noindent Next we introduce some terminology, and we give some
preliminary lemmata, that lead towards the main goal of this
section.

Let $H$ be an ASC-hypergraph and let $M\subseteq H$. We say that a
subset $S$ of $M$ is an \emph{$M$-antichain} when $|S|\geq 2$ and
for every two distinct members $X$ and $Y$ of $S$ neither
$X\subseteq Y$ nor $Y\subseteq X$.

An $M$-antichain $S$ \emph{misses} $H$ when $\bigcup S\not\in H$.
(This notion may be found in Definition 2.7 of \cite{FK04}.) An
$M$-antichain is \emph{pairwise disjoint} when every pair of
distinct members of it is disjoint.

For the series of lemmata that follows we make all the time the
following assumption:
\begin{itemize}
\item[(M)]{$M$ is a subset of the ASC-hypergraph $H$ such that
every $M$-antichain misses~$H$.}
\end{itemize}
We can prove the following.

\prop{Lemma 6.3}{Every $M$-antichain is pairwise disjoint.}

\dkz Suppose that for some distinct members $X$ and $Y$ of an
$M$-antichain we had $X\cap Y\neq \emptyset$. From $X,Y\in H$ and
the assumption that $H$ is saturated we could then infer that
$X\cup Y\in H$. So $\{X,Y\}$ would make an $M$-antichain that does
not miss $H$, which contradicts~(M).\qed

We can also infer that if (M) holds, then every pair of members
$X$ and $Y$ of $M$ is non-overlapping, which means that either
$X\cap Y= \emptyset$ or $X\subseteq Y$ or $Y\subseteq X$, and
non-adjacent, which means that if $X\cap Y= \emptyset$, then
$X\cup Y\not\in H$ (see Lemmata A1 and A2 of Appendix~A). Note
however that non-overlapping and non-adjacency are binary, and are
tied to tubings, which are subsets of constructions of
\emph{graphs} (see Appendix~A). These two properties do not
suffice for the constructions of hypergraphs in general.

For $X\in M$, an element $x$ of $X$ is called
\emph{$X$-superficial} when for every $Y$ in $M$ that is a proper
subset of $X$ we have $x\not\in Y$. The notion of $X$-superficial
element is relative to $M$, but we need not mention that when $M$
is fixed, as it will be for us most of the time. This notion
corresponds to the notion of root of an f-construction (see the
end of Section~3). We can prove the following.

\prop{Lemma 6.4}{Every $X$ in $M$ has at least one $X$-superficial
element.}

\dkz Suppose that some $X$ in $M$ has no $X$-superficial element.
So for every element $y$ of $X$ there is at least one $Y$ in $M$
that is a proper subset of $X$ such that $y\in Y$. Let
$Y_1,\ldots,Y_m$ be all these sets $Y$, for all the elements $y$
of $X$. Here $m\geq 1$, because $X$ is nonempty (it is an element
of the hypergraph $H$), and we cannot have that $m=1$, because
$Y_1$ is a proper subset of $X$; so $m\geq 2$.

Eliminate from $Y_1,\ldots,Y_m$ every $Y_i$ that is a proper
subset of another $Y_j$ in the sequence, and let the resulting
sequence be $X_1,\ldots,X_n$. Here $n\geq 2$ for the same reasons
that gave $m\geq 2$ above. Since $\{X_1,\ldots,X_n\}$ is an
$M$-antichain, by (M) it should miss $H$, but $X_1\cup\ldots\cup
X_n=X$, and $X\in H$, which is contradictory.
\mbox{\hspace{.5em}}\qed

We have established with this lemma that there is a map $\varphi$
that assigns to every $X$ in $M$ a nonempty set $\varphi(X)$ of
$X$-superficial elements. Next we establish the following.

\prop{Lemma 6.5}{For every $X$ and $Y$ in $M$, if $X\neq Y$, then
$\varphi(X)\cap\varphi(Y)=\emptyset$.}

\dkz Suppose that for $X$ and $Y$ in $M$ we have an $x$ in
$\bigcup H$ such that $x\in\varphi(X)\cap\varphi(Y)$. If we have
that $X\neq Y$, then we have the $M$-antichain $\{X,Y\}$; we
cannot have that $X\subseteq Y$, because then $x$ would not be
$Y$-superficial, and analogously we cannot have that $Y\subseteq
X$.

The $M$-antichain $\{X,Y\}$ does not miss $H$. This is because
from $\varphi(X)\cap\varphi(Y)\neq\emptyset$ we infer $X\cap
Y\neq\emptyset$, and hence $X\cup Y\in H$, by the saturation of
$H$. This contradicts~(M).\qed

Let $\varphi(M)$ be $\{\varphi(X)\mid X\in M\}$. With Lemmata 6.4
and 6.5 we have established that $\varphi(M)$ is a family of
nonempty pairwise disjoint subsets of $\bigcup H$.

From Lemmata 6.4 and 6.5 we also infer the following.

\prop{Lemma 6.6}{The map $\varphi$ is one-one.}

We need the following two general remarks for the proof of the two
lemmata that follow them. Consider a family $\Phi$ of nonempty
pairwise disjoint subsets of a set $W$. We have the following.

\propo{Remark 6.7}{$|\Phi|\leq |W|$.}

\prop{Remark 6.8}{For $W$ finite, we have that $|\Phi|=|W|$ iff
$\Phi=\{\{w\}\mid w\in W\}$.}

\noindent Note that Remark 6.8 does not hold for $W$ infinite. If
$W$ is the set of natural numbers \textbf{N}, then $\Phi$ may be
$\{\{0,1\},\{2\},\{3\},\ldots\}$ or $\{\{0\},\{2\},\{3\},\ldots\}$
with $|\Phi|=|\textbf{N}|$, but $\Phi\neq\{\{n\}\mid n\in
\textbf{N}\}$.

Assume (M) as above, and let the map $\varphi$ and $\varphi(M)$ be
defined as above. We can prove the following.

\prop{Lemma 6.9}{$|M|\leq |\bigcup H|$.}

\dkz As we said above, with Lemmata 6.4 and 6.5 we have
established that $\varphi(M)$ is a family of nonempty pairwise
disjoint subsets of $\bigcup H$. So, by Remark 6.7, we have that
$|\varphi(M)|\leq |\bigcup H|$. By using Lemma 6.6, we have that
$|M|=|\varphi(M)|$, and the lemma follows.\qed

\oprop{Lemma 6.10}{We have that $|M|=|\bigcup H|$ iff
$\varphi(M)=\{\{x\}\mid x\in\bigcup H\}$.}

\dkz With Lemmata 6.4 and 6.5 we have established that
$\varphi(M)$ is a family of nonempty pairwise disjoint subsets of
$\bigcup H$. By using Lemma 6.6, we have that $|M|=|\bigcup H|$
iff $|\varphi(M)|=|\bigcup H|$. Since moreover $\bigcup H$ is
finite, we obtain the desired conclusion by Remark~6.8.\qed

The following proposition gives our noninductive characterization
of constructions of ASC-hypergraphs.

\prop{Proposition 6.11}{We have \mbox{\rm (M)} and $|M|=|\bigcup
H|$ iff $M$ is a construction of~$H$.}

\dkz From left to right we proceed as follows. If $|\bigcup H|=0$,
then $H=\emptyset$ and $M=\emptyset$. If $|\bigcup H|=1$, then
$H=\{\{x\}\}$ and $M=\{\{x\}\}$.

We establish next that if $|\bigcup H|\geq 1$, then $\bigcup H\in
M$. For $|\bigcup H|=1$, we established that in the preceding
paragraph. Suppose $|\bigcup H|=k\geq 2$ and $\bigcup H\not\in M$.
Then by omitting from $M$ members that are proper subsets of other
members we can obtain an $M$-antichain $Y_1,\ldots,Y_m$. We must
have that $m\geq 2$ because of Lemma 6.10 and $\bigcup H\not\in
M$. However, by Lemma 6.10 and the saturation and connectedness of
$H$, we have that $Y_1\cup\ldots\cup Y_m=\bigcup H\in H$, and
hence our $M$-antichain does not miss $H$, which contradicts~(M).

We proceed next by induction on $|\bigcup H|$, the basis being the
case above when $|\bigcup H|=1$. So let $|\bigcup H|\geq 2$, and
let $x$ be the $\bigcup H$-superficial element, which exists by
Lemma 6.10. Then for $\{H_1,\ldots,H_n\}$, where $n\geq 1$, being
the finest hypergraph partition of $H_{\ccup H-\{x\}}$ we have
that $M=M_1\cup\ldots\cup M_n\cup\{\bigcup H\}$, where for every
$i\in\{1,\ldots,n\}$ the set $M_i$ is $M_{\ccup H_i}$. We show
that
\begin{itemize}
\item[(M$_i$)]{$M_i$ is a subset of the ASC-hypergraph $H_i$ such
that every $M_i$-antichain misses~$H_i$.}
\end{itemize}
It follows immediately that $M_i\subseteq H_i$, and $H_i$ is an
ASC-hypergraph by Remark 6.1. Every $M_i$-antichain is an
$M$-antichain, and misses $H$ by assumption; hence it must also
miss $H_i$. So we have~(M$_i$).

We have that $|M_i|=|\bigcup H_i|$ for the following reason. For
every $i\in\{1,\ldots,n\}$ we must have that $|M_i|\leq|\bigcup
H_i|$ by Lemma 6.9, and if for some $j\in\{1,\ldots,n\}$ we had
$|M_i|<|\bigcup H_i|$, then we could not secure that
\[|M_1|\pl\ldots\pl|M_n|\pl 1=|M|=\mbox{$|\bigcup H|=|\bigcup
H_1|\pl\ldots\pl|\bigcup H_n|\pl 1.$}\]

By the induction hypothesis, we may conclude that $M_i$ is a
construction of the ASC-hypergraph $H_i$, and it follows that $M$
is a construction of~$H$.

From right to left we proceed by induction on the size of the
construction $K$ of the ASC-hypergraph $H$. Consider what $K$ may
be according to Proposition 6.2. If $K=\emptyset$ or
$K=\{\{x\}\}$, then it is trivial that every $K$-antichain misses
$H$ (there are no $K$-antichains) and that $|K|=|\bigcup H|$. If
we have (K) as in Proposition 6.2, then the induction hypothesis
applies to the constructions $K_i$ for every $i\in\{1,\ldots,n\}$,
and by it we conclude that
\begin{itemize}
\item[(K$_i$)]{$K_i$ is a subset of the ASC-hypergraph $H_i$ such
that every $K_i$-antichain misses~$H_i$}
\end{itemize}
and $|K_i|=|\bigcup H_i|$.

If a $K$-antichain is a $K_i$-antichain, then (K$_i$) applies to
it; hence it misses $H_i$, and it follows that it misses~$H$.

Suppose we have a $K$-antichain $S$ that is not a $K_i$-antichain
for any $i\in\{1,\ldots,n\}$. The $\bigcup H$-superficial element
$x$ relative to $K$ cannot be in $\bigcup S$, because otherwise
$\bigcup H$ would have to belong to $S$, and this is impossible
(every member of $S$ is a subset of~$\bigcup H$).

Since $S$ is not a $K_i$-antichain, there must be two distinct
elements $y$ and $z$ of $\bigcup S$ such that $y\in \bigcup H_p$
and $z\in \bigcup H_q$ for two distinct members $H_p$ and $H_q$ of
the finest hypergraph partition $\{H_1,\ldots,H_n\}$ of $H_{\ccup
H-\{x\}}$, which we have according to (K). We conclude that
$\bigcup S\not\in H_{\ccup H-\{x\}}$, and since $x\not\in\bigcup
S$, and hence $\bigcup S\subseteq \bigcup H\mn\{x\}$, it follows
that $\bigcup S\not\in H$.

So every $K$-antichain misses $H$, and $|K|=|\bigcup H|$, by
Remark 5.1.\qed

The following proposition will help us to obtain a
characterization of the notion of construct in the style of
Proposition~6.11.

\prop{Proposition 6.12}{We have \mbox{\rm (M)} iff for some
construction $K$ of $H$ we have that $M\subseteq K$.}

\dkz The direction from right to left is obtained easily as
follows from Proposition 6.11. Suppose for some construction $K$
of $H$ we have that $M\subseteq K$. By Proposition 6.11 from right
to left we have that every $K$-antichain misses $H$. Every
$M$-antichain is however a $K$-antichain.

For the other direction, suppose (M). If we could prove that
\begin{itemize}
\item[(R)]{there is a subset $K$ of $H$ such that $M\subseteq K$,
$|K|=|\bigcup H|$ and every $K$-antichain misses $H$,}
\end{itemize}
then by Proposition 6.11 from left to right we would have that $K$
is a construction of $H$, and we would obtain the right-hand side
of the proposition we are proving. The remainder of our proof is
an inductive proof of~(R).

We have that $|M|\leq|\bigcup H|$ by Lemma 6.9. Our proof of (R)
will proceed by induction on $|\bigcup H|\mn |M|$. In the basis,
when this number is zero, and hence $|M|=|\bigcup H|$, we take
$K=M$.

Suppose for the induction step that $|M|<|\bigcup H|$. Let
$M^+=M\cup\{\bigcup H\}$. From (M) we easily infer that every
$M^+$-antichain misses $H$, since $\bigcup H$ is not a member of
any $M^+$-antichain, and hence $M^+$-antichains are
$M$-antichains. If $\bigcup H\not\in M$, then $|M|<|M^+|\leq
|\bigcup H|$, and we may apply the induction hypothesis to $M^+$;
namely, we have (R) for $M$ replaced by $M^+$. For the set $K$
that this yields we have that $M^+\subseteq K$, and hence
$M\subseteq K$, which gives~(R).

Suppose that $\bigcup H\in M$. For every $x$ in $\bigcup H$ we
have then a member of $M$, namely $\bigcup H$, to which $x$
belongs. We easily infer that there is hence a member $X$ of $M$
such that $x$ is $X$-superficial. If that is not $\bigcup H$, then
we pass to a proper subset $Y$ of $\bigcup H$ such that $Y\in M$
and $x\in Y$, and continue in this manner until we reach $X$
(officially, a trivial induction on the number of members of $M$
to which $x$ belongs is here at work; this number is, of course,
finite). So for every $x$ in $\bigcup H$ there is an $X$ in $M$
such that $x\in\varphi(X)$.

Since $|M|<|\bigcup H|$, we have that $\bigcup H$ is not empty,
and for an $x$ in $\bigcup H$ there is an $X$ in $M$ such that
$x\in\varphi(X)$ and $|\varphi(X)|\geq 2$. If we always had
$|\varphi(X)|=1$, then, by Lemma 6.10, we would obtain
$|M|=|\bigcup H|$, which contradicts $|M|<|\bigcup H|$. So we have
that $x,y\in\varphi(X)$ and $x\neq y$.

Let $\{H_1,\ldots,H_n\}$, for $n\geq 1$, be the finest hypergraph
partition of $H_{X-\{x\}}$. This partition is nonempty, because
for some $i\in\{1,\ldots,n\}$ we must have that $y\in\bigcup H_i$.
Note that $\bigcup H_i\not\in M$; otherwise, $y\not\in\varphi(X)$,
i.e.\ $y$ would not be $X$-superficial. Let $M'$ be
$M\cup\{\bigcup H_i\}$. Since $\bigcup H_i\not\in M$, we have that
$|M|<|M'|$.

We prove next that every $M'$-antichain misses $H$. (This will
occupy us for most of the remainder of the proof.) Suppose there
is an $M'$-antichain $S$ that does not miss $H$. We must have that
$\bigcup H_i\in S$; otherwise, $S$ would be an $M$-antichain, and
we would contradict~$M$.

Let $S'=\{Y\in S\mid Y\not\subseteq X\}$. If $S'\neq\emptyset$,
then we show that (M) does not hold. Take $S'\cup\{X\}$; this is
an $M$-antichain, because, first, all its members are from $M$ (we
have that $\bigcup H_i\subseteq X$), and, secondly, $X$ cannot be
a subset of a $Y$ in $S'$; otherwise, we would have that $\bigcup
H_i\subseteq Y$, and $S$ would not be an $M'$-antichain.

We show next that $\bigcup(S'\cup\{X\})\in H$, which will
contradict (M). Since every $Z$ in $S\mn S'$ is a subset of $X$,
we have that $\bigcup(S'\cup\{X\})=\bigcup(S\cup\{X\})=\bigcup
S\cup X$. We have that $y\in\bigcup H_i\subseteq\bigcup S\cap X$,
and so by the saturation of $H$, we have that $\bigcup S\cup X\in
H$; so $\bigcup(S'\cup\{X\})\in H$, which contradicts~(M).

So we have that $S'=\emptyset$, and hence $\bigcup S\subseteq X$.
We have that $X\not\in S$; otherwise, $S$ would not be an
$M'$-antichain (we have that $\bigcup H_i\subseteq X$). We also
have that $x\not\in\bigcup S$. Otherwise, since $x\not\in\bigcup
H_i$, we would have a $Z$ in $M$ such that $Z\in S$ and $x\in Z$.
But $Z\subseteq X$ and $Z\neq X$, as we have just shown, and hence
$x$ is not $X$-superficial relative to $M$, which is a
contradiction.

So $\bigcup S\subseteq X\mn\{x\}$, and since $\bigcup S\in H$, we
have that $\bigcup S\in H_{X-\{x\}}$. Since  $\{H_1,\ldots,H_n\}$
is the finest hypergraph partition of $H_{X-\{x\}}$ and $\bigcup
H_i\in S$, we obtain $\bigcup S\in H_i$. Hence $\bigcup
H_i\subseteq\bigcup S$ and $\bigcup S\subseteq\bigcup H_i$, which
means that $\bigcup S=\bigcup H_i$, and since $\bigcup H_i\in S$,
this contradicts the assumption that $S$ is an $M'$-antichain.

So every $M'$-antichain misses $H$, and we may apply the induction
hypothesis to $M'$; namely, we have (R) for $M$ replaced by $M'$.
For the set $K$ that this yields we have that $M'\subseteq K$, and
hence $M\subseteq K$, which gives~(R).

So we have established (R), and, as explained above, with that we
may end our proof.\qed

As an immediate corollary of Proposition 6.12 we have the
following.

\prop{Proposition 6.13}{We have \mbox{\rm (M)} and $\bigcup H\in
M$ iff $M$ is a construct of~$H$.}

We conclude this section with the following technical lemma, which
we need for the proof of Lemma 9.5 in Section~9.

\prop{Lemma 6.14}{If $K$ is a construction of the ASC-hypergraph
$H$ and $Y\in H\mn K$, then there is an $X$ in $K$ such that
$Y\subset X$ and for $x$ being $X$-superficial we have $x\in Y$.}

\dkz If $K=\emptyset$ or $K=\{\{x\}\}$, then $H\mn K$ is empty,
and the lemma holds trivially.

Suppose $K$ is of the form specified by (K) of Proposition 6.2 and
$Y\in H\mn K$. Then either we find in $Y$ the $\bigcup
H$-superficial element $x$, and we are done, or there is an
$i\in\{1,\ldots,n\}$ such that $Y\in H_i\mn K_i$; otherwise,
$\{H_1,\ldots,H_n\}$ would not be the finest hypergraph partition
of $H_{\ccup H-\{x\}}$. We may then proceed by induction.\qed

\section{\large \textbf{Continuations of constructions}}
In this section we describe the vertices of the abstract polytope
of a hypergraph in terms of a partial operation on constructions
that we call continuation, which may be understood intuitively as
indeed the continuation of one construction by another. This
description leads to a characterization of all the faces of
abstract polytopes of hypergraphs, and in particular of their
facets. This characterization of facets will play an important
role in the next section, where we prove that abstract polytopes
of hypergraphs are indeed abstract polytopes. The results of this
section are closely related to the results of \cite{Z06}, though
the presentation is different. (In Definition 3.1 of \cite{Z06},
which introduces a notion that plays a role analogous to our $\uhy
L$ below, one should require in (3.3) that $I$ is nonempty.)

We start with some preliminary matters. Assume for the proposition
below that $L$ is a construction of the ASC-hypergraph $H$ and
that $Y\in L$. Then we can prove the following.

\prop{Proposition 7.1}{The set $L_Y$ is a construction of the
ASC-hypergraph~$H_Y$.}

\dkz We know that $H_Y$ is an ASC-hypergraph by Remark 6.1. It is
clear that $L_Y\subseteq H_Y$. We show next that every
$L_Y$-antichain misses $H_Y$ (see Section~6). If for an
$L_Y$-antichain $S$ we had that $\bigcup S\in H_Y$, then we would
also have $\bigcup S\in H$. Since $S$ is also an $L$-antichain, we
would have, according to Proposition 6.11 from right to left, that
$L$ is not a construction of $H$, which contradicts our
assumption. So every $L_Y$-antichain misses~$H_Y$.

It remains to show that $|L_Y|=|\bigcup H_Y|$ in order to obtain
the proposition by applying Proposition 6.11 from left to right.
We know by Lemma 6.10 and $|L|=|\bigcup H|$ that
$\varphi(L)=\{\{x\}\mid x\in\bigcup H\}$. Let
$\varphi(L_Y)=\{\varphi(X)\mid X\in L_Y\}$. It is easy to see that
$\varphi(L_Y)=\{\{y\}\mid y\in Y\}$ since $\bigcup H_Y=Y$. So
$|\varphi(L_Y)|=|Y|$, and, since $\varphi$ is one-one, we have
$|L_Y|=|Y|=|\bigcup H_Y|$. \qed

For $H$ a hypergraph and $X\subseteq\bigcup H$ let
\[
X\cup_H Y=\left\{
\begin{array}{ll}
X\cup Y & {\mbox{\rm if }} X\cup Y\in H,
\\
X & {\mbox{\rm if }} X\cup Y\not\in H.
\end{array}
\right.
\]
Then, with assumptions as before the preceding proposition, we
have the following lemma.

\prop{Lemma 7.2}{If $X\in L\mn L_Y$, then $(X\mn Y)\cup_H Y=X$.}

\dkz If $Y\subseteq X$, then we clearly have
\[(X\mn Y)\cup_H Y=(X\mn Y)\cup Y=X.\]
Suppose not $Y\subseteq X$. It is impossible that $X\subseteq Y$,
because $X\in L\mn L_Y$. So $\{X,Y\}$ is an $L$-antichain (see
Section~6). Since $X\cap Y=\emptyset$ (see Lemma 6.3), we have
that $X\mn Y=X$, and since $X\cup Y$ cannot belong to $H$ (our
$L$-antichain misses $H$; see Section~6), we have that $X\cup_H
Y=X$. \qed

For $F\subseteq\P C$ and $Z\subseteq C$ let
\[_Z F=_{df}\{X\cap Z\mid X\in F\; \& \; X\cap Z\neq\emptyset\}.\]
If $F$ is a hypergraph $H$, and $Z\subseteq \bigcup H$, then it is
easy to check that $_Z H$ is a hypergraph on $Z$. We have always
that $H_Z$ (which is defined at the beginning of Section~3) is a
subset of $_Z H$. The converse need not hold, and here is a
counterexample for that.

Consider the ASC-hypergraph
\[H=\{\{x\},\{y\},\{z\},\{u\},\{x,y,z\},\{y,z,u\},\{x,y,z,u\}\}\]
on $\{x,y,z,u\}$, and let $Z=\{x,y,z\}$ (this hypergraph $H$ is
the hypergraph $H_{4021}$ of Appendix~B). Then we have
\begin{tabbing}
\hspace{10.5em}\=$H_Z\:$\=$=\{\{x\},\{y\},\{z\},\{x,y,z\}\}$,
\\
\>$_Z H$\>$=H_Z\cup\{\{y,z\}\}$.
\end{tabbing}
There are simpler counterexamples when $H$ is not an
ASC-hypergraph, but we wanted to have a counterexample with such a
hypergraph.

We can verify the following in a straightforward manner.

\prop{Remark 7.3}{If $H$ is an ASC-hypergraph and
$Z\subseteq\bigcup H$, then $_Z H$ is an ASC-hypergraph on $Z$.}

\noindent This remark should be compared with Remark 6.1

Let us assume, as before Proposition 7.1, that $L$ is a
construction of the ASC-hypergraph $H$ and that $Y\in L$. It is
then easy to see that $\uhy L={\uhy (L\mn L_Y)}$, and we will rely
on this equation without notice in the rest of this section. We
have the following.

\prop{Proposition 7.4}{The set $\uhy L$ is a construction of the
ASC-hypergraph $\uhy H$.}

\dkz We know that $\uhy H$ is an ASC-hypergraph by Remark 7.3. It
is clear that $\uhy L\subseteq\: \uhy H$. We show next that every
${\uhy L}$-antichain misses $\uhy H$. Suppose for such an
antichain $\{X_1, \ldots, X_k\}$, where $k\geq 2$, we had that
$X_1\cup\ldots\cup X_k\in\: \uhy H$.

For every $Z\in\: \uhy L$ we have that $Z=X\mn Y$ for some $X\in
L\mn L_Y$. Then by Lemma 7.2 we have that $Z\cup_H Y=X$.

Consider then the set $S=\{X_1\cup_H Y,\ldots, X_k\cup_H Y\}$.
Since for every $i\in\{1,\ldots,k\}$ we have that $X_i\in\: \uhy
L$, we can conclude as above that $X_i\cup_H Y=X$ for some $X\in
L\mn L_Y$, and hence $X_i\cup_H Y\in L$. It follows that $S$ is an
$L$-antichain.

We have assumed above that $X_1\cup\ldots\cup X_k=W-Y$ for some
$W\in H$. Suppose (1) for every $i$ we have that $X_i\cup_H Y=X_i$
and (2) $W\mn Y=W$. Then $S$ is an $L$-antichain that does not
miss $H$, which together with Proposition 6.11 from right to left
contradicts the assumption that $L$ is a construction of $H$.

Suppose (1) and not (2), i.e. $W\mn Y\neq W$. Then
\begin{tabbing}
\hspace{10.5em}$\bigcup S \cup Y$ \= $=X_1\cup\ldots\cup X_k\cup
Y$
\\
\>$=(W\mn Y)\cup Y$
\\
\>$=W\cup Y$.
\end{tabbing}
Since $W\cap Y\neq\emptyset$, as we supposed above, we have that
$W\cup Y\in H$, because $H$ is saturated. Then $S\cup\{Y\}$ is an
$L$-antichain that does not miss $H$, which is contradictory, as
above.

Suppose not (1), i.e.\ for some $i$ we have that $X_i\cup_H
Y=X_i\cup Y\in H$. Then $\bigcup S=W\cup Y$. Since $X_i\subseteq
W$ and $X_i\neq\emptyset$, we have that $W\cap(X_i\cup
Y)\neq\emptyset$. So $W\cup(X_i\cup Y)\in H$, because $H$ is
saturated, and hence $W\cup Y\in H$. So $S$ is an $L$-antichain
that does not miss $H$, which is contradictory, as above. Thereby,
we have shown that every $\uhy L$-antichain misses $\uhy H$.

We have shown in the proof of Proposition 7.1 that $|L_Y|=|Y|$. It
remains to show that $|\uhy L|=|\bigcup H\mn Y|$ (we have that
$\bigcup H\mn Y=\bigcup\: \uhy H$) in order to obtain the
proposition by applying Proposition 6.11 from left to right. We
have that $| \uhy L|=|L\mn L_Y|=|L|\mn |L_Y|=|\bigcup H|-|Y|$,
because $L$ is a construction of $H$, and hence $|L|=|\bigcup H|$
(see Remark 5.1), and because $|L_Y|=|Y|$, as we showed in the
proof of Proposition 7.1. We have that $|\bigcup H|\mn|Y|=|\bigcup
H\mn Y|$, which concludes our proof. \mbox{\hspace{.5em}}\qed

Let $H$ be an ASC-hypergraph. For $Y\in H$, let $K$ and $J$ be
respectively constructions of the ASC-hypergraphs $H_Y$ and $\uhy
H$. Then we define the \emph{continuation} $K\ast J$ of $K$ by $J$
in the following way:
\[K\ast J=K\cup\{X\cup_H Y\mid X\in J\}.\]
Here is an example of continuation.

Take the ASC-hypergraph
\[
\bar{A}=\{\{x\},\{y\},\{z\},\{u\},\{x,y\}, \{y,z\},\{z,u\},
\{x,y,z\},\{y,z,u\},\{x,y,z,u\}\},\] which is the saturated
closure of the hypergraph $A$ of Section~3. (This hypergraph is
called $H'_{4321}$ in Appendix~B.) Let $K$ be the construction
$\{\{u\},\{z,u\}\}$ of the ASC-hypergraph
$\bar{A}_{\{z,u\}}=\{\{z\},\{u\},\{z,u\}\}$, and let $J$ be the
construction $\{\{x\},\{x,y\}\}$ of the ASC-hypergraph
$_{\{x,y\}}\bar{A}=\:
_{\{x,y,z,u\}-\{z,u\}}\bar{A}=\{\{x\},\{y\},\{x,y\}\}$. Then the
continuation $K\ast J$ is $\{\{u\},\{z,u\}, \{x\},\{x,y,z,u\}\}$.
(The construction $K\ast J$ is the construction $N$ of Section~5,
which corresponds to the s-construction $y(x\pl(zu))$; with this
s-construction we have labelled one of the vertices of the
associahedron $K_5$ in Appendix~B.)

Another example of continuation with the same ASC-hypergraph
$\bar{A}$ is obtained by taking $Y$ to be $\{y,z,u\}$, with $K$
being the construction $\{\{y\},\{u\},$ $\{y,z,u\}\}$ of the
ASC-hypergraph
\[\bar{A}_{\{y,z,u\}}=\{\{y\},\{z\},\{u\},\{y,z\},\{z,u\},\{y,z,u\}\},\]
and $J$ being the construction $\{\{x\}\}$ of the ASC-hypergraph
\[_{\{x\}}\bar{A}=\: _{\{x,y,z,u\}-\{y,z,u\}}\bar{A}=\{\{x\}\}.\]
Then $K\ast J$ is $\{\{y\},\{u\},\{y,z,u\},\{x,y,z,u\}\}$. (This
is the construction $M$ of Section~3, which corresponds to the
s-construction $xz(y\pl u)$.)

A third example of continuation is provided by taking $\bar{A}$,
$Y$ and $K$ as in the first example and $J$ as
$\{\{y\},\{x,y\}\}$. Then $K\ast J$ is
$\{\{u\},\{z,u\},\{y,z,u\},$ $\{x,y,z,u\}\}$. (This is the
construction $L$ of Section~3, which corresponds to the
s-construction $xyzu$.)

A fourth and final example is with everything being as in the
second example except for $K$, which is now
$\{\{u\},\{z,u\},\{y,z,u\}\}$. Then $K\ast J$ is equal to the
$K\ast J$ of the preceding example (namely, to the
construction~$L$).

Note that if $Y=\bigcup H$, then $\bigcup H\mn Y=\emptyset$, and
we have that $_\emptyset H=\emptyset$ and $J=\emptyset$. In that
case $K\ast\emptyset=K$. We shall however need $\ast$ mostly when
$J$ is not~$\emptyset$.

It is easy to see that $|K\ast J|=|K|\pl|J|$. This matches the
fact that $|\bigcup H|=|Y|\pl|\bigcup H\mn Y|$, as the
propositions below will show. We can then prove the following.

\prop{Proposition 7.5}{The set $K\ast J$ is a construction
of~$H$.}

\dkz We prove first that $K\ast J\subseteq H$. Let $X\in K\ast J$.
If $X\in K$, then $X\in H_Y$, and hence $X\in H$. Suppose
$X=X'\cup_H Y$ for $X'\in J$. If $X'\cup Y\in H$, then we are
done. If $X'\cup Y\not\in H$, then we reason as follows.

Since $X'\in\: \uhy H$, we have that $X'=W\mn Y$ for some $W\in
H$. If $W\mn Y\neq W$, then $W\cap Y\neq\emptyset$, and $W\cup
Y\in H$, because $H$ is saturated; so $X'\cup Y=W\cup Y$, which
contradicts the assumption that $X'\cup Y\not\in H$. So we must
have that $W\mn Y=W$, and hence $X'=W$. From $X'\cup Y\not\in H$,
we have that $X'\cup_H Y=X'$, and so $X\in H$. This proves that
$K\ast J\subseteq H$.

We show next that every $K\ast J$-antichain misses $H$. Suppose
not, and let $S$ be a $K\ast J$-antichain such that $\bigcup S\in
H$.

If $|\uhy S|=0$, then $S$ is a $K$-antichain that does not miss
$H_Y$, which together with Proposition 6.11 from right to left
contradicts our assumption that $K$ is a construction of~$H_Y$.

If $|\uhy S|\geq 2$, then $\uhy S$ is a $J$-antichain that does
not miss $\uhy H$, because $\bigcup(\uhy S)=\bigcup S\cap(\bigcup
H\mn Y)$. To show that $\uhy S$ is a $J$-antichain it is
sufficient to show that if $X_1\subseteq X_2$, then $X_1\cup_H
Y\subseteq X_2\cup_H Y$. This holds because it is impossible that
$X_1\cup_H Y=X_1\cup Y$, while $X_2\cup_H Y=X_2$; if $X_1\cup Y\in
H$, then $X_2\cup Y\in H$, since $X_1\subseteq X_2$ and $H$ is
saturated.

If $|\uhy S|=1$, then $S$ would not be a $K\ast J$-antichain for
the following reason. There would be a unique member $X$ of $S$ of
the form $X'\cup_H Y$ for $X'\in J$. We wish to show that
$X'\cup_H Y=X'\cup Y$, and that will be the case when $X'\cup Y\in
H$. We must have that $S\mn\{X\}\neq \emptyset$, because $S$ must
have at least two members. Since $\bigcup(S\mn\{X\})$ is a
nonempty subset of $Y$, we must have that $(\bigcup(S\mn\{X\})\cup
X)\cap Y\neq\emptyset$ while $\bigcup(S\mn\{X\})\cup X=\bigcup
S\in H$ and $Y\in H$. We then obtain that $\bigcup(S\mn\{X\})\cup
X\cup Y=X\cup Y\in H$, because $H$ is saturated. So $(X'\cup_H
Y)\cup Y=X'\cup Y\in H$, and so $X'\cup_H Y=X'\cup Y$. Every
member $Z$ of $S\mn\{X\}$ is a subset of $Y$, and so $Z\subseteq
X'\cup Y$. So $S$ is not a $K\ast J$-antichain.

It remains to appeal to $|K\ast J|=|\bigcup H|$. This follows from
the equations we mentioned before stating the proposition,
together with $|K|=|Y|$ and $|J|=|\bigcup H\mn Y|$, which hold
because $K$ and $J$ are constructions of $H_Y$ and $\uhy H$
respectively. We conclude that $K\ast J$ is a construction by
applying Proposition 6.11 from left to right. \qed

With the assumptions stated before Proposition 7.1 we have the
following lemma.

\prop{Lemma 7.6}{We have that $L_Y\ast\: \uhy L=L$.}

\dkz If $Y=\bigcup H$, then $L_Y=L$ and $\uhy L=\emptyset$. We
then have that $L\ast\emptyset=L$.

Suppose $Y\subset\bigcup H$. To show that $L_Y\ast\: \uhy
L\subseteq L$, suppose $X\in L_Y\ast \uhy\: (L\mn L_Y)$. If $X\in
L_Y$, then $X\in L$. If, on the other hand, $X=X'\cup_H Y$ for
some $X'\in\: \uhy L$, then $X'=X''\mn Y\neq\emptyset$ for some
$X''\in L\mn L_Y$. So $X=(X''\mn Y)\cup_H Y$, and, by Lemma 7.2,
we obtain that $X=X''$. So $X\in L\mn L_Y$, and hence $X\in L$.
Therefore we have established that $L_Y\ast\: \uhy L\subseteq L$.

For the converse inclusion, suppose $X\in L$. Then either $X\in
L_Y$ or $X\in L\mn L_Y$. If $X\in L_Y$, then $X\in L_Y\ast\: \uhy
L$. If $X\in L\mn L_Y$, then our purpose is to show that
$X=X'\cup_H Y$ for some $X'\in\: \uhy L$, which means that
$X'=X''\mn Y\neq\emptyset$ for some $X''\in L'\mn L_Y$. So it is
enough to establish that $X=(X\mn Y)\cup_H Y$, which we have by
Lemma~7.2. \mbox{\hspace{.5em}}\qed

With the assumptions stated before Proposition 7.5 we easily prove
the following.

\prop{Lemma 7.7}{We have that $(K\ast J)_Y=K$ and $\uhy (K\ast
J)=J$.}

\noindent For the second equation we rely on the equation
$(X\cup_H Y)\mn Y=X$ for $X\subseteq \bigcup H\mn Y$.

We can then prove the following.

\prop{Proposition 7.8}{For a given construction $L$ of the
ASC-hypergraph $H$, and a given $Y\in L$, the constructions $L_Y$
and $\uhy L$ are the unique constructions $K$ and $J$ of the
ASC-hypergraphs $H_Y$ and $\uhy H$ respectively such that $K\ast
J=L$.}

\dkz We rely on Propositions 7.1 and 7.4, and on Lemma 7.6, to
obtain that for the constructions $K=L_Y$ and $J=\: \uhy L$ of
$H_Y$ and $\uhy H$ respectively we have that $K\ast J=L$. For
uniqueness suppose that $K\ast J= K'\ast J'$. Then, by relying on
Lemma 7.7, we obtain that $K=K'$ and $J=J'$. \qed

For every ASC-hypergraph $H$, every facet of the abstract polytope
$\A(H)$ of $H$ (see Section~5) is of the form $\{Y,\bigcup H\}$
for $Y$ a member of $H\mn\{\bigcup H\}$, provided the rank of
$\A(H)$ is at least 1. Take the section $\{Y,\bigcup H\}/F_{-1}$
of $\A(H)$, i.e.\ the set of all the faces of $\A(H)$ below, i.e.\
including, $\{Y,\bigcup H\}$ and above, i.e.\ included in,
$F_{-1}$. The least face $F_{-1}$ of $\A(H)$ is $\bar{H}^*$, but
here $\bar{H}=H$, since $H$ is saturated.

The vertices of $\A(H)$ in $\{Y,\bigcup H\}/F_{-1}$ are all the
constructions of $H$ in which $Y$ is a member. These are the
vertices in which the facet $\{Y,\bigcup H\}$ is included, i.e.\
the vertices incident with this facet. Each such vertex $L$ is
equal to $L_Y\ast\: \uhy L$, and, according to Proposition 7.8,
this is the only way to represent $L$ as a continuation of
constructions of $H_Y$ and $\uhy H$.

The faces of $\A(H)$ in $\{Y,\bigcup H\}/F_{-1}$ are all the
constructs of $H$ in which $Y$ is a member, together with
$\bar{H}^*$ as an additional face. These are the faces in which
the facet $\{Y,\bigcup H\}$ is included, i.e.\ the faces incident
with this facet.

The abstract polytope $\{Y,\bigcup H\}/F_{-1}$ is isomorphic to
$\A(H_Y)\otimes \A( \uhy H)$ (see Section~5 for the product
$\otimes$). We may take $\{Y,\bigcup H\}/F_{-1}$ of $\A(H)$ as
being the result of a partial operation $\ast$ applied to
$\A(H_Y)$ and $\A( \uhy H)$, akin to $\otimes$, but different from
it.

We define $\A(H_Y)\ast\A( \uhy H)$ as the set of all the
constructs of $H$ in which $Y$ is a member, together with
$\bar{H}^*$ as an additional element. The reason for writing this
operation $\ast$ is that every construction $L$ in $\A(H_Y)\ast\A(
\uhy H)$ is of the form $L_Y\ast\: \uhy L$ for $L_Y$ in $\A(H_Y)$
and $\uhy L$ in $\A( \uhy H)$; this continuation $\ast$ operation
on constructions induces analogous continuation $\ast$ operations
for all the other constructs in $\A(H_Y)\ast\A( \uhy H)$. Each of
these constructs $C$ is the result of applying a $\ast$ to a
construct $C_1$ in $\A(H_Y)$ and a construct $C_2$ in $\A( \uhy
H)$. The presence of $Y$ in $C$ guarantees the existence of $C_1$,
and the presence of $\bigcup H$ in $C$ guarantees the existence of
$C_2$, in which we have $\bigcup H\mn Y$. (The set $Y$, as a
member of $H$, is nonempty, and since it is a proper subset of
$\bigcup H$, we have that $\bigcup H\mn Y$ is nonempty too.)

For an ASC-hypergraph $H$ with $\A(H)$ of rank $r\geq 1$, we can
construct $\A(H)$ inductively in terms of abstract polytopes of
lower rank in the following manner. For every member $Y$ of
$H\mn\{\bigcup H\}$, take the abstract polytopes $\A(H_Y)$ and
$\A( \uhy H)$, which are of rank lower than $r$; the rank of
$\A(H_Y)$ is $|Y|\mn 1$, and the rank of $\A( \uhy H)$ is
$|\bigcup H\mn Y|\mn 1$. Then take the union of all the abstract
polytopes $\A(H_Y)\ast\A( \uhy H)$ (note that they all have the
same least face $F_{-1}$, which is $\bar{H}^*$), and add as a new
face $\{\bigcup H\}$ as the greatest face. This is $\A(H)$.

The basis of this induction is given by the abstract polytopes
$\A(H)=\{\{x\},\{x,*\}\}$ of rank 0, where $H$ is $\{\{x\}\}$.
Note that (as we said at the end of Section~5) if $H=\emptyset$,
then $\A(H)=\{\emptyset,\{*\}\}$, which is also of rank 0, but is
not needed for the basis of our induction, because $\emptyset$
cannot be a member of a hypergraph.

For every atomic hypergraph $H$, we can obtain $\A(H)$ as
$\A(H_1)\otimes\ldots\otimes\A(H_n)$, for $n\geq 1$, and
$\{H_1,\ldots,H_n\}$ being the finest hypergraph partition of the
saturated closure of $H$. (If $n=1$, then
$\A(H_1)\otimes\ldots\otimes\A(H_n)$ is, of course, just
$\A(H_1)$.) The hypergraphs $H_1,\ldots, H_n$ are ASC-hypergraphs.

So with the help of the operation $\otimes$, and of the related
operation $\ast$, on hypergraphs we can define inductively the
hypergraph $\A(H)$ for an arbitrary hypergraph $H$. If $n>1$, then
for every $i\in\{1,\ldots,n\}$, the rank of $\A(H_i)$ is strictly
smaller than the rank of $\A(H)$. We deal with this matter in more
detail in the next section.

\section{\large \textbf{Abstract polytopes of hypergraphs are
abstract polytopes}}
In this section we verify that for an
arbitrary atomic hypergraph $H$ we have that the abstract polytope
$\A(H)$ is indeed an abstract polytope of rank $|\bigcup H|\mn n$,
where $n$ is the connectedness number of $H$ (see the beginning of
Section~5). This result follows from the geometric representation
of $\A(H)$ as a convex polytope in Euclidean space, with which we
deal in Section~9, but, since that section is rather concise, we
prefer to give an independent direct proof in the present section,
and make our treatment of the abstract polytopes $\A(H)$
self-contained.

Our general notion of abstract polytope is as in \cite{MS02}
(Section 2A). The definition given there, which we summarize now,
says that an abstract polytope of rank $r$, for $r\geq -1$, is a
partially ordered set $\langle P,\leq\rangle$, with \emph{carrier}
$P$, that satisfies four properties (P1), $\ldots$~,(P4).

Property (P1) is that $P$ has a least and a greatest \emph{face},
i.e.\ element, denoted by $F_{-1}$ and $F_r$ (they need not be
distinct).

Property (P2) is that each \emph{flag} of $P$, i.e.\ maximal
linearly ordered subset of $P$, contains exactly $r\pl 2$ faces
including $F_{-1}$ and $F_r$. A partially ordered set $P$ that
satisfies (P1) and (P2) is said to be a partially ordered set of
rank $r$.

Both (P1) and (P2) are easily verified for $\A(H)$. The
verification of the remaining two properties (P3) and (P4) is a
more difficult task, for which we will first reformulate the
properties within an inductive definition of abstract polytope of
rank $r$, the induction being based on $r$. Let us first state the
remaining properties as they are formulated in \cite{MS02}
(Section 2A).

For property (P3) we need some preliminary notions. A partially
ordered set $P$ of rank $r$ is said to be \emph{connected} when
either $r\leq 1$, or $r\geq 2$ and for any two \emph{proper} faces
$F$ and $G$ of $P$, i.e.\ faces distinct from $F_{-1}$ and $F_r$,
there exists a finite sequence of consecutively incident proper
faces connecting $F$ with $G$, i.e.\ a sequence, starting with $F$
and terminating with $G$ such that every two consecutive faces in
the sequence are \emph{incident}, i.e.\ either in the relation
$\leq$ or in the converse relation (either $C_i\leq C_{i+1}$ or
$C_{i+1}\leq C_i$).

For $F$ and $G$ two faces of $P$ such that $F\leq G$, the
\emph{section} $G/F$ of $P$ is the set of faces $H$ of $P$ such
that $F\leq H\leq G$. Note that $P$ itself is a section of $P$,
and note that each section of $P$ is a partially ordered set of
rank $r'\leq r$.

We can then state the next property:
\begin{itemize}
\item[(P3)] $P$ is \emph{strongly connected}, which means that
every section of $P$ is a connected partially ordered set of rank
$r'\leq r$.
\end{itemize}

The \emph{rank $r'$ of a face} $F$ of $P$ is the rank of the
section $F/F_{-1}$, understood as a partially ordered set of rank
$r'\leq r$. We can then state the last defining property of
abstract polytopes:
\begin{itemize}
\item[(P4)] For every $i\in\{0,\ldots,r\mn 1\}$, if for $F$ and
$G$ of ranks $i\mn 1$ and $i\pl 1$ respectively we have that
$F\leq G$, then there are exactly two faces $H_1$ and $H_2$ of
rank $i$ such that $F\leq H_1\leq G$ and $F\leq H_2\leq G$.
\end{itemize}

As we announced above, we will not check properties (P3) and (P4)
for $\A(H)$ directly, but note that verifying (P4) when $i>0$ is
quite easy. This is because then $G\mn F=\{a_1,a_2\}$ with
$a_1\neq a_2$, and we have that $H_1=F\cup\{a_1\}$ and
$H_2=F\cup\{a_2\}$.

When sometimes, for the sake of brevity, we say in the remainder
of this section just polytope, we mean abstract polytope. We state
now the clauses and notions we need for our inductive definition
of an abstract polytope of rank $r\geq -1$:
\begin{itemize}
\item[(-1)] $P_{-1}=\langle \{F_{-1}\},\leq_{P_{-1}}\rangle$,
where $\leq_{P_{-1}}=\{(F_{-1},F_{-1})\}$, is a polytope of rank
$-1$, with the unique face $F_{-1}$ of $P_{-1}$ being of rank
$-1$. \item[(0)] For any object $a\neq F_{-1}$, we have that
$P_0^a=\langle\{F_{-1},a\},\leq_{P_0^a}\rangle$, where
$\leq_{P_0^a}= \{(F_{-1},F_{-1}),(F_{-1},a), (a,a)\}$, is a
polytope of rank 0, with the faces $F_{-1}$ and $a$ being
respectively of rank $-1$ and 0.
\end{itemize}
For a polytope of rank $r\geq 0$, its faces of rank $r\mn 1$ are
its \emph{facets}. (So $P_{-1}$ has no facets.)
\begin{itemize}
\item[(P$3'$)] Two distinct polytopes of the same rank are
\emph{close neighbours} when they have a common facet. A set of
polytopes all of the same rank is \emph{closely connected} when
each pair of distinct polytopes in it is connected by a finite
sequence of consecutively close neighbours (in other words, they
are connected by the transitive closure of the close neighbours
relation). \item[(P$4'$)] A set of polytopes is \emph{bivalent}
when every facet of a polytope in it belongs to exactly two
polytopes in this set.
\end{itemize}

Our inductive definition of an abstract polytope has clause (-1)
in the basis, and it has the following inductive clause:
\begin{itemize}
\item[] If $S$ is a closely connected bivalent set of polytopes of
rank $r$, let $P=\langle\bigcup S\cup\{F_{r+1}^P\},\leq_P\rangle$,
where
\\[1.5ex]
\mbox{\hspace{5em}}$\leq_P=\{(x,y)\mid (\exists Q\in S) x\leq_Q
y\; \mbox{ or }\; y=F_{r+1}^P\}$,
\\[1.5ex]
be a polytope of rank $r\pl 1$; the face $F_{r+1}^P$, which is not
a face of any polytope in $S$, is the unique face of $P$ of rank
$r\pl 1$, and the remaining faces have in $P$ the rank they had in
the polytopes of $S$.
\end{itemize}
This concludes the definition.

In $P$ the greatest faces of the polytopes of $S$ become facets,
and the facets of these polytopes become \emph{ridges}, i.e.\
faces of the rank of the polytope minus~2. Clause (0) is obtained
by applying our inductive clause to clause (-1) in a trivial
manner. (There are no two distinct polytopes in $S=\{P_{-1}\}$,
and there are no facets in $P_{-1}$.) We have however stated (0)
for the sake of clarity.

We have indeed obtained in this manner clause (0), because it is
possible that two different polytopes $P_1$ and $P_2$ of rank
$r\pl 1$ arise out of the same set $S$ of polytopes of rank $r$.
The faces $F_{r+1}^{P_1}$ and $F_{r+1}^{P_2}$ are then different
objects.

To prevent confusion, we should understand our inductive clause in
such a manner that if the sets $S_1$ and $S_2$ of polytopes of
rank $r$ are different, then the faces $F_{r+1}^{P_1}$ and
$F_{r+1}^{P_2}$ of the polytopes $P_1$ and $P_2$ of rank $r\pl 1$
with the carriers $S_1\cup\{F_{r+1}^{P_1}\}$ and $S_2\cup
\{F_{r+1}^{P_2}\}$ must be different.

All our abstract polytopes have the same least face $F_{-1}$, but
this is not an essential matter.

Our task now is to verify that our inductive definition is
equivalent with the definition in terms of (P1), $\ldots$~,(P4).
(An equivalent inductive definition, different from ours, is
mentioned in \cite{MS02}, Section 2A.)

It is very easy to see that every polytope of rank $r$ defined by
our inductive definition is a partial order of rank $r$, i.e.\
that it satisfies clauses (P1) and (P2). To show that it satisfies
(P3), we proceed by induction on~$r$.

If $r=-1$, then the clause is satisfied trivially. As a matter of
fact, it is satisfied trivially, by the definition of
connectedness, for every $r\leq 1$.

If $r>1$, then we have in $P$ old sections, which occur in a
polytope of the set $S$ used for defining $P$ inductively, which
are connected as before, and new sections $F_r^P/H$. For $F$ and
$G$ two proper faces in this new section we find the facets $F'$
and $G'$ of $P$ with which they are respectively incident. We find
that $F'$ and $G'$ are connected in $P$ by applying (P$3'$). (The
facets of $P$ are the greatest elements of the polytopes in $S$.)
So $F$ and $G$ are connected in~$P$.

To show that our inductively defined polytopes of rank $r$ satisfy
(P4) we proceed again by induction. The cases where $r\leq 0$ are
trivial (then $\{0,\ldots, r\mn 1\}$ is the empty set). For $r=1$
we have that all our polytopes have the following Hasse diagram:
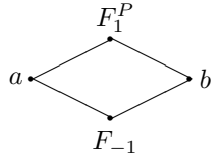
\begin{center}
\begin{picture}(60,50)
\put(30,10){\line(-2,1){30}} \put(30,10){\line(2,1){30}}
\put(30,40){\line(-2,-1){30}} \put(30,40){\line(2,-1){30}}

\put(30,10){\makebox(0,0){\circle*{2}}}
\put(0,25){\makebox(0,0){\circle*{2}}}
\put(60,25){\makebox(0,0){\circle*{2}}}
\put(30,40){\makebox(0,0){\circle*{2}}}

\put(31,5){\makebox(0,0)[t]{$F_{-1}$}}
\put(-4,25){\makebox(0,0)[r]{$a$}}
\put(63,25){\makebox(0,0)[l]{$b$}}
\put(31,45){\makebox(0,0)[b]{$F_1^P$}}
\end{picture}
\end{center}
(The facet of $P_0^a$ and $P_0^b$ is $F_{-1}$.)

For $r>1$, we have as in the paragraph above old sections $G/F$ of
$P$, which are taken care of by the induction hypothesis, and new
sections. The new sections are of the form $F_r^P/F$, and then for
them we have by (P$4'$) that exactly two facets $F_{r-1}^{Q_1}$
and $F_{r-1}^{Q_2}$ of $P$ are incident with a face $G$ of rank
$r\mn 2$ (which is an old facet).

With that we have finished verifying that the inductively defined
polytopes are abstract polytopes according to the old definition,
in terms of (P1), $\ldots$~,(P4). For the converse we have the
following.

Note first that the old definition allows for our polytope
$P_{-1}$ and for another polytope $P'_{-1}$ isomorphic to
$P_{-1}$, which is in all respects like $P_{-1}$ save that
$F_{-1}$ of $P_{-1}$ is replaced by an object $F'_{-1}$ different
from $F_{-1}$. With our clause (-1) we have provided only for a
single polytope $P_{-1}$. Two courses are now open to us. We may
either consider that $P_{-1}$ and $P'_{-1}$, since they are
isomorphic, are in fact the same, or we may consider that $F_{-1}$
in (-1) is in fact a variable. We will follow the first course,
but this is not an essential matter.

To verify \emph{close connectedness}, which corresponds to the
property defined in (P$3'$), for a polytope $P$ of rank $r\pl 1$
defined in the old way means verifying that for two distinct
facets (i.e.\ faces of rank $r$) $F'$ and $F''$ of $P$ there is a
finite sequence of \emph{close neighbours} connecting $F'$ with
$F''$; two close neighbours being now two facets incident with a
common ridge of $P$, i.e.\ face of rank $r\mn 1$. We will show
that close connectedness is a consequence of essentially (P3).

Let a \emph{proper path} of $P$ from $F'$ to $F''$ be, as in the
definition of connectedness above, a finite sequence of
consecutively incident proper faces of $P$ connecting $F'$ with
$F''$.

For $k_{l-2}$ being the number of faces of rank $r\mn l$, for
$l\geq 2$, in a proper path $\Pi$ let the \emph{weight} of $\Pi$
be the ordinal
\[k_0\omega^0\pl\ldots\pl k_n\omega^n,\;\;\mbox{\rm for }\; n=r\mn 2.\]
If $k_0=\ldots=k_n=0$, then the weight of $\Pi$ is 0, and in $\Pi$
we have only facets and ridges, as required by close
connectedness. The weight is an ordinal less than $\omega^\omega$,
which is infinite if one of $k_1,\ldots, k_n$ is greater than 0.
We can prove the following.

\prop{Lemma 8.1}{For every proper path $\Pi$ from $F'$ to $F''$ of
weight $w$ greater than $0$, there is a proper path $\Pi'$ from
$F'$ to $F''$ of weight strictly less than $w$.}

\dkz Let $\Pi$ be the sequence $F_1\ldots F_m$, with $F_1$ being
$F'$ and $F_m$ being $F''$; since $w>0$, we must have that $m\geq
5$. Let $\rho(F_i)$ be the rank of $F_i$ in $P$, and consider a
subsequence $F_i\ldots F_{i+2}$ of $\Pi$ such that
$\rho(F_i)=\rho(F_{i+1})\pl 1=\rho(F_{i+2})\leq r\mn 1$. Such a
subsequence must exist because $w>0$. For $G$ being the greatest
element of $P$ consider the section $G/F_{i+1}$. Then we have that
either $F_i=F_{i+2}$, or for $F_i\neq F_{i+2}$ by (P3) there is a
proper path $F_i G_1\ldots G_q F_{i+2}$, for $q\geq 1$, of the
polytope $G/F_{i+1}$.

For $\Gamma$ being either empty or $G_1\ldots G_q F_{i+2}$, let
$\Pi'$ be $F_1\ldots F_i\Gamma F_{i+3}\ldots F_m$. The weight of
$\Pi'$ is less than $w$. \qed

\noindent From this lemma we can infer close connectedness by
induction.

Suppose that for $P$ of rank $r$ a face $F$ of rank $r\mn 2$ does not
respect bivalence, which means that it is not incident with
exactly two faces of $P$ of rank $r\mn 1$, i.e.\ two facets of
$P$. It is clear that the section $F_r^P/F$ does not then respect
(P4).

With that we have finished verifying that the old definition and
the new inductive definition are equivalent. We will verify below
that $\A(H)$ is an abstract polytope by relying on the new
inductive definition. For that verification we need some
preliminary matters.

For $P_1$ and $P_2$ abstract polytopes, consider the partial order
$P_1\cdot P_2$ with carrier
$((P_1\mn\{F_{-1}\})\times(P_2\mn\{F_{-1}\}))\cup\{F_{-1}\}$,
where $(F_1,F_2)\leq_{P\cdot Q}(G_1,G_2)$ iff $F_1\leq_{P_1} G_1$
and $F_2\leq_{P_2} G_2$, and moreover $F_{-1}$ is the least
element with respect to $\leq_{P_1\cdot P_2}$. Note that
$\A(H_1)\otimes\A(H_2)$, $\A(H_1)\ast\A(H_2)$ and
$\A(H_1)\cdot\A(H_2)$ are all isomorphic, provided the first two
products are defined (for $\otimes$ we must have $\bigcup H_1$ and
$\bigcup H_2$ disjoint, and for $\ast$, the hypergraphs $H_1$ and
$H_2$ should be $H_Y$ and $\uhy H$ for $Y\in H$).

Then we can prove the following.

\prop{Proposition 8.2}{For $P_1$ and $P_2$ abstract polytopes of
ranks $r_1$ and $r_2$ respectively, $P_1\cdot P_2$ is an abstract
polytope of rank $r_1\pl r_2$.}

\dkz With the old definition of an abstract polytope, it is easy
to check (P1) and (P2). For (P3) we make our connection via
$\leq_{P\cdot Q}$ by moving first in one coordinate and then in
the other. To check (P4) suppose that, for $\rho(F)$ being the
rank of $F$ in the appropriate partial order, we have that
$(F_1,F_2)\leq_{P\cdot Q}(G_1,G_2)$ and
$\rho(G_1,G_2)\mn\rho(F_1,F_2)=2$. We have the last equation only
if one of the following cases obtains:
\begin{tabbing}
\hspace{1.5em}\=(1)\quad\=$\rho(G_1)\mn\rho(F_1)=2$,\quad\=
$\rho(G_2)=\rho(F_2)$,
\\
\>(2)\>$\rho(G_1)=\rho(F_1)$,\> $\rho(G_2)\mn\rho(F_2)=2$,
\\
\>(3)\>$\rho(G_1)\mn\rho(F_1)=1$,\> $\rho(G_2)\mn\rho(F_2)=1$.
\end{tabbing}

In cases (1) and (2) we just rely on (P4) for $P_1$ and $P_2$
respectively. In case (3) we have $(F_1,G_2)$ for $H_1$ and
$(F_2,G_1)$ for $H_2$. That the rank of $P_1\cdot P_2$ is the sum
of the ranks $r_1$ and $r_2$ of $P_1$ and $P_2$ is clear from the
fact that we go up now in two coordinates, $r_1$ steps in the
first coordinate, and $r_2$ steps in the second. \qed

We can then prove the following.

\prop{Proposition 8.3}{For every atomic hypergraph $H$ we have
that $\A(H)$ is an abstract polytope of rank $|\bigcup H|\mn n$,
for $n$ the connectedness number of $H$.}

\dkz We rely on the inductive definition of $\A(H)$ introduced at
the end of Section~7. In the basis of our induction, it is clear
that $\{\{x\},\{x,*\}\}$ is an abstract polytope of rank 0.

Suppose for the induction step that, for an ASC-hypergraph $H$ and
for $Y\in H\mn\{\bigcup H\}$, we have that $\A(H_Y)$ and $\A(\uhy
H)$ are abstract polytopes of ranks $|Y|\mn 1$ and $|\bigcup H\mn
Y|\mn1$ respectively. To have both of these ranks at least 0, we
must have that $|\bigcup H|\geq 2$. Then we have that
$\A(H_Y)\ast\A(\uhy H)$, which is isomorphic to
$\A(H_Y)\cdot\A(\uhy H)$, is an abstract polytope of rank
$|\bigcup H|\mn 2$, by Proposition 8.2.

Let $S$ be the set $\{\A(H_Y)\ast\A(\uhy H)\mid Y\in H\mn\{\bigcup
H\}\}$. We show that this set is closely connected. Take the
polytopes $P_Y=\A(H_Y)\ast\A(\uhy H)$ and
$P_Z=\A(H_Z)\ast\A(_{\raisebox{-.5pt}{$\cup$} H-Z} H)$ from $S$.
Suppose that $Y\cup Z=\bigcup H$.

If $|\bigcup H|=2$, then $Y=\{y\}$ and $Z=\{z\}$. In that case
$P_Y$ and $P_Z$ are the only polytopes in $S$. They are closely
connected because they have $F_{-1}$ as a common facet, and $S$ is
bivalent.

If $|\bigcup H|>2$, then either $|Y|\geq 2$ or $|Z|\geq 2$.
Suppose $|Y|\geq 2$ (in case $|Z|\geq 2$ we proceed analogously).
Then $Y\not\subseteq Z$, since $Y$ and $Z$ are both proper subsets
of $\bigcup H$, and $Y\cup Z=\bigcup H$. So for some $y\in Y$ we
have that $y\not\in Z$. We have that $\{y\}\in H$, since $H$ is
atomic, and $y\in\bigcup H$, and so
$P_{\{y\}}=\A(H_{\{y\}})\ast\A(_{\raisebox{-.5pt}{$\cup$} H-\{y\}}
H)$ is a polytope in $S$. The polytopes $P_Y$ and $P_{\{y\}}$ have
as a common facet $\{\{y\},Y,\bigcup H\}$; the polytopes
$P_{\{y\}}$ and $P_Z$ have as a common facet $\{\{y\},Z,\bigcup
H\}$ if $Z\cup\{y\}\not\in H$, and they have as a common facet
$\{Z,Z\cup\{y\},\bigcup H\}$ if $Z\cup\{y\}\in H$. From that we
may conclude that $S$ is closely connected and bivalent.

Suppose that $Y\cup Z\subset \bigcup H$; in that case $|\bigcup
H|>2$. The following cases may arise.

\vspace{1ex} (1) If $Y\subseteq Z$ or $Z\subseteq Y$, then $P_Y$
and $P_Z$ have as a common facet $\{Y,Z,\bigcup H\}$.

\vspace{1ex} (2) Suppose neither $Y\subseteq Z$ nor $Z\subseteq
Y$.

\vspace{1ex} (2.1) If $Y\cup Z\in H$, then since $Y\cup Z\subset
\bigcup H$, we have $P_{Y\cup Z}=\A(H_{Y\cup
Z})\ast\A(_{\raisebox{-.5pt}{$\cup$} H-(Y\cup Z)} H)$ as another
polytope in $S$. The polytopes $P_Y$ and $P_{Y\cup Z}$ have as a
common facet $\{Y,Y\cup Z,\bigcup H\}$, and $P_{Y\cup Z}$ and
$P_Z$ have as a common facet $\{Z,Y\cup Z,\bigcup H\}$.

\vspace{1ex} (2.2) If $Y\cup Z\not\in H$ and $Y\cap Z=\emptyset$,
then, as for (1), the facet $\{Y,Z,\bigcup H\}$ is common to $P_Y$
and $P_Z$.

\vspace{1ex}\noindent Note that because of the saturation of $H$
it is impossible that $Y\cup Z\not\in H$ and $Y\cap
Z\neq\emptyset$. From all that we conclude that $S$ is closely
connected. That $S$ is bivalent follows from the fact that, if
$|\bigcup H|>2$, every facet of a polytope in $S$ is of the form
$\{Y,Z,\bigcup H\}$ for two polytopes $P_Y$ and $P_Z$ as above.

We add to $S$ the new face $\{\bigcup H\}$ and obtain the polytope
$\A(H)$ of rank $|\bigcup H|\mn 1$, in accordance with the
inductive clause of the new definition of an abstract polytope.

We have said at the end of Section~7 that for every atomic
hypergraph $H$ we can obtain $\A(H)$ as
$\A(H_1)\otimes\ldots\otimes\A(H_n)$, for $n\geq 1$, where
$\{H_1,\ldots,H_n\}$ is the finest hypergraph partition of the
saturated closure of $H$. The hypergraphs $H_1,\ldots,H_n$ are
ASC-hypergraphs, and so $\A(H_1),\ldots,\A(H_n)$ are, according to
the proof above, abstract polytopes of ranks $|\bigcup H_1|\mn
1,\ldots,|\bigcup H_n|\mn 1$ respectively. That
$\A(H_1)\otimes\ldots\otimes\A(H_n)$ is an abstract polytope of
rank $|\bigcup H|\mn n$ follows then by Proposition 8.2, and by
the isomorphism of the products $\cdot$
and~$\otimes$.\mbox{\hspace{.5em}} \qed

The proof we have just given shows that the facets of the abstract
polytope $\A(H)$ for $H$ atomic and connected are never very far
from each other. They either share a common ridge, or there is in
between them a facet with which each of them shares a ridge.

\section{\large \textbf{Realizations}}
In this section we turn towards a geometrical approach to
hypergraph polytopes. Our main concern in this work was the
abstract approach of the preceding sections, and so in this
section, which is in a related, but nevertheless different field,
our exposition will be less detailed. We will strive to be concise
in order not to make an already long text still longer. We will
give no examples in this section; they can be worked out from
Appendix~B. We presume the reader is acquainted with some basic
notions of the geometrical theory of polytopes, which may all be
found in \cite{Z95}, whose terminology we will follow.

We define for every atomic hypergraph $H$ a convex polytope in
$\R^n$, for whose face lattice we will prove that it is isomorphic
(as a partial order) to $\A(H)$. Most of the section is devoted to
proving that for ASC-hypergraphs. The proof for the remaining
atomic hypergraphs will then follow easily.

Let $H$ be an ASC-hypergraph on the carrier $\bigcup
H=\{1,\ldots,d\pl 1\}$ for $d\geq 0$. We take now the elements of
$\bigcup H$ to be positive integers, because we want them to
function as indices, but we may, however, take $\bigcup H$ to be
an arbitrary finite nonempty set. We deal separately below with
the case of the empty hypergraph $\emptyset$ on the carrier
$\emptyset$.

For $\cal S$ being the set of hyperplanes in $\R^{d+1}$, let the
map $\pi\!:H\rightarrow{\cal S}$ be defined by
\[\pi_X=\{(x_1,\ldots,x_{d+1})\in\R^{d+1}\mid
\mbox{$\sum_{i\in X}$} x_i=3^{|X|}\},\] where $\pi_X$ stands for
$\pi(X)$. Note that $\pi$ is one-one. The function $f(x)=3^x$ in
the definition of $\pi_X$ is not the only one that could be
chosen. Any function on natural numbers that would enable us to
prove an analogue of Lemma 9.1 below would do. The function
$f(x)=3^x$ is one of the ``suitable'' functions introduced in
\cite{S97a} (Appendix~B). Intuitively, the choice of $3^x$ may be
explained by the wish not to truncate too much. In a very simple
case, this means that after truncating a one-dimensional edge at
both ends, we are left with something; hence we divide the edge in
three parts.

Consider the closed halfspace
\[\pi_X^+=\{(x_1,\ldots,x_{d+1})\in\R^{d+1}\mid \mbox{$\sum_{i\in X}$}
x_i\geq 3^{|X|}\},\] whose boundary hyperplane is $\pi_X$. Let
$\G(H)$ be the polytope
\[\mbox{$\bigcap\{\pi_X^+\mid X\in H\mn\{\bigcup H\}\}\cap \pi_{\ccup H}$}\]
in the hyperplane $\pi_{\ccup H}$.

That $\G(H)$ is indeed an $\cal H$-polytope, and not just an $\cal
H$-polyhedron (see \cite{Z95}, Lecture~0, Definition 0.1), is
guaranteed by the atomicity of $H$. The set $\G(H)$ is bounded by
the $d$-dimensional simplex $\bigcap\{\pi_{\{i\}}^+\mid i\in
\bigcup H\}\cap \pi_{\ccup H}$. Intuitively, we may assume that
the polytope $\G(H)$ is obtained by truncating this simplex, which
is a limit case with no truncation; in that case, $H$ is just
$\{\{i\}\mid i\in\bigcup H\}\cup\{\bigcup H\}$. The limit case at
the other end, with all possible truncations, is with $H$ being
the set of all nonempty subsets of $\bigcup H$; in that case we
obtain the $d$-dimensional permutohedron.

In the case of the empty hypergraph $\emptyset$ we define
$\G(\emptyset)$ to be the polytope $\R^0=\{*\}$ in $\R^0$, whose
face lattice is $\{\{*\},\emptyset\}$, with $\leq$ being
$\subseteq$. This face lattice is isomorphic to
$\A(\emptyset)=\{\emptyset,\{*\}\}$ (see the end of Section~5) by
the bijection that assigns the vertex $\emptyset$ of
$\A(\emptyset)$ to the vertex $\{*\}$ of $\G(\emptyset)$, and
$\{*\}$, which is the face $F_{-1}$ in $\A(\emptyset)$, to
$\emptyset$, which is the face $F_{-1}$ for $\G(\emptyset)$; the
partial order $\leq$ in $\A(\emptyset)$ is the converse of
$\subseteq$.

To make $\G(\emptyset)$ a limit case of the definition given
above, take first $\{1,\ldots,d\pl 1\}$ to be the empty set when
$d=-1$. In that case, we have $\R^{d+1}=\R^0=\{*\}$, and
$\bigcap\{\pi_X^+\mid X\in\emptyset\}=\{*\}$. The face $F_{-1}$ of
any abstract polytope should be mapped to the empty subset of
$\R^n$ (the empty set is a face of any geometrical polytope; see
\cite{Z95}, Lectures 0 and~2, and \cite{MS02}, Section 5a). So
$\G(\emptyset)$ has one vertex $\{*\}$, and the face $\emptyset$,
as the image of $F_{-1}$. We may extend $\pi$ and $\pi^+$ to
$\bar{H}^*=H\cup\{*\}$ so that $\pi_*=\pi_*^+=\emptyset$; then
$\bigcap\{\pi_X^+\mid X\in \bar{H}^*\mn\{\bigcup H\}\}=\emptyset$.

We have however separated the case $\G(\emptyset)$ from the rest
because it is degenerate; as $\A(\emptyset)$, the polytope
$\G(\emptyset)$ has no important role to play. (Note that in the
inductive definition of $\A(H)$, at the end of Section~7, we do
not have $\A(H)$ in the basis.)

We will prove that $\G(H)$ is a polytope each of whose vertices
lies in exactly $d\pl 1$ boundary hyperplanes; these $d\pl 1$
hyperplanes are the elements of $\{\pi_X\mid X\in K\}$ for some
construction $K$ of $H$ (for a precise statement see Proposition
9.3). This will imply that $\G(H)$ is a simple polytope (which
means that each of its vertex figures, which are figures obtained
by truncating a vertex, is a simplex; see \cite{Z95}, Section 2.5,
Proposition 2.16). On the other hand, we will prove that for every
construction $K$ of $H$ the intersection of the hyperplanes in
$\{\pi_X\mid X\in K\}$ is a vertex of $\G(H)$ (see Proposition
9.6). From all that we will conclude that $\G(H)$ has a face
lattice isomorphic to $\A(H)$.

The set ${\cal P}^+(\bigcup H)$ of all the nonempty subsets of
$\bigcup H$ is an ASC-hypergraph on $\bigcup H$. The set ${\cal
P}^+(\bigcup H)$ is a subset of itself, and we can consider all
the ${\cal P}^+(\bigcup H)$-antichains, in accordance with our
definition of an $M$-antichain in Section~6. We prove the
following.

\prop{Lemma 9.1}{Let $\{X_1,\ldots,X_n\}$ be a ${\cal P}^+(\bigcup
H)$-antichain. If for every $i\in X_1\cup\ldots\cup X_n$ we have
that $x_i\geq 0$, and for every $j\in\{1,\ldots,n\}$ we have that
$\sum_{i\in X_j}x_i\leq 3^{|X_j|}$, then $\sum_{i\in
X_1\cup\ldots\cup X_n}x_i< 3^{|X_1\cup\ldots\cup X_n|}$.}

\dkz We proceed by induction on $n\geq 2$. If $n=2$, then
\begin{tabbing}
\hspace{1.5em}$\sum_{i\in X_1\cup X_2}x_i$ \= $\leq \sum_{i\in
X_1}x_i + \sum_{i\in X_2}x_i$,\quad since $x_i\geq 0$
\\*[.5ex]
\> $\leq 3^{|X_1|}+ 3^{|X_2|}$,\quad by the assumption
\\
\> $< 3^{\,\mbox{\scriptsize max}(|X_1|,|X_2|)+1}$
\\
\> $\leq 3^{|X_1\cup X_2|}$,\quad since $X_1,X_2\subset X_1\cup
X_2$.
\end{tabbing}

If $n>2$, then let $X=X_1\cup\ldots\cup X_{n-1}$. By the induction
hypothesis, we have that $\sum_{i\in X}x_i<3^{|X|}$. If
$X_n\subseteq X$, i.e.\ $X_1\cup\ldots\cup X_n=X$, then we are
done. If $X_n\not\subseteq X$, then, since $\{X_1,\ldots, X_n\}$
is a ${\cal P}^+(\bigcup H)$-antichain, we cannot have that
$X\subseteq X_n$; hence $\{X,X_n\}$ is a ${\cal P}^+(\bigcup
H)$-antichain, and we may apply the induction hypothesis to it.
\qed

\noindent The main arithmetical idea of this lemma is based on the
following. If for every $i\in\{1,\ldots,k\}$ we have that $m_i<m$,
then we have the inequality
\begin{tabbing}
\centerline{$\sum_{i=1}^k(k\pl 1)^{m_i}<(k\pl 1)^m$,}
\end{tabbing}
which yields as a particular case:
\begin{tabbing}
\centerline{if $m_1<m$ and $m_2<m$, then $3^{m_1}\pl
3^{m_2}<3^m$.}
\end{tabbing}
\noindent The idea of this inequality may be gathered from
\cite{G38} (see also \cite{DM79}).

\prop{Lemma 9.2}{For every $M\subseteq H$, if $\bigcap\{\pi_X\mid
X\in M\}\cap \G(H)\neq\emptyset$, then every $M$-antichain misses
$H$.}

\dkz Suppose $\{X_1,\ldots,X_n\}$ is an $M$-antichain such that
$X_1\cup\ldots\cup X_n\in H$, and let
\begin{tabbing}
\hspace{7em}$(x_1,\ldots,x_{d+1})\in\bigcap\{\pi_X\mid X\in
M\}\cap \G(H)$.
\end{tabbing}
Since $H$ is atomic, we have that
\begin{tabbing}
\hspace{7em}$(x_1,\ldots,x_{d+1})\in \G(H)\subseteq
\bigcap\{\pi_{\{i\}}^+\mid i\in \{1,\ldots,d\pl 1\}\}$,
\end{tabbing}
and hence $x_1,\ldots,x_{d+1}\geq 0$. Since
\begin{tabbing}
\hspace{7em}$(x_1,\ldots,x_{d+1})\in\bigcap\{\pi_{X_j}\mid j\in
\{1,\ldots,n\}\}$,
\end{tabbing}
we have for every $j\in\{1,\ldots,n\}$ that $\sum_{i\in
X_j}x_i=3^{|X_j|}\leq 3^{|X_j|}$. Then by Lemma 9.1 we have that
$\sum_{i\in X_1\cup\ldots\cup X_n}x_i<3^{|X_1\cup\ldots\cup
X_n|}$, and hence
$(x_1,\ldots,x_{d+1})\not\in\pi_{X_1\cup\ldots\cup X_n}^+$ (we
have that $\pi_{X_1\cup\ldots\cup X_n}^+$ is defined since
$X_1\cup\ldots\cup X_n\in H$). From
$\G(H)\subseteq\pi_{X_1\cup\ldots\cup X_n}^+$ we conclude that
$(x_1,\ldots,x_{d+1})\not\in \G(H)$, which is a contradiction.
\qed

The following proposition shows that $\G(H)$ is a polytope whose
vertices correspond to constructions of $H$.

\prop{Proposition 9.3}{For every vertex $\{v\}$ of $\G(H)$ there
is a construction $K$ of $H$ such that $\{v\}=\bigcap\{\pi_X\mid
X\in K\}$ and for every $Y$ in $H\mn K$ we have that
$v\not\in\pi_Y$.}

\dkz Since $\{v\}$ is a vertex of $\G(H)$, there are at least $d$
members $X_1,\ldots X_d$ of $H\mn\{\bigcup H\}$ such that
\[\pi_{X_1}\cap\ldots\cap\pi_{X_d}\cap\pi_{\ccup H}=\{v\}\subseteq\G(H).\]
Let $K=\{X_1,\ldots,X_d,\bigcup H\}$. By Lemma 9.2, every
$K$-antichain misses $H$, and so, by Proposition 6.11, we have
that $K$ is a construction of $H$. It remains to establish that if
$Y\in H\mn K$, then $v\not\in\pi_Y$.

Let $Y\in H\mn K$. By Lemma 6.9, for $M=K\cup\{Y\}$ there is an
$M$-antichain that does not miss $H$. By Lemma 9.2, we have that
$\bigcap\{\pi_X\mid X\in M\}\cap\G(H)=\emptyset$, and hence
$v\not\in\pi_Y$ since $v\in\bigcap\{\pi_X\mid X\in K\}\cap\G(H)$.
\qed

From this proposition it is easy to derive the following two
corollaries.

\propo{Corollary 9.31}{For every vertex $\{v\}$ of $\G(H)$ there
is a unique construction $K$ of $H$ such that
$\{v\}=\bigcap\{\pi_X\mid X\in K\}$.}

\prop{Corollary 9.32}{For every vertex $\{v\}$ of $\G(H)$ there
are exactly $d$ halfspaces $\pi_X^+$ such that $X\in H\mn\{\bigcup
H\}$ and $v\in\pi_X$.}

\noindent For Corollary 9.32 we rely on the fact that every
construction of $H$ has $d\pl 1$ members, one of which is $\bigcup
H$, and that $\pi_{X_1}=\pi_{X_2}$ implies
$\pi_{X_1}^+=\pi_{X_2}^+$.

Now we need to show in the converse direction that for every
construction $K$ of $H$ there is a vertex of $\G(H)$ such that
$\{v\}=\bigcap\{\pi_X\mid X\in K\}$. This will be a consequence of
the following two lemmata.

\prop{Lemma 9.4}{For every construction $K$ of $H$ there is a
unique function $x\!:\bigcup H\rightarrow \R$ such that all the
equations in the set $\{\sum_{i\in X}x(i)=3^{|X|}\mid X\in K\}$
hold. Moreover, if for $X\in K$ we have that $s$ is
$X$-superficial, then $x(s)>3^{|X|-1}$.}

\dkz By induction on $d$.

If $d=0$, then $H=\{\{1\}\}$, and the only construction of $H$ is
$H$ itself. Then we have a single equation $x(1)=3$ that defines
the function $x$, and $x(1)=3>3^{1-1}=1$.

Suppose $d>0$. We have that there is an $s\in \bigcup H$ such that
$H_1,\ldots,H_n$, for $n\geq 1$, is the finest hypergraph
partition of $H_{\ccup H-\{s\}}$, and $K=K_1\cup\ldots\cup
K_n\cup\{\bigcup H\}$, for $K_j$ being a construction of the
ASC-hypergraph $H_j$. By the induction hypothesis, for every
$j\in\{1,\ldots,n\}$ there is a unique function $x^j\!:\bigcup
H_j\rightarrow\R$ such that all the equations in the following set
$\{\sum_{i\in X}x^j(i)=3^{|X|}\mid X\in K_j\}$ hold. Since
$\bigcup H_j\in K_j$, we must have that
\begin{tabbing}
\hspace{1.5em}(a)\hspace{1.5em}$\sum_{i\in \ccup
H_j}x^j(i)=3^{|\ccup H_j|}$.
\end{tabbing}
To obtain $x\!:\bigcup H\rightarrow\R$ with the desired
properties, we form the union of the functions $x^j$ (which is a
function since the sets $\bigcup H_j$ are disjoint), and it
remains to find the unique value of $x(s)$. (Note that $x(s)$
figures only in the equation $\sum_{i\in \ccup H}x(i)=3^{|\ccup
H|}$.) So, for $i\in\bigcup H_j$, let $x(i)=x^j(i)$, and let
\begin{tabbing}
\centerline{$x(s)=3^{|\ccup H|}-\sum_{i\in\ccup H-\{s\}}x(i)$.}
\end{tabbing}

We have by the induction hypothesis that if $s_j$ is the
$X$-superficial element for $X\in K_j$, then
$x(s_j)=x^j(s_j)>3^{|X|-1}$. In particular, we have for every
$i\in\bigcup H\mn\{s\}$ that $x(i)\geq 0$. It remains to check the
analogous inequality concerning the $\bigcup H$-superficial
element $s$; namely $x(s)>3^{|\ccup H|-1}$. In the case $n=1$, we
have by (a) that
\begin{tabbing}
\centerline{$\sum_{i\in\ccup H-\{s\}}x(i)=\sum_{i\in\ccup
H_1}x(i)=3^{|\ccup H_1|}=3^{|\ccup H|-1}$,}
\end{tabbing}
and hence $x(s)=3^{|\ccup H|}-3^{|\ccup H|-1}>3^{|\ccup H|-1}$.

If $n>1$, then we apply Lemma 9.1 to the $\P^+(\bigcup
H)$-antichain $\{\bigcup H_1,\ldots,$ $\bigcup H_n\}$, relying on
(a) and on the fact that $x(i)\geq 0$ for every $i\in\bigcup
H\mn\{s\}$, in order to obtain
\begin{tabbing}
\centerline{$\sum_{i\in\ccup H-\{s\}}x(i)=\sum_{i\in(\ccup
H_1)\cup\ldots\cup(\ccup H_n)}x(i)<3^{|(\ccup
H_1)\cup\ldots\cup(\ccup H_n)}|=3^{|\ccup H|-1}$.}
\end{tabbing}
With this we have that
\begin{tabbing}
\centerline{$x(s)=3^{|\ccup H|}\mn \sum_{i\in\ccup
H-\{s\}}x(i)>3^{|\ccup H|}\mn 3^{|\ccup H|-1}>3^{|\ccup
H|-1}$.}\hspace{-.7em}$\dashv$
\end{tabbing}

This lemma says that the system of equations $\{\sum_{i\in
X}x(i)=3^{|X|}\mid X\in K\}$ has a unique solution, and this
solution is the unique element of $\bigcap\{\pi_X\mid X\in K\}$.
Moreover, it says something about the location of the coordinates
of this solution, and hence about the location of this solution in
the interior of $\bigcap\{\pi_Y^+\mid Y\in H\mn K\}$, which will
serve to determine that the solution is a vertex of $\G(H)$.

For $K$ a construction of $H$ and $x\!:\bigcup H\rightarrow \R$
the function obtained by Lemma 9.4, we can prove the following.

\prop{Lemma 9.5}{If $v=(x(1),\ldots,x(d\pl 1))\in\R^{d+1}$, then
$\{v\}$ is a vertex of $\G(H)$.}

\dkz By Lemma 9.4, we have that $\{v\}=\bigcap\{\pi_X\mid X\in
K\}$. So for the rest of the proof it is sufficient to show that
for every $Y\in H\mn K$ we have that $v\in\pi_Y^+$.

Let $Y\in H\mn K$. By Lemma 6.14, there is an $X\in K$ such that
$Y\subset X$, and the $X$-superficial element $s$ is in $Y$. By
Lemma 9.4, we have that $x(s)>3^{|X|-1}\geq 3^{|Y|}$, and since
every other $x(i)\geq 0$ (see the proof of Lemma 9.4) we have that
$\sum_{i\in Y}x(i)\geq x(s)>3^{|Y|}$, and hence $v\in
\pi_Y^+\mn\pi_Y$. \qed

As a corollary of Lemmata 9.4 and 9.5, we have the following.

\prop{Proposition 9.6}{For every construction $K$ of $H$ there is
a unique vertex $\{v\}$ of $\G(H)$ such that
$\{v\}=\bigcap\{\pi_X\mid X\in K\}$.}

For the following two lemmata, which are not about $\G(H)$
specifically, but are more general, we have the following
assumptions. Let $d\geq 1$ and let $\pi$ be a hyperplane in
$\R^{d+1}$. Let $S^+$ be a set of closed halfspaces in $\R^{d+1}$
whose boundary hyperplanes are collected in a set $S$. Let
$P=(\bigcap S^+)\cap\pi\neq\emptyset$ be a polytope. For every
vertex $\{v\}$ of $P$ let there be exactly $d$ halfspaces from
$S^+$ such that $\{v\}$ is contained in their boundary
hyperplanes. Then we can prove the following (for the notion of
simple polytope see \cite{Z95}, Section 2.5, Proposition 2.16).

\prop{Lemma 9.7}{$P$ is a simple $d$-dimensional polytope.}

\dkz Since $P\neq\emptyset$, the set of vertices of $P$ is not
empty. Let $\{v\}$ be a vertex of $P$. We show that a vertex
figure (see \cite{Z95}, Section 2.1) of $P$ at $v$ is a $(d\mn
1)$-dimensional simplex, from which the proposition follows.

Since $\{v\}$ is a vertex, there is a set $V^+$ of $d$ halfspaces
from $S^+$, whose boundary hyperplanes are collected in a set $V$,
such that $(\bigcap V)\cap\pi=\{v\}$. By our assumption, for every
halfspace in $S^+\mn V^+$ (we have that $S^+\mn V^+\neq\emptyset$
since $d\geq 1$ and $P$ is a polytope, and not just an $\cal
H$-polyhedron) we have that $v$ lies in the interior of this
halfspace.

Since $S^+\mn V^+$ is finite, we have that there is an open
neighbourhood $U$ of $v$ in $\R^{d+1}$ such that $(\bigcap
V^+)\cap\pi\cap U\subseteq P$. Hence a vertex figure of $P$ at $v$
is a $(d\mn 1)$-dimensional simplex. The following picture
illustrates the case when $d=3$:
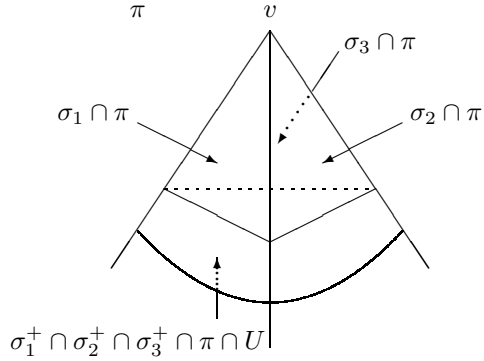
\begin{center}
\begin{picture}(60,130)(0,5)
\put(40,120){\line(-2,-3){60}} \put(40,120){\line(2,-3){60}}
\put(40,120){\line(0,-1){120}}

\multiput(0,60)(4,0){20}{\line(1,0){1}}
\put(40,40){\line(-2,1){40}} \put(40,40){\line(2,1){40}}

\qbezier(-10,44)(40,-10)(90,44)

\put(40,125){\makebox(0,0)[b]{$v$}}

\put(90,85){\vector(-2,-1){30}} \put(-10,85){\vector(2,-1){30}}
\put(56,97){\line(2,3){10}}
\multiput(45,80)(2,3){6}{\makebox(0,0){\circle*{.5}}}
\put(45,80){\vector(-2,-3){2}}

\put(20,11){\line(0,1){10}}
\multiput(20.5,22)(0,2){6}{\makebox(0,0){\circle*{.5}}}
\put(20,32){\vector(0,1){2}}

\put(93,88){\makebox(0,0)[l]{$\sigma_2\cap\pi$}}
\put(-13,88){\makebox(0,0)[r]{$\sigma_1\cap\pi$}}
\put(68,113){\makebox(0,0)[bl]{$\sigma_3\cap\pi$}}
\put(-10,8){\makebox(0,0)[t]{$\sigma_1^+\cap\sigma_2^+\cap\sigma_3^+\cap\pi\cap
U$}} \put(-10,125){\makebox(0,0)[b]{$\pi$}}
\end{picture}
\end{center}
\mbox{\hspace{.5em}}\qed

\oprop{Lemma 9.8}{There is a bijection $\beta$ from the set of all
the facets of $P$ to the set of all the pairs $\{\sigma,\pi\}$,
where $\sigma$ is an element of $S$ that contains a vertex of $P$.
This bijection is such that for every facet $\varphi$ and every
vertex $\{v\}$ of $P$ we have that $\{v\}\subseteq\varphi$ iff
$\{v\}\subseteq\bigcap\beta(\varphi)$.}

\dkz Let $\beta$ be a relation between the set of all the facets
of $P$ and the set of all the pairs $\{\sigma,\pi\}$, where
$\sigma$ is an element of $S$ that contains a vertex of $P$,
defined by
\[
(\varphi,\{\sigma,\pi\})\in\beta\; \mbox{when the affine hull of
$\varphi$ is $\sigma\cap\pi$}.
\]
First we show that $\beta$ is a function between these two sets.
Let $\varphi$ be a facet of $P$. Since $d\geq 1$, there is a
vertex $\{v\}$ of $P$ incident with $\varphi$. By reasoning as in
the proof of Lemma 9.7, there must be a $\sigma\in S$ such that
for $\mbox{aff}(\varphi)$ being the affine hull of $\varphi$ we
have that $\mbox{aff}(\varphi)=\sigma\cap\pi$; so $v\in\sigma$,
and we can take $\beta(\varphi)=\{\sigma,\pi\}$.

If $\mbox{aff}(\varphi)=\sigma'\cap\pi$, then $v\in\sigma'$, and,
since there are exactly $d$ hyperplanes from $S$ containing
$\{v\}$, it must be that $\sigma'=\sigma$; otherwise, we would
have that
$\sigma'\cap\sigma\cap\ldots\cap\pi=\{v\}=\sigma\cap\ldots\cap\pi$,
and hence we would have $d$ hyperplanes in $\R^{d+1}$ intersecting
in a point, which is impossible. So $\beta$ is a function.

Since for all facets $\varphi_1$ and $\varphi_2$ of $P$ we have
that $\mbox{aff}(\varphi_1)= \mbox{aff}(\varphi_2)$ iff
$\varphi_1=\varphi_2$, we have that $\beta$ is one-one. It remains
to show that $\beta$ is onto.

Let $\sigma\in S$, and let $\{v\}$ be a vertex of $P$ such that
$v\in\sigma$. By reasoning as in the proof of Lemma 9.7, we have
that $\sigma\cap P$ is a facet of $P$ whose affine hull is
$\sigma\cap\pi$. So this is the facet mapped by $\beta$ to
$\{\sigma,\pi\}$.

The equivalence of the proposition from left to right is trivial.
For the other direction, we rely on the fact that for every facet
$\varphi$ of $P$ we have that $\mbox{aff}(\varphi)\cap P=\varphi$.
\qed

For $d\geq 1$, by Corollary 9.32 and by Proposition 9.6 (which
gives that $\G(H)$ is not empty), we have that $\G(H)$ satisfies
the conditions of Lemmata 9.7 and 9.8 with $\pi$, $S$ and $S^+$
being respectively $\pi_{\ccup H}$, $\{\pi_X\mid X\in
H\mn\{\bigcup H\}\}$ and $\{\pi_X^+\mid X\in H\mn\{\bigcup H\}\}$.
So there is a bijection from the set of facets of $\G(H)$ to the
set $\Psi$ of all pairs $\{\pi_X,\pi_{\ccup H}\}$ such that
\begin{tabbing}
\hspace{1.5em}(1)\quad$X\in H\mn\{\bigcup H\}$ and $\pi_X$
contains a vertex of $\G(H)$.
\end{tabbing}
Since for every $X$ in $H$ there is a construction of $H$ to which
it belongs (see Proposition 4.8), and since $\bigcup H$ belongs to
every construction of $H$, from (1) we infer that
\begin{tabbing}
\hspace{1.5em}(2)\quad$X\!\neq\bigcup H$ and there is a
construction $K$ of $H$ such that $\{X,\bigcup H\}\!\subseteq\!
K$.
\end{tabbing}
(We could make the same inference with the help of Proposition 9.3
too.) We use Proposition 9.6 to show that (1) follows from (2),
and hence $\Psi=\{\{\pi_X,\pi_{\ccup H}\}\mid \mbox{(2) holds}\}$.

Since $\pi$ is one-one, we obtain a bijection $\beta^\pi$ from the
set of facets of $\G(H)$ to the set of facets of $\A(H)$. From
Corollary 9.31 and Proposition 9.6 we infer that there is a
bijection $\gamma$ from the set of vertices of $\G(H)$ to the set
of vertices of $\A(H)$. The bijection $\gamma$ may be explicitly
defined by $\gamma(\{v\})=\{X\in H\mid \{v\}\subseteq\pi_X\}$. For
$\{v\}$ a vertex and $\varphi$ a facet of $\G(H)$, these two
bijections satisfy the following.

\prop{Lemma 9.9}{We have that $\{v\}\subseteq\varphi$ iff
$\beta^\pi(\varphi)\subseteq\gamma(\{v\})$.}

\dkz Let $\beta(\varphi)=\{\pi_X,\pi_{\ccup H}\}$. We have that
\begin{tabbing}
\hspace{1.5em}$\{v\}\subseteq\varphi$\quad\= iff\quad\=
$\{v\}\subseteq\bigcap\beta(\varphi) =\pi_X\cap\pi_{\ccup
H}$,\quad by Lemma 9.8
\\
\> iff\> $\beta^\pi(\varphi)= \{X,\bigcup
H\}\subseteq\gamma(\{v\})$,\quad by the definition of $\gamma$.
\`\qed
\end{tabbing}

Let $\X=\{X_1,\ldots,X_n\}$, for $n\geq 1$. Consider a hypergraph
$\H=\{\X_1,\ldots,\X_m\}$ on $\X$. Here we must have $m\geq 1$,
since $n\geq 1$. Let $\cal L$ be the lattice whose elements are in
the set
\[
\P(\X_1)\cup\ldots\cup\P(\X_m)\cup(\X\cup\{*\}).
\]
The join of $\cal L$ is set intersection $\cap$, while the meet
$\wedge$ for $A,B\in{\cal L}$ is defined by
\[
A\wedge B=\left\{
\begin{array}{ll}
A\cup B & {\mbox{\rm if }} (\exists i\in\{1,\ldots,m\}) A\cup
B\subseteq\X_i,
\\
\X\cup\{*\} & {\mbox{\rm otherwise}}.
\end{array}
\right.
\]
This lattice has a greatest element, namely $\emptyset$, and a
least element, namely $\X\cup\{*\}$. (A more natural lattice is
the dual lattice---namely, the same lattice taken upside
down---but the present lattice $\cal L$ is analogous to the
lattice $\A(H)$.) Consider a hypergraph
$\H'=\{\X'_1,\ldots,\X'_m\}$ on $\X'=\{X'_1,\ldots,X'_n\}$, and
let ${\cal L}'$ be defined exactly as $\cal L$; we just add the
primes. Suppose that for every $i\in\{1,\ldots,n\}$ and every
$j\in\{1,\ldots,m\}$ we have
\begin{tabbing}
\hspace{1.5em}(coinc)\quad$X_i\in\X_j$\quad iff\quad
$X'_i\in\X'_j$.
\end{tabbing}
Then it is trivial that the lattices $\cal L$ and ${\cal L}'$ are
isomorphic.

This is the situation we have with, on the one hand, $\X$ being
the set of facets of $\G(H)$, and the hypergraph $\H$ on $\X$
being the set obtained from the set of vertices of $\G(H)$ by
representing each vertex with the set of facets of $\G(H)$ in
which this vertex lies. It follows from Proposition 2.16 of
\cite{Z95} (Section 2.5) and Birkhoff's representation theorem for
finite Boolean algebras that for every $\X_i$ in $\H$ we have that
$\P(\X_i)$ is isomorphic to the Boolean algebra of all the faces
in which the vertex of $\G(H)$ corresponding to $\X_i$ lies. The
set $\H$ above is clearly in a bijection with the set of vertices
of $\G(H)$. The lattice $\cal L$ is isomorphic to the face lattice
of $\G(H)$.

On the other hand, take that $\X'$ is $H\mn\{\bigcup H\}$ for an
ASC-hypergraph $H$. This set is in a bijection with the set of
facets of $\A(H)$ because each facet of $\A(H)$ is of the form
$\{X,\bigcup H\}$ for an $X\in H\mn\{\bigcup H\}$. Let the
hypergraph $\H'$ on $\X'$ be the set of vertices of $\A(H)$ with
$\bigcup H$ removed. Each vertex of $\A(H)$ is a construction of
$H$, and we consider here the set $\H'=\{K\mn\{\bigcup H\}\mid K\;
\mbox{is a construction of}\; H\}$, which happens to be a
hypergraph on $H\mn\{\bigcup H\}$. It is clear that $\H'$ is in a
bijection with the set of vertices of $\A(H)$ (just remove
$\bigcup H$ from every vertex of $\A(H)$). It is easy to see that
the lattice ${\cal L}'$ is isomorphic to $\A(H)$ (just remove
$\bigcup H$ from every member of $\A(H)$).

The bijections $\beta^\pi$ and $\gamma$ above deliver two
bijections between $\X$ and $\X'$ and $\H$ and $\H'$ respectively,
and, due to Lemma 9.9, we have (coinc). So $\cal L$ and ${\cal
L}'$ are isomorphic, and hence the face lattice of $\G(H)$ is
isomorphic to $\A(H)$.

The foregoing covers the case when $d\geq 1$, and for $d<1$ we
obtain our isomorphism trivially. We have dealt with the case
$d=-1$, i.e.\ the case when $H=\emptyset$, at the beginning of the
section. When $d=0$, i.e.\ when $H$ is the singleton $\{\{1\}\}$,
then $\bigcap\{\pi_X^+\mid X\in\emptyset\}=\R^1$, and
$\pi_{\{1\}}=\{x_1\in\R^1\mid x_1=3\}=\{3\}$; so $\G(H)=\{3\}$.
The face lattice of $\G(H)$ has one vertex $\{3\}$, and it has
$\emptyset$ as the image of $F_{-1}$. The abstract polytope
$\A(H)$ is now $\{\{\{1\}\},\{\{1\},*\}\}$, with the vertex
$\{\{1\}\}$ being $F_0$ and $\{\{1\},*\}$ being $F_{-1}$. It is
clear that $\G(H)$ and $\A(H)$ are isomorphic.

When $H$ is atomic and saturated, but not connected, and
$\{H_1,\ldots,H_n\}$, for $n\geq 2$, is the finest hypergraph
partition of $H$, then $\G(H)$ is defined as
$\G(H_1)\times\ldots\times\G(H_n)$, with $\times$ being cartesian
product (polytopes are closed under $\times$; see \cite{Z95},
Lecture~0, pp.\ 9-10). The face lattice of $\G(H)$ so defined is
again isomorphic to $\A(H)$, which is defined as
$\A(H_1)\otimes\ldots\otimes\A(H_n)$ (see Sections 5 and~7).

When $H$ is atomic, but not saturated, we define the polytope
$\G(H)$ as the polytope $\G(\bar{H})$, for $\bar{H}$ being the
saturated closure of $H$, and again we obtain that the face
lattice of $\G(H)$ is isomorphic to $\A(H)$. So we may conclude
the following.

\prop{Proposition 9.10}{For every atomic hypergraph $H$ the face
lattice of $\G(H)$ is isomorphic to $\A(H)$.}

\appendix

\section{\large \textbf{Constructs and tubings}}
A notion which, as we shall see, is closely related to the notion
of construct of Section~5 was introduced under the name ``tubing''
in \cite{CD06}; it was modified in \cite{D09}, and modified
further in \cite{DF08} (which is posterior to \cite{D09}). In this
section, we will determine with the help of the results of
Section~6 the exact relationship between the tubings of
\cite{DF08} and constructs. For graphs, but not for hypergraphs in
general, the two notions happen to be equivalent. Along the way,
we obtain simpler characterizations of tubings than that given by
the definition of \cite{DF08}.

Problems arise for the tubings of \cite{CD06} and \cite{D09} with
connected graphs like
\begin{center}
\begin{picture}(100,20)

\put(0,10){\line(1,0){50}}

\put(50,10){\line(1,0){50}}

\put(1,10){\makebox(0,0){\circle*{2}}}
\put(51,10){\makebox(0,0){\circle*{2}}}
\put(101,10){\makebox(0,0){\circle*{2}}}

\put(1,7){\makebox(0,0)[t]{$x$}}

\put(51,7){\makebox(0,0)[t]{$y$}}

\put(101,7){\makebox(0,0)[t]{$z$}}
\end{picture}
\end{center}
We have that $\{\{x\},\{y,z\}\}$ is a tubing because $\{x\}$ and
$\{y,z\}$ are disjoint and are not adjacent (their union
$\{x,y,z\}$ is not a tube in the sense of \cite{CD06} and
\cite{D09}, because it is not a proper subset of vertices). This
tubing is however rejected in Fig. 1(b) of \cite{D09}. A simpler
problematic example is with the graph
\begin{center}
\begin{picture}(50,20)

\put(0,10){\line(1,0){50}}

\put(1,10){\makebox(0,0){\circle*{2}}}
\put(51,10){\makebox(0,0){\circle*{2}}}

\put(1,7){\makebox(0,0)[t]{$x$}}

\put(51,7){\makebox(0,0)[t]{$y$}}
\end{picture}
\end{center}
and the tubing $\{\{x\},\{y\}\}$ (in the sense of \cite{CD06} and
\cite{D09}). It is presumably because of such problems that the
modifications of \cite{DF08} were introduced in the definition of
a tubing, and we will consider here only this last modified
definition.

A tubing is defined relative to a \emph{graph}, which we will
identify with an ASC-hypergraph $G$ that is the saturated closure
$\bar{H}$ of a nonempty atomic hypergraph $H$ such that $\bigcup
H\in H$, and every member of $H$ that is neither a singleton nor
$\bigcup H$ is a two-element set; these two-element sets are the
edges of the graph $G$, and its vertices are the members of
$\bigcup G$, which is equal to $\bigcup H$. The remaining members
of $G$ are the connected subsets of the graph as they are usually
conceived, except for $\bigcup G$, which is in $G$ even if the
underlying graph is not connected in the usual sense. We put
always $\bigcup G$ in $G$ to match what is in \cite{DF08}; this
assumption, which yields connectedness, enables us also to
identify a graph with a kind of ASC-hypergraph. (The nonemptiness
of $G$, which means that $\bigcup G\neq\emptyset$, is not stated
explicitly in \cite{DF08}, but it is a common assumption for
graphs, which seems to be made because of the treatment of tubes
in \cite{DF08}, as we shall see in a moment; cf.\ the comment
after Remark 2.1).

A \emph{tube} of a graph $G$ is a member of $G$. Members of $G$
are always nonempty, and in \cite{CD06} tubes are said to be
nonempty. In \cite{DF08} this is not stated explicitly, but may be
taken to follow from the nonemptiness of graphs.

Two nonempty sets $X$ and $Y$ are \emph{overlapping} when $X\cap
Y\neq\emptyset$ and neither $X\subseteq Y$ nor $Y\subseteq X$ (in
\cite{DF08} one finds ``intersect'' instead of ``overlap''); they
are said to be \emph{adjacent} relative to a hypergraph $H$ when
$X\cap Y=\emptyset$ and $X\cup Y\in H$ (see the comments after
Lemma 6.3).

We will say that an ASC-hypergraph is \emph{loose} when
$H\mn\{\bigcup H\}$ is not a connected hypergraph, which means
that $\bigcup H$ is not dispensable in $H$ (see Section~4). Loose
graphs would normally be considered unconnected, but as
hypergraphs they are connected.

For a loose hypergraph $H$, let $\{H_1,\ldots,H_n\}$, with $n\geq
2$, be the finest hypergraph partition of $H\mn\{\bigcup H\}$, and
let $V_H$ be $\{\bigcup H_1,\ldots,\bigcup H_n\}$, which is a
partition of $\bigcup H$.

According to \cite{DF08}, a \emph{tubing} of a graph $G$ is a set
$T$ of tubes of $G$ (i.e.\ a subset of $G$) such that every pair
of tubes in $T$ is neither overlapping nor adjacent relative to
$G$; moreover, if $G$ is loose, then $V_G$ is not a subset of $T$,
and, finally, $\bigcup G\in T$.

Consider the assumption stated in Section~6:
\begin{itemize}
\item[(M)]{$M$ is a subset of the ASC-hypergraph $H$ such that
every $M$-antichain misses~$H$.}
\end{itemize}
We can prove the following.

\prop{Lemma A1}{If \mbox{\rm (M)} and $X,Y\in M$, then $X$ and $Y$
are not overlapping.}

\dkz If $X$ and $Y$ were overlapping, then $\{X,Y\}$ would be an
$M$-antichain that would not miss $H$, because of the saturation
of $H$.\qed

\oprop{Lemma A2}{If \mbox{\rm (M)} and $X,Y\in M$, then $X$ and
$Y$ are not adjacent.}

\dkz This holds simply because if $X\cap Y=\emptyset$, then
$\{X,Y\}$ is an $M$-antichain, because $X$ and $Y$ are
nonempty.\qed

\oprop{Lemma A3}{If \mbox{\rm (M)} and $H$ is loose, then $V_H$ is
not a subset of $M$.}

\dkz If $H\mn\{\bigcup H\}$ is not connected and $V_H\subseteq M$,
then $V_H$ is an $M$-antichain. This $M$-antichain does not miss
$H$, since $\bigcup V_H=\bigcup H\in H$, which contradicts
(M).\qed

As a corollary of these three lemmata we have the following.

\prop{Proposition A4}{If \mbox{\rm (M)} and $\bigcup H\in M$, then
$M$ is a tubing of $H$.}

We can also prove the converse for graphs $H$. For that we need
the following lemmata.

\prop{Lemma A5}{For $H$ a graph that is not loose and
$\{X_1,\ldots,X_n\}\subseteq H$, with $n\geq 2$, if for every
$i,j\in\{1,\ldots,n\}$ such that $i\neq j$ we have that $X_i\cup
X_j\not\in H$, then $X_1\cup\ldots\cup X_n\not\in H$.}

\dkz Since $H$ is a graph and is not loose, it is equal to the
saturated closure of an atomic hypergraph $H'$ in which as members
besides singletons we have only two-element sets. Suppose
$X_1\cup\ldots\cup X_n\in H$. Choose two distinct elements $x$ and
$y$ from respectively two different members of
$\{X_1,\ldots,X_n\}$; these elements exist because $n\geq 2$.
There is a path of $H$ connecting $x$ with $y$, and since $H$ and
$H'$ are cognate, by Proposition 4.6, there is a path
$Y_1,\ldots,Y_m$ of $H'$ connecting $x$ with $y$. Since $x$ and
$y$ are distinct, all the members of $Y_1,\ldots,Y_m$ may be taken
as two-element sets, and for each of the two elements there is an
$i\in\{1,\ldots,n\}$ such that this element belongs to $X_i$.

Since $x$ and $y$ are from two different members of
$\{X_1,\ldots,X_n\}$, for some $l\in\{1,\ldots,m\}$ we have that
$Y_l=\{x_i,x_j\}$ for $x_i\in X_i$ and $x_j\in X_j$, with
$i,j\in\{1,\ldots,n\}$ and $i\neq j$. So we have that
$X_i,Y_l,X_j\in H$, together with $X_i\cap Y_l\neq\emptyset$,
$X_j\cap Y_l\neq\emptyset$ and $Y_l\subseteq X_i\cup X_j$. Since
$H$ is saturated, we may conclude that $X_i\cup X_j\in H$, which
contradicts our assumption.\qed

\noindent This lemma is trivial when $n=2$. What we need for the
application in the next lemma are the cases $n\geq 3$.

\prop{Lemma A6}{For $H$ a graph that is not loose, an
$H$-antichain such that every pair of its members is neither
overlapping nor adjacent relative to $H$ misses $H$.}

\dkz Take an $H$-antichain $\{X_1,\ldots,X_n\}$, for $n\geq 2$. If
$n=2$, then we have that $X_1$ and $X_2$ are not overlapping, and
since neither $X_1\subseteq X_2$ nor $X_2\subseteq X_1$, we have
that $X_1\cap X_2=\emptyset$. Since $X_1$ and $X_2$ are not
adjacent, it follows that $X_1\cup X_2\not\in H$.

For $n\geq 3$, we have for every $i,j\in\{1,\ldots,n\}$ such that
$i\neq j$ that $\{X_i, X_j\}$ is an $H$-antichain. According to
what we have just proved, $X_i\cup X_j\not\in H$. By Lemma A5, we
may then conclude that $X_1\cup\ldots\cup X_n\not\in H$.\qed

We can then prove for graphs the following converse of Proposition
A4.

\prop{Proposition A7}{For $H$ a graph, if $M$ is a tubing of $H$,
then every $M$-antichain misses $H$ and $\bigcup H\in M$.}

\dkz Suppose $\{H_1,\ldots,H_n\}$, for $n\geq 1$, is the finest
hypergraph partition of $H\mn\{\bigcup H\}$. Take an $M$-antichain
$S$ of the form $S_1\cup\ldots\cup S_n$ such that for every
$i\in\{1,\ldots,n\}$ we have that $S_i\subseteq H_i$.

Note first that it is impossible that every member of
$\{S_1,\ldots,S_n\}$ is empty. If just one of these members is
nonempty, then $S\subseteq H_i$ for some $i\in\{1,\ldots,n\}$.
Then we have that either $n\geq 2$ and $H_i$ is a graph that is
not loose, or $n=1$ and $H_i=H_1=H\mn\{\bigcup H\}$; in the second
case we have that $S\subseteq H$ and $H$ is not loose. So $S$ is
either an $H_i$-antichain or an $H$-antichain. In both cases, by
Lemma A6, we conclude that $\bigcup S\not\in H$.

If at least two members of $\{S_1,\ldots,S_n\}$ are nonempty and
at least one member is empty, then $\bigcup S\subset\bigcup H$,
and we may conclude that $\bigcup S\not\in H$. It remains to
consider the case when all the members of $\{S_1,\ldots,S_n\}$ are
nonempty and $n\geq 2$. In that case $H$ is loose.

Since for every $i$ we have that $S_i\subseteq H_i$, we must have
that $\bigcup S_i\subseteq \bigcup H_i$. Suppose that for every
$i$ we have that $\bigcup S_i=\bigcup H_i$. We can conclude that
$S_i=\{\bigcup H_i\}$. Otherwise, $|S_i|\geq 2$, and so $S_i$ is
an $H_i$-antichain; by Lemma A6, it misses $H_i$, which
contradicts $\bigcup S_i=\bigcup H_i$. So $S=\{\bigcup
H_1,\ldots,\bigcup H_n\}=V_H$, and this contradicts the assumption
that $V_H$ is not a subset of $M$, which we have when $H$ is
loose. So for some $i$ we have that $\bigcup S_i\subset\bigcup
H_i$, and then $\bigcup S\not\in H$.

The condition $\bigcup H\in M$ is assumed for tubings.\qed

As a corollary of Propositions A4 and A7, we have that for graphs
$H$ a subset $M$ of $H$ is a tubing of $H$ iff every $M$-antichain
misses $H$ and $\bigcup H\in M$. This gives an alternative,
simpler, definition of a tubing.

With this characterization of tubings, we can immediately infer
from Proposition 6.13 that for graphs $H$ a subset $M$ of $H$ is a
tubing of $H$ iff $M$ is a construct of $H$. We can infer from
Proposition 6.11 that for graphs $H$ a subset $M$ of $H$ is a
tubing of $H$ and $|M|=|\bigcup H|$ iff $M$ is a construction of
$H$.

Note that we have these characterizations of constructs and
constructions in terms of tubings only for graphs, i.e.\ for
specific hypergraphs with which we have identified graphs. We do
not have them for hypergraphs in general.

Another difference of our approach with the approach through
tubings is that for us constructs, which for graphs amount to
tubings, are a derived, secondary, notion. The basic, primary,
notion is the notion of construction.

\section{\large \textbf{Hypergraph polytopes of dimension 3 and lower}}

In this appendix we survey the abstract polytopes $\A(H)$ of
atomic saturated hypergraphs $H$ with carriers $\bigcup H$ having
no more than four elements. In this survey, we deal with the main
types of these polytopes, and do not distinguish polytopes that
would differ only up to renaming the elements of the carrier.

We name the hypergraphs in our survey by adding to $H$ subscripts
according to the following system (sometimes we also have
superscripts). In $H_{i_1,\ldots,i_k}$, the subscript $i_j$ for
$j\in\{1,\ldots,k\}$ is the number of $j$-element members of
$H_{i_1,\ldots,i_k}$. So, for example, $H_{21}$ below has 2
singletons and one pair. Since our hypergraphs are atomic, we
always have that the first subscript $i_1$ is the cardinality of
the carrier, while the last subscript $i_k$ can be either 0 or 1.

\vspace{1ex}

\hspace{-2.5em} $(H_0,H_1)$\quad Let $H_0$ be the empty hypergraph
$\emptyset$; then
$\A(H_0)=\{H_0,H_0\cup\{*\}\}=\{\emptyset,\{*\}\}$. Let $H_1$ be
the hypergraph $\{\{x\}\}$; then $\A(H_1)=\{H_1,H_1\cup\{*\}\}$.
With that we have surveyed all we have with the carrier of the
hypergraph having not more than one element.

\vspace{1ex}

\hspace{-2.5em} $(H_{20},H_{21})$\quad With the carrier having two
elements, we have two atomic hypergraphs: $H_{20}=\{\{x\},\{y\}\}$
and $H_{21}=H_{20}\cup\{\{x,y\}\}$. Analogously to what we had in
$(H_0,H_1)$, we have that $\A(H_{20})=\{H_{20},H_{20}\cup\{*\}\}$,
and $\A(H_{21})$ has the following structure:
\begin{tabbing}
\hspace{1.5em}\=$F_1$ (edge)\hspace{10em}$\{\{x,y\}\}$\\*[1ex]

\>vertices\hspace{5.5em}$\{\{x\},\{x,y\}\}$\hspace{3em}$\{\{y\},\{x,y\}\}$\\*[1ex]

\>$F_{-1}$\hspace{7.2em}$H_{21}\cup\{*\}=\{\{x\},\{y\},\{x,y\},*\}$
\end{tabbing}
A realization of $\A(H_{21})$ may be pictured as
\begin{center}
\begin{picture}(50,18)

\put(0,10){\line(1,0){50}}

\put(1,7){\makebox(0,0)[t]{$\{x\}$}}

\put(51,7){\makebox(0,0)[t]{$\{y\}$}}
\end{picture}
\end{center}
where in the labels of the vertices we have omitted $\{x,y\}$ and
the outermost braces.

\vspace{1ex}

We pass next to atomic saturated hypergraphs whose carrier has
three elements.

\vspace{1ex}

\hspace{-2.5em} $(H_{300},H_{310})$\quad We have first
$H_{300}=\{\{x\},\{y\}\}$, with
$\A(H_{300})=\{H_{300},H_{300}\cup\{*\}\}$. Next we have
$H_{310}=H_{300}\cup\{\{x,y\}\}$, with $\A(H_{310})$ being
isomorphic to (i.e.\ being in an order-preserving bijection with)
$\A(H_{21})$, which we write $\A(H_{310})\cong\A(H_{21})$.

Cases like these with $H_{300}$ and $H_{310}$, where the
hypergraph is not connected (which happens when there is more than
one subscript, and the last is 0), will be called
\emph{degenerate}. (So $H_{20}$ above is degenerate.) In general,
in a degenerate case the rank of $\A(H_{k\ldots l})$ is lower than
$k\mn 1$. So the rank of $\A(H_{310})$ is 1, while in the four
non-degenerate cases with $k=3$, which follow, it will be~2.

\vspace{1ex}

\hspace{-2.5em} $(H_{301},H_{311},H_{321},H_{331})$\quad As
non-degenerate cases with the carrier having three elements, we
have the hypergraphs $H$ on the left, with the corresponding
realizations of $\A(H)$ pictured on the right:
\begin{center}
\begin{picture}(0,70)(-50,-5)
\put(30,0){\line(-3,5){30}} \put(30,0){\line(3,5){30}}
\put(0,50){\line(1,0){60}} \put(12,25){\makebox(0,0)[r]{$\{x\}$}}
\put(48,25){\makebox(0,0)[l]{$\{y\}$}}
\put(30,55){\makebox(0,0)[b]{$\{z\}$}}

\put(-203,35){\makebox(0,0)[l]{$H_{301}=H_{300}\cup\{\{x,y,z\}\}$}}
\end{picture}
\end{center}
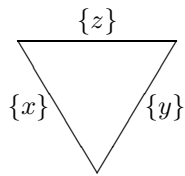
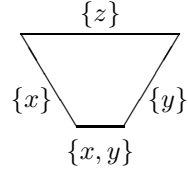
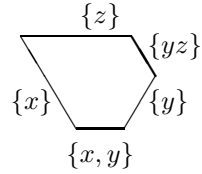
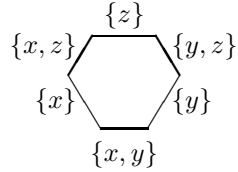
\begin{center}
\begin{picture}(0,80)(-50,-5)
\put(21,15){\line(-3,5){21}} \put(39,15){\line(3,5){21}}
\put(21,15){\thicklines \line(1,0){18}} \put(0,50){\line(1,0){60}}
\put(12,25){\makebox(0,0)[r]{$\{x\}$}}
\put(48,25){\makebox(0,0)[l]{$\{y\}$}}
\put(30,55){\makebox(0,0)[b]{$\{z\}$}}
\put(30,10){\makebox(0,0)[t]{$\{x,y\}$}}

\put(-203,35){\makebox(0,0)[l]{$H_{311}=H_{301}\cup\{\{x,y\}\}$}}
\end{picture}
\end{center}
\begin{center}
\begin{picture}(0,80)(-50,-5)
\put(21,15){\line(-3,5){21}} \put(39,15){\line(3,5){12}}
\put(51,35){\thicklines\line(-3,5){9}} \put(21,15){\thicklines
\line(1,0){18}} \put(0,50){\line(1,0){42}}
\put(12,25){\makebox(0,0)[r]{$\{x\}$}}
\put(48,25){\makebox(0,0)[l]{$\{y\}$}}
\put(30,55){\makebox(0,0)[b]{$\{z\}$}}
\put(30,10){\makebox(0,0)[t]{$\{x,y\}$}}
\put(48,46){\makebox(0,0)[l]{$\{yz\}$}}

\put(-203,35){\makebox(0,0)[l]{$H_{321}=H_{311}\cup\{\{y,z\}\}$}}
\end{picture}
\end{center}
\begin{center}
\begin{picture}(0,80)(-50,-5)
\put(21,15){\line(-3,5){12}} \put(39,15){\line(3,5){12}}
\put(9,35){\thicklines\line(3,5){9}}
\put(51,35){\thicklines\line(-3,5){9}} \put(21,15){\thicklines
\line(1,0){18}} \put(18,50){\line(1,0){24}}
\put(12,25){\makebox(0,0)[r]{$\{x\}$}}
\put(48,25){\makebox(0,0)[l]{$\{y\}$}}
\put(30,55){\makebox(0,0)[b]{$\{z\}$}}
\put(30,10){\makebox(0,0)[t]{$\{x,y\}$}}
\put(48,46){\makebox(0,0)[l]{$\{y,z\}$}}
\put(12,46){\makebox(0,0)[r]{$\{x,z\}$}}

\put(-203,35){\makebox(0,0)[l]{$H_{331}=H_{321}\cup\{\{x,z\}\}$}}
\end{picture}
\end{center}

For the labels of the edges in the pictures on the right we have
the convention that $\{x,y,z\}$ and the outermost braces are
omitted; when they are restored, we obtain the edges of $\A(H)$. A
vertex of $\A(H)$ is obtained by taking the set made of the labels
of the edges that are incident with this vertex plus $\{x,y,z\}$.
Finally, $F_2$ is here always $\{\{x,y,z\}\}$, which corresponds
to the whole polygon.

Without all these abbreviations, the first picture would be
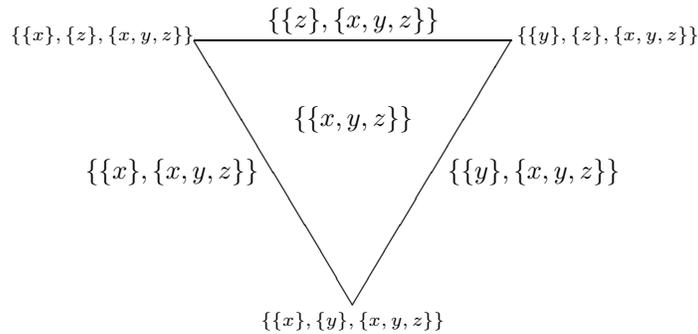
\begin{center}
\begin{picture}(120,125)(0,-10)
\put(60,0){\line(-3,5){60}} \put(60,0){\line(3,5){60}}
\put(0,100){\line(1,0){120}}
\put(24,50){\makebox(0,0)[r]{$\{\{x\},\{x,y,z\}\}$}}
\put(96,50){\makebox(0,0)[l]{$\{\{y\},\{x,y,z\}\}$}}
\put(60,105){\makebox(0,0)[b]{$\{\{z\},\{x,y,z\}\}$}}
\put(60,70){\makebox(0,0){$\{\{x,y,z\}\}$}}
\put(60,-3){\makebox(0,0)[t]{\scriptsize$\{\{x\},\{y\},\{x,y,z\}\}$}}
\put(0,100){\makebox(0,0)[br]{\scriptsize$\{\{x\},\{z\},\{x,y,z\}\}$}}
\put(122,100){\makebox(0,0)[bl]{\scriptsize$\{\{y\},\{z\},\{x,y,z\}\}$}}
\end{picture}
\end{center}
where the labels are the members of $\A(H_{301})$ without
$F_{-1}$, which is
$\bar{H}^*_{301}=\{\{x\},\{y\},\{z\},\{x,y,z\},*\}$ (we have
$\bar{H}_{301}=H_{301})$. This is, of course, the picture of the
two-dimensional simplex.

Without the abbreviations, for the vertices of the last,
hexagonal, picture we would have the labels below, for which we
also write underneath the corresponding s-constructions:
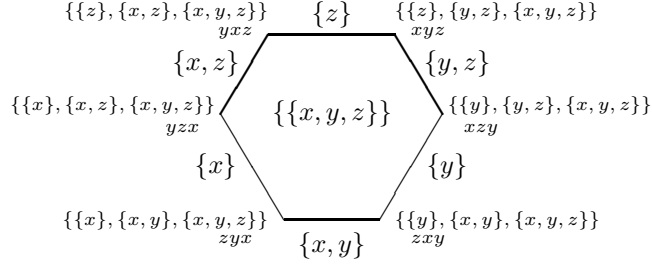
\begin{center}
\begin{picture}(120,110)(0,10)
\put(42,30){\line(-3,5){24}} \put(78,30){\line(3,5){24}}
\put(18,70){\thicklines\line(3,5){18}}
\put(102,70){\thicklines\line(-3,5){18}} \put(42,30){\thicklines
\line(1,0){36}} \put(36,100){\line(1,0){48}}
\put(24,50){\makebox(0,0)[r]{$\{x\}$}}
\put(96,50){\makebox(0,0)[l]{$\{y\}$}}
\put(60,105){\makebox(0,0)[b]{$\{z\}$}}
\put(60,25){\makebox(0,0)[t]{$\{x,y\}$}}
\put(95,89){\makebox(0,0)[l]{$\{y,z\}$}}
\put(25,89){\makebox(0,0)[r]{$\{x,z\}$}}

\put(36,107){\makebox(0,0)[br]{\scriptsize$\{\{z\},\{x,z\},\{x,y,z\}\}$}}
\put(30,100){\makebox(0,0)[br]{\scriptsize$yxz$}}

\put(84,107){\makebox(0,0)[bl]{\scriptsize$\{\{z\},\{y,z\},\{x,y,z\}\}$}}
\put(90,100){\makebox(0,0)[bl]{\scriptsize$xyz$}}

\put(16,73){\makebox(0,0)[r]{\scriptsize$\{\{x\},\{x,z\},\{x,y,z\}\}$}}
\put(10,64){\makebox(0,0)[r]{\scriptsize$yzx$}}

\put(104,73){\makebox(0,0)[l]{\scriptsize$\{\{y\},\{y,z\},\{x,y,z\}\}$}}
\put(110,64){\makebox(0,0)[l]{\scriptsize$xzy$}}

\put(36,33){\makebox(0,0)[tr]{\scriptsize$\{\{x\},\{x,y\},\{x,y,z\}\}$}}
\put(30,24){\makebox(0,0)[tr]{\scriptsize$zyx$}}

\put(84,33){\makebox(0,0)[tl]{\scriptsize$\{\{y\},\{x,y\},\{x,y,z\}\}$}}
\put(90,24){\makebox(0,0)[tl]{\scriptsize$zxy$}}

\put(60,70){\makebox(0,0){$\{\{x,y,z\}\}$}}
\end{picture}
\end{center}
These s-constructions correspond, of course, to the six
permutations of $x$, $y$ and~$z$.

For the third, pentagonal, picture, most of the labels would be
the same; the difference would be only in the left upper corner,
where we have a vertex labelled $\{\{x\},\{z\},\{x,y,z\}\}$, with
the corresponding s-construction being $y(x\pl z)$.

This third picture is the picture of the two-dimensional
associahedron, also known as Mac Lane's pentagon, and the last,
fourth, picture is the picture of the two-dimensional
permutohedron, also known as Mac Lane's hexagon (see \cite{ML63}
and \cite{ML98}, Sections VII.1 and VII.7, for Mac Lane's pentagon
and hexagon; see $(H'_{4321})$ and $(H_{4641})$ below for
references concerning associahedra and permutohedra). This
comment, and the connection with the labels for vertices written
as s-constructions, are explained in \cite{DP06} and \cite{DP10a}.

We pass from the triangle to the quadrilateral, the pentagon and
the hexagon by truncating the vertices in succession. This
truncating is explained in Section~9. Whereas here we truncate, in
\cite{DP10a} we find the converse operation of collapsing several
vertices into one. So the starting point would be not the simplex,
but the permutohedron, and the direction would be in this case
from the hexagon towards the triangle. Although the direction is
reversed, this does not differ essentially from what we have here.

\vspace{1ex}

In the remainder of this survey we have atomic hypergraphs whose
carrier has four elements.

\vspace{1ex}

\hspace{-2.5em} $(H_{4000},\ldots,H_{4310})$\quad As degenerate
cases, we have first the following:
\begin{tabbing}
\hspace{1.5em}\=$H_{4000}\:$\=$=\{\{x\},\{y\},\{z\},\{u\}\}$,
\hspace{4em}\=$\A(H_{4000})\:$\=$=\{H_{4000},
H_{4000}\cup\{*\}\}$,\\

\>$H_{4100}\:$\>$=H_{4000}\cup\{\{x,y\}\}$,
\>$\A(H_{4100})\:$\=$\cong\A(H_{21})$,\\

\>$H_{4010}\:$\>$=H_{4000}\cup\{\{x,y,z\}\}$,
\>$\A(H_{4010})\:$\=$\cong\A(H_{301})$,\\

\>$H_{4110}\:$\>$=H_{4010}\cup\{\{x,y\}\}$,
\>$\A(H_{4110})\:$\=$\cong\A(H_{311})$,\\

\>$H_{4210}\:$\>$=H_{4100}\cup\{\{y,z\},\{x,y,z\}\}$,
\>$\A(H_{4210})\:$\=$\cong\A(H_{321})$,\\

\>$H_{4310}\:$\>$=H_{4210}\cup\{\{x,z\}\}$,
\>$\A(H_{4310})\:$\=$\cong\A(H_{331})$
\end{tabbing}
(the case with $H_{4310}$ was investigated as Example 5.2 in
\cite{DP10a}).

\vspace{1ex}

\hspace{-2.5em} $(H_{4200})$\quad As the last degenerate case, we
have $H_{4200}=H_{4100}\cup\{\{z,u\}\}$, and a realization of
$\A(H_{4200})$ is pictured by
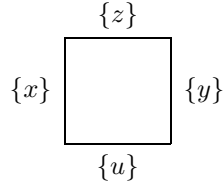
\begin{center}
\begin{picture}(40,60)
\put(0,10){\line(1,0){40}} \put(0,50){\line(1,0){40}}
\put(0,10){\line(0,1){40}} \put(40,10){\line(0,1){40}}

\put(-5,30){\makebox(0,0)[r]{$\{x\}$}}
\put(45,30){\makebox(0,0)[l]{$\{y\}$}}
\put(20,6){\makebox(0,0)[t]{$\{u\}$}}
\put(20,56){\makebox(0,0)[b]{$\{z\}$}}

\end{picture}
\end{center}
The polytope $\A(H_{4200})$ is obtained as the product with
$\otimes$ of two copies of $\A(H_{21})$ (see Section~5). We pass
next to non-degenerate cases.

\vspace{1ex}

\hspace{-2.5em} $(H_{4001})$\quad As the first non-degenerate case
with four elements in the carrier we have
$H_{4001}=H_{4000}\cup\{\{x,y,z,u\}\}$, with $\A(H_{4001})$ being
realized as the tetrahedron, i.e.\ the three-dimensional simplex.
In general, for every $k\geq 3$ we have that $\A(H_{k0\ldots01})$
may be realized as the ${(k\mn 1)}$-dimensional simplex (see the
realization of $\A(H_{301})$ above); $\A(H_{21})$ is realized as
the one-dimensional simplex, and $\A(H_{1})$ as the
zero-dimensional simplex (see above). The tetrahedron is pictured
by
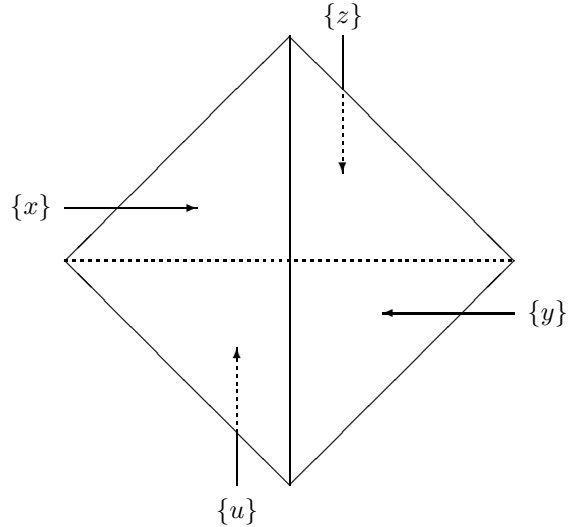
\begin{center}
\begin{picture}(170,190)(0,-10)
\put(85,0){\line(-1,1){85}} \put(85,0){\line(1,1){85}}
\put(85,170){\line(-1,-1){85}} \put(85,170){\line(1,-1){85}}
\multiput(0,85)(3,0){57}{\line(1,0){1}}
\put(85,0){\line(0,1){170}}

\put(0,105){\vector(1,0){50}}
\put(-5,105){\makebox(0,0)[r]{$\{x\}$}}

\put(170,65){\vector(-1,0){50}}
\put(175,65){\makebox(0,0)[l]{$\{y\}$}}

\put(65,0){\line(0,1){20}}
\multiput(65,20)(0,3){10}{\line(0,1){1}}
\put(65,50){\vector(0,1){2}}
\put(65,-5){\makebox(0,0)[t]{$\{u\}$}}

\put(105,170){\line(0,-1){20}}
\multiput(105,120)(0,3){10}{\line(0,1){1}}
\put(105,120){\vector(0,-1){2}}
\put(105,175){\makebox(0,0)[b]{$\{z\}$}}
\end{picture}
\end{center}

We label here only the facets of the tetrahedron, with
$\{x,y,z,u\}$ and the outermost braces omitted. The edges and
vertices may be reconstructed out of these labels. We just look
what facets are incident with the edge or the vertex. For example,
the north-west edge is $\{\{x\},\{z\},\{x,y,z,u\}\}$, and the
north vertex is $\{\{x\},\{y\},\{z\},\{x,y,z,u\}\}$. The whole
tetrahedron corresponds to $\{\{x,y,z,u\}\}$.

\vspace{1ex}

\hspace{-2.5em} $(H_{4011})$\quad In this case we truncate a
vertex. We have $H_{4011}=H_{4001}\cup\{\{x,y,z\}\}$, with a
realization of $\A(H_{4011})$ obtained from our tetrahedron by
truncating the vertex $\{\{x\},\{y\},\{z\},\{x,y,z,u\}\}$, which
is pictured by
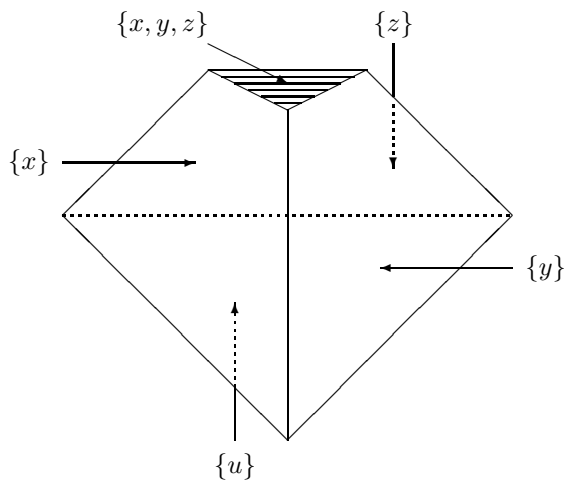
\begin{center}
\begin{picture}(170,175)(0,-10)
\put(85,0){\line(-1,1){85}} \put(85,0){\line(1,1){85}}
\put(55,140){\line(-1,-1){55}} \put(115,140){\line(1,-1){55}}
\multiput(0,85)(3,0){57}{\line(1,0){1}}
\put(85,0){\line(0,1){125}}

\put(55,140){\line(1,0){60}} \put(55,140){\line(2,-1){30}}
\put(115,140){\line(-2,-1){30}}

\put(60,137.5){\line(1,0){50}} \put(65,135){\line(1,0){40}}
\put(70,132.5){\line(1,0){30}} \put(75,130){\line(1,0){20}}
\put(80,127.5){\line(1,0){10}}

\put(0,105){\vector(1,0){50}}
\put(-5,105){\makebox(0,0)[r]{$\{x\}$}}

\put(170,65){\vector(-1,0){50}}
\put(175,65){\makebox(0,0)[l]{$\{y\}$}}

\put(65,0){\line(0,1){20}}
\multiput(65,20)(0,3){10}{\line(0,1){1}}
\put(65,50){\vector(0,1){2}}
\put(65,-5){\makebox(0,0)[t]{$\{u\}$}}

\put(125,150){\line(0,-1){20}}
\multiput(125,105)(0,3){8}{\line(0,1){1}}
\put(125,105){\vector(0,-1){2}}
\put(125,155){\makebox(0,0)[b]{$\{z\}$}}

\put(55,150){\vector(2,-1){30}}
\put(55,155){\makebox(0,0)[br]{$\{x,y,z\}$}}
\end{picture}
\end{center}
This polytope may also be realized as a three-sided prism.

\vspace{1ex}

\hspace{-2.5em}$(H_{4021},H_{4031},H_{4041})$\quad Next we have
$H_{4021}=H_{4011}\cup\{\{y,z,u\}\}$,
$H_{4031}=H_{4021}\cup\{\{x,z,u\}\}$ and
$H_{4041}=H_{4031}\cup\{\{x,y,u\}\}$, with $\A(H_{4021})$,
$\A(H_{4031})$ and $\A(H_{4041})$ being realized as the
tetrahedron in which we have truncated two, three and four
vertices respectively. (None of the last four cases is covered by
the approach of \cite{CD06} and \cite{D09}, which is based on
graphs, as explained in Appendix~A; we will call such cases
\emph{essentially hypergraphical}.)

\vspace{1ex}

\hspace{-2.5em} $(H_{4101})$\quad In this case we truncate an
edge. We have $H_{4101}=H_{4001}\cup\{\{x,y\}\}$, with a
realization of $\A(H_{4101})$ obtained from our tetrahedron by
truncating the edge $\{\{x\},\{y\},\{x,y,z,u\}\}$, which is
pictured by
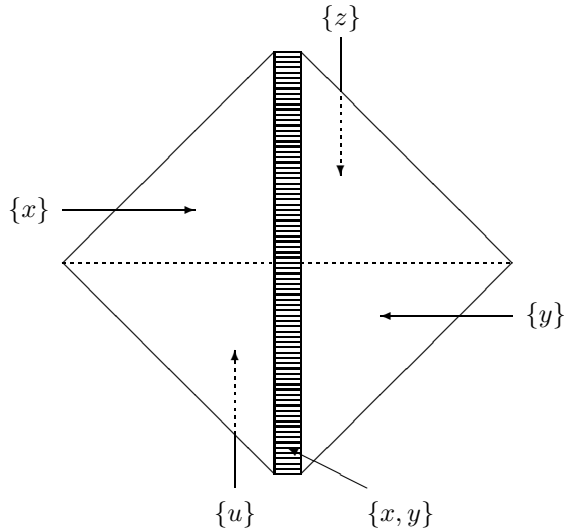
\begin{center}
\begin{picture}(170,190)(0,-10)
\put(80,5){\line(-1,1){80}} \put(90,5){\line(1,1){80}}
\put(80,165){\line(-1,-1){80}} \put(90,165){\line(1,-1){80}}
\multiput(0,85)(3,0){57}{\line(1,0){1}}
\put(80,5){\line(0,1){160}} \put(90,5){\line(0,1){160}}
\put(80,5){\line(1,0){10}} \put(80,165){\line(1,0){10}}

\multiput(80,5)(0,2){80}{\line(1,0){10}}

\put(0,105){\vector(1,0){50}}
\put(-5,105){\makebox(0,0)[r]{$\{x\}$}}

\put(170,65){\vector(-1,0){50}}
\put(175,65){\makebox(0,0)[l]{$\{y\}$}}

\put(65,0){\line(0,1){20}}
\multiput(65,20)(0,3){10}{\line(0,1){1}}
\put(65,50){\vector(0,1){2}}
\put(65,-5){\makebox(0,0)[t]{$\{u\}$}}

\put(105,170){\line(0,-1){20}}
\multiput(105,120)(0,3){10}{\line(0,1){1}}
\put(105,120){\vector(0,-1){2}}
\put(105,175){\makebox(0,0)[b]{$\{z\}$}}

\put(115,0){\vector(-2,1){30}}
\put(115,-5){\makebox(0,0)[tl]{$\{x,y\}$}}

\end{picture}
\end{center}
This polytope, as the preceding one, viz.\ $\A(H_{4011})$, may be
realized as a three-sided prism.

\vspace{1ex}

\hspace{-2.5em} $(H_{4201})$\quad With two opposite edges
truncated, we have $H_{4201}=H_{4101}\cup\{\{z,u\}\}$, where a
realization of $\A(H_{4201})$ is pictured by
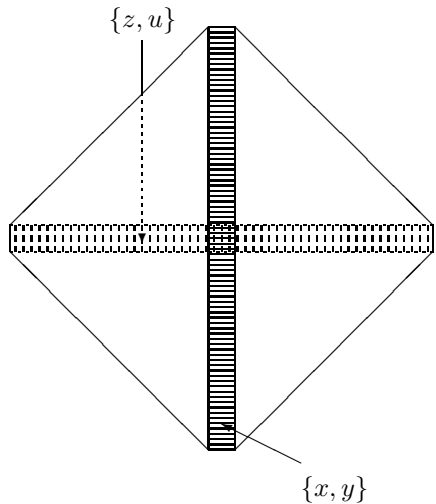
\begin{center}
\begin{picture}(170,190)(0,-10)
\put(80,5){\line(-1,1){75}} \put(90,5){\line(1,1){75}}
\put(80,165){\line(-1,-1){75}} \put(90,165){\line(1,-1){75}}
\multiput(5,90)(3,0){54}{\line(1,0){1}}
\multiput(5,80)(3,0){54}{\line(1,0){1}}
\put(80,5){\line(0,1){160}} \put(90,5){\line(0,1){160}}
\put(80,5){\line(1,0){10}} \put(80,165){\line(1,0){10}}
\put(5,80){\line(0,1){10}} \put(165,80){\line(0,1){10}}

\multiput(80,5)(0,2){80}{\line(1,0){10}}
\multiput(7,80)(3,0){53}{\line(0,1){2}}
\multiput(7,84)(3,0){53}{\line(0,1){2}}
\multiput(7,88)(3,0){53}{\line(0,1){2}}

\put(55,160){\line(0,-1){20}}
\multiput(55,85)(0,3){20}{\line(0,1){1}}
\put(55,85){\vector(0,-1){2}}
\put(55,165){\makebox(0,0)[b]{$\{z,u\}$}}

\put(115,0){\vector(-2,1){30}}
\put(115,-5){\makebox(0,0)[tl]{$\{x,y\}$}}

\end{picture}
\end{center}
with the labels for facets $\{x\}$, $\{y\}$, $\{z\}$ and $\{u\}$
omitted; they will mostly be omitted from now on. This polytope
may also be realized as a cube.

\vspace{1ex}

\hspace{-2.5em} $(H_{4111})$\quad With one edge and one incident
vertex truncated, we have $H_{4111}=H_{4011}\cup H_{4101}$, where
a realization of $\A(H_{4111})$ is pictured by
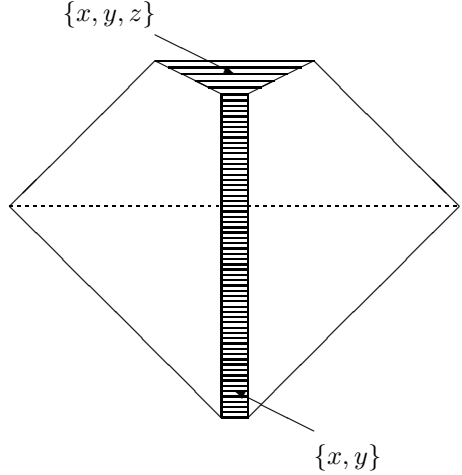
\begin{center}
\begin{picture}(170,175)(0,-10)
\put(80,5){\line(-1,1){80}} \put(90,5){\line(1,1){80}}
\put(55,140){\line(-1,-1){55}} \put(115,140){\line(1,-1){55}}
\multiput(0,85)(3,0){57}{\line(1,0){1}}
\put(80,5){\line(0,1){122.5}} \put(90,5){\line(0,1){122.5}}
\put(80,5){\line(1,0){10}}

\put(55,140){\line(1,0){60}} \put(55,140){\line(2,-1){25}}
\put(115,140){\line(-2,-1){25}} \put(80,127.5){\line(1,0){10}}

\multiput(80,5)(0,2){62}{\line(1,0){10}}

\put(60,137.5){\line(1,0){50}} \put(65,135){\line(1,0){40}}
\put(70,132.5){\line(1,0){30}} \put(75,130){\line(1,0){20}}

\put(55,150){\vector(2,-1){30}}
\put(55,155){\makebox(0,0)[br]{$\{x,y,z\}$}}

\put(115,0){\vector(-2,1){30}}
\put(115,-5){\makebox(0,0)[tl]{$\{x,y\}$}}

\end{picture}
\end{center}
This polytope, as the preceding one, viz.\ $\A(H_{4201})$, may be
realized as a cube.

\vspace{1ex}

\hspace{-2.5em} $(H'_{4111})$\quad With one edge and one
non-incident vertex truncated, we have
$H'_{4111}=H_{4101}\cup\{\{y,z,u\}\}$, with the picture of a
realization of $\A(H'_{4111})$ obtained from that given for
$\A(H_{4101})$ by truncating the east, i.e.\ right, vertex
$\{\{y\},\{z\},\{u\},$ $\{x,y,z,u\}\}$.

\vspace{1ex}

\hspace{-2.5em} $(H_{4121})$\quad With one edge and two incident
vertices truncated, we have $H_{4121}=H_{4111}\cup\{\{x,y,u\}\}$,
where a realization of $\A(H_{4121})$ is pictured by
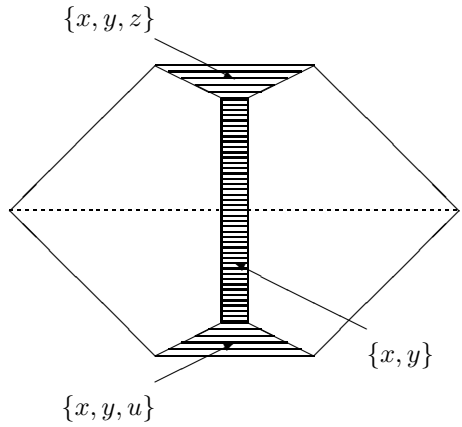
\begin{center}
\begin{picture}(170,150)(0,10)
\put(55,30){\line(-1,1){55}} \put(115,30){\line(1,1){55}}
\put(55,140){\line(-1,-1){55}} \put(115,140){\line(1,-1){55}}
\multiput(0,85)(3,0){57}{\line(1,0){1}}
\put(80,42.5){\line(0,1){85}} \put(90,42.5){\line(0,1){85}}
\put(80,42.5){\line(1,0){10}}

\put(55,140){\line(1,0){60}} \put(55,140){\line(2,-1){25}}
\put(115,140){\line(-2,-1){25}} \put(80,127.5){\line(1,0){10}}

\put(55,30){\line(1,0){60}} \put(55,30){\line(2,1){25}}
\put(115,30){\line(-2,1){25}} \put(80,42.5){\line(1,0){10}}

\multiput(80,43)(0,2){43}{\line(1,0){10}}

\put(60,137.5){\line(1,0){50}} \put(65,135){\line(1,0){40}}
\put(70,132.5){\line(1,0){30}} \put(75,130){\line(1,0){20}}

\put(60,32.5){\line(1,0){50}} \put(65,35){\line(1,0){40}}
\put(70,37.5){\line(1,0){30}} \put(75,40){\line(1,0){20}}

\put(55,150){\vector(2,-1){30}}
\put(55,155){\makebox(0,0)[br]{$\{x,y,z\}$}}

\put(55,20){\vector(2,1){30}}
\put(55,15){\makebox(0,0)[tr]{$\{x,y,u\}$}}

\put(135,40){\vector(-2,1){50}}
\put(135,35){\makebox(0,0)[tl]{$\{x,y\}$}}

\end{picture}
\end{center}
This polytope may also be realized as a five-sided prism.

\vspace{1ex}

\hspace{-2.5em} $(H'_{4121})$\quad As a case where we truncate one
edge and two vertices, one incident and the other not, we have
$H'_{4121}=H_{4111}\cup H'_{4111}$. The picture of a realization
of $\A(H'_{4121})$ is obtained from that given for $H_{4111}$ by
truncating the right vertex.

\vspace{1ex}

\hspace{-2.5em} $(H''_{4121},H_{4131},H'_{4131},H_{4141})$\quad As
a case where we truncate one edge and two vertices, none of them
incident, we have $H''_{4121}=H_{4111}\cup\{\{x,z,u\}\}$. Next we
have two cases where we truncate one edge and three vertices:
$H_{4131}=H'_{4111}\cup H_{4121}$ and $H'_{4131}=H_{4111}\cup
H''_{4121}$, and one case where we truncate one edge and all the
four vertices: $H_{4141}=H_{4131}\cup H'_{4041}$. In all these
cases, it should be clear by now how to picture a realization of
$\A(H)$ starting from previous pictures. (The last eight cases are
essentially hypergraphical.)

\vspace{1ex}

\hspace{-2.5em} $(H_{4211})$\quad If we truncate two edges of our
tetrahedron that are not opposite, but are incident with a common
vertex, then we must truncate this vertex too (which is something
related to saturation). This happens with $H_{4211}=H_{4111}\cup
\{\{y,z\}\}$, with the picture of a realization of $\A(H_{4211})$
being
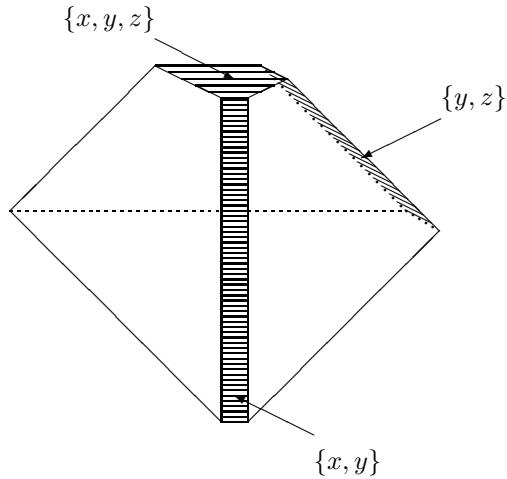
\begin{center}
\begin{picture}(170,170)(0,-10)
\put(80,5){\line(-1,1){80}} \put(90,5){\line(1,1){72.5}}
\put(55,140){\line(-1,-1){55}}

\multiput(97,138)(3,-3){18}{\circle*{.5}}
\put(105,135){\line(1,-1){57.5}}

\multiput(0,85)(3,0){50}{\line(1,0){1}}
\multiput(150,85)(2.5,-1.5){5}{\circle*{.5}}
\put(80,5){\line(0,1){122.5}} \put(90,5){\line(0,1){122.5}}
\put(80,5){\line(1,0){10}}

\put(55,140){\line(1,0){40}} \put(55,140){\line(2,-1){25}}
\put(95,140){\line(2,-1){10}} \put(105,135){\line(-2,-1){15}}
\put(80,127.5){\line(1,0){10}}

\multiput(80,5)(0,2){62}{\line(1,0){10}}

\multiput(95.5,139.5)(3,-3){19}{\line(2,-1){10}}

\put(60,137.5){\line(1,0){40}} \put(65,135){\line(1,0){40}}
\put(70,132.5){\line(1,0){30}} \put(75,130){\line(1,0){20}}

\put(55,150){\vector(2,-1){30}}
\put(55,155){\makebox(0,0)[br]{$\{x,y,z\}$}}

\put(115,0){\vector(-2,1){30}}
\put(115,-5){\makebox(0,0)[tl]{$\{x,y\}$}}

\put(163,120){\vector(-2,-1){30}}
\put(163,125){\makebox(0,0)[bl]{$\{y,z\}$}}

\end{picture}
\end{center}
This polytope, as $\A(H_{4121})$, may be realized as a five-sided
prism.

\vspace{1ex}

\hspace{-2.5em} $(H'_{4211})$\quad As a case where we truncate two
opposite edges and a vertex incident with just one of them, we
have $H'_{4211}=H_{4111}\cup\{\{z,u\}\}$. The picture of a
realization of $\A(H'_{4211})$ is
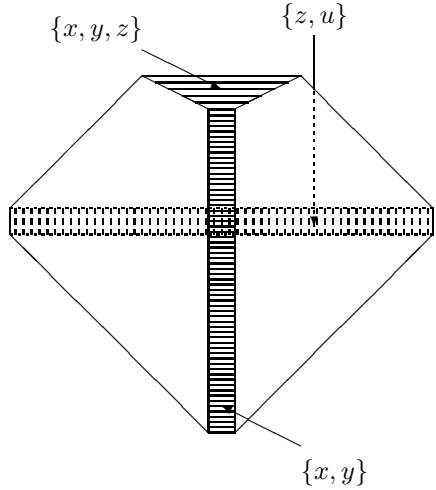
\begin{center}
\begin{picture}(170,175)(0,-10)
\put(80,5){\line(-1,1){75}} \put(90,5){\line(1,1){75}}
\put(55,140){\line(-1,-1){50}} \put(115,140){\line(1,-1){50}}
\multiput(5,90)(3,0){54}{\line(1,0){1}}
\multiput(5,80)(3,0){54}{\line(1,0){1}}
\put(80,5){\line(0,1){122.5}} \put(90,5){\line(0,1){122.5}}
\put(80,5){\line(1,0){10}} \put(5,80){\line(0,1){10}}
\put(165,80){\line(0,1){10}}

\put(55,140){\line(1,0){60}} \put(55,140){\line(2,-1){25}}
\put(115,140){\line(-2,-1){25}} \put(80,127.5){\line(1,0){10}}

\multiput(80,5)(0,2){62}{\line(1,0){10}}
\multiput(7,80)(3,0){53}{\line(0,1){2}}
\multiput(7,84)(3,0){53}{\line(0,1){2}}
\multiput(7,88)(3,0){53}{\line(0,1){2}}

\put(60,137.5){\line(1,0){50}} \put(65,135){\line(1,0){40}}
\put(70,132.5){\line(1,0){30}} \put(75,130){\line(1,0){20}}

\put(55,150){\vector(2,-1){30}}
\put(55,155){\makebox(0,0)[br]{$\{x,y,z\}$}}

\put(115,0){\vector(-2,1){30}}
\put(115,-5){\makebox(0,0)[tl]{$\{x,y\}$}}

\put(120,155){\line(0,-1){20}}
\multiput(120,85)(0,3){19}{\line(0,1){1}}
\put(120,85){\vector(0,-1){2}}
\put(120,160){\makebox(0,0)[b]{$\{z,u\}$}}

\end{picture}
\end{center}
This polytope, as the preceding one, viz.\ $\A(H_{4211})$, and as
$\A(H_{4121})$, may be realized as a five-sided prism.

\vspace{1ex}

\hspace{-2.5em} $(H_{4221},\ldots,H_{4241})$\quad As remaining
cases with two edges truncated we have
\begin{tabbing}
\hspace{1.5em}\=$H_{4221}\:$\=$=H_{4211}\cup\{\{y,z,u\}\}$,\\
\>$H'_{4221}$\>$=H_{4211}\cup\{\{x,z,u\}\}$,\\
\>$H_{4231}$\>$=H_{4221}\cup H'_{4221}$,\\
\>$H'_{4231}$\>$=H_{4221}\cup H_{4121}$,\\
\>$H_{4241}$\>$=H_{4231}\cup H'_{4231}$,
\end{tabbing}
with realizations of the corresponding polytopes having pictures
easily obtained from the preceding ones. (The last six cases are
essentially hypergraphical.)

\vspace{1ex}

\hspace{-2.5em} $(H_{4311})$\quad If we truncate three edges and a
common vertex with which they are all incident, we have the case
of $H_{4311}=H_{4211}\cup \{\{x,z\}\}=H_{4310}\cup
\{\{x,y,z,u\}\}$, with the picture of a realization of
$\A(H_{4311})$ being
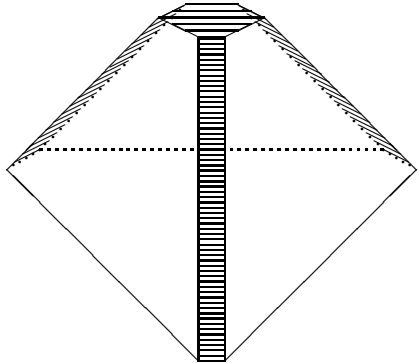
\begin{center}
\begin{picture}(170,147)
\put(80,5){\line(-1,1){72.5}} \put(90,5){\line(1,1){72.5}}

\multiput(73,138)(-3,-3){18}{\circle*{.5}}
\put(65,135){\line(-1,-1){57.5}}

\multiput(97,138)(3,-3){18}{\circle*{.5}}
\put(105,135){\line(1,-1){57.5}}

\multiput(20,85)(3,0){44}{\line(1,0){1}}
\multiput(20,85)(-2.5,-1.5){5}{\circle*{.5}}
\multiput(150,85)(2.5,-1.5){5}{\circle*{.5}}
\put(80,5){\line(0,1){122.5}} \put(90,5){\line(0,1){122.5}}
\put(80,5){\line(1,0){10}}

\put(75,140){\line(1,0){20}} \put(65,135){\line(2,-1){15}}
\put(65,135){\line(2,1){10}} \put(95,140){\line(2,-1){10}}
\put(105,135){\line(-2,-1){15}} \put(80,127.5){\line(1,0){10}}

\multiput(80,5)(0,2){62}{\line(1,0){10}}

\multiput(74.5,139.5)(-3,-3){19}{\line(-2,-1){10}}
\multiput(95.5,139.5)(3,-3){19}{\line(2,-1){10}}

\put(70,137.5){\line(1,0){30}} \put(65,135){\line(1,0){40}}
\put(70,132.5){\line(1,0){30}} \put(75,130){\line(1,0){20}}

\end{picture}
\end{center}
We have here omitted all the labels. This polytope may be realized
also as a six-sided prism.

\vspace{1ex}

\hspace{-2.5em} $(H_{4321},H_{4331},H_{4341})$\quad Next we have
\begin{tabbing}
\hspace{1.5em}\=$H_{4321}\:$\=$=H_{4311}\cup\{\{y,z,u\}\}$,\\
\>$H_{4331}$\>$=H_{4321}\cup\{\{x,z,u\}\}$,\\
\>$H_{4341}$\>$=H_{4331}\cup\{\{x,y,u\}\}$,
\end{tabbing}
with the pictures of realizations of the corresponding polytopes
obtained from the preceding picture by truncating the edges the
base hexagon shares with the shaded quadrilaterals, so as to
obtain one, two or three additional quadrilaterals. These edges
originate from vertices, and we have in fact truncated these
vertices. (The last three cases are essentially hypergraphical.)

\vspace{1ex}

\hspace{-2.5em} $(H'_{4321})$\quad If we truncate the tetrahedron
along a path of three edges and two vertices so as to obtain
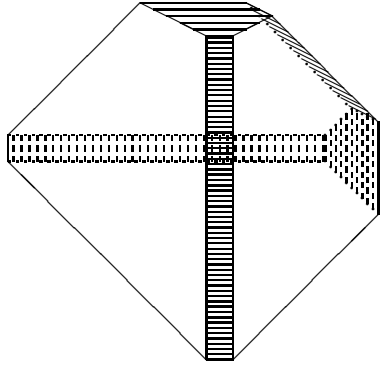
\begin{center}
\begin{picture}(170,147)
\put(80,5){\line(-1,1){75}} \put(90,5){\line(1,1){55}}
\put(55,140){\line(-1,-1){50}}

\multiput(97,138)(3,-3){13}{\circle*{.5}}
\put(105,135){\line(1,-1){40}}

\multiput(125,90)(2,2){5}{\circle*{.5}}
\multiput(125,80)(2,-2){10}{\circle*{.5}}
\multiput(135,100)(2,-1){5}{\circle*{.5}}
\put(145,60){\line(0,1){35}}

\multiput(5,90)(3,0){40}{\line(1,0){1}}
\multiput(5,80)(3,0){40}{\line(1,0){1}}
\put(80,5){\line(0,1){122.5}} \put(90,5){\line(0,1){122.5}}
\put(80,5){\line(1,0){10}} \put(5,80){\line(0,1){10}}
\multiput(125,80)(0,4){3}{\thicklines\line(0,1){2}}

\put(55,140){\line(1,0){40}} \put(55,140){\line(2,-1){25}}
\put(95,140){\line(2,-1){10}} \put(105,135){\line(-2,-1){15}}
\put(80,127.5){\line(1,0){10}}

\multiput(80,5)(0,2){62}{\line(1,0){10}}
\multiput(7,80)(3,0){40}{\line(0,1){2}}
\multiput(7,84)(3,0){40}{\line(0,1){2}}
\multiput(7,88)(3,0){40}{\line(0,1){2}}

\multiput(95,140)(3,-3){14}{\line(2,-1){10}}

\put(60,137.5){\line(1,0){40}} \put(65,135){\line(1,0){40}}
\put(70,132.5){\line(1,0){30}} \put(75,130){\line(1,0){20}}

\multiput(143,62)(0,4){9}{\line(0,1){2}}
\multiput(141,64)(0,4){9}{\line(0,1){2}}
\multiput(139,66)(0,4){8}{\line(0,1){2}}
\multiput(137,68)(0,4){8}{\line(0,1){2}}
\multiput(135,70)(0,4){8}{\line(0,1){2}}
\multiput(133,72)(0,4){7}{\line(0,1){2}}
\multiput(131,74)(0,4){6}{\line(0,1){2}}
\multiput(129,76)(0,4){5}{\line(0,1){2}}
\multiput(127,78)(0,4){4}{\line(0,1){2}}

\end{picture}
\end{center}
we have the picture of a realization of $\A(H'_{4321})$ for
$H'_{4321}=H_{4221}\cup\{\{u,z\}\}$. This polytope is the
three-dimensional associahedron $K_5$ (see \cite{S97a},
\cite{S97b} and \cite{T97}), and $H'_{4321}$ is the saturated
closure~of
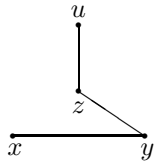
\begin{center}
\begin{picture}(50,65)
\put(0,10){\line(1,0){50}} \put(50,10){\line(-3,2){25}}
\put(25,26.5){\line(0,1){25}}

\put(1,10){\makebox(0,0){\circle*{2}}}
\put(51,10){\makebox(0,0){\circle*{2}}}
\put(26,27){\makebox(0,0){\circle*{2}}}
\put(26,52){\makebox(0,0){\circle*{2}}}

\put(1,7){\makebox(0,0)[t]{$x$}} \put(51,7){\makebox(0,0)[t]{$y$}}
\put(25,23){\makebox(0,0)[t]{$z$}}
\put(25,55){\makebox(0,0)[b]{$u$}}
\end{picture}
\end{center}
(see  $A''$ in Section~3). It is also the saturated closure of the
hypergraph $A$ of Section~3.

This truncation should be compared with Example 5.15 of
\cite{DP10a}, where one reaches $K_5$ not by truncating the
tetrahedron, but from the other end, by collapsing the vertices of
the three-dimensional permutohedron (see $(H_{4641})$ below).

Note that we have truncated in the tetrahedron a chain made of
three edges and two vertices. There is in the tetrahedron a
complementary chain of exactly the same kind, with three edges and
two vertices. Since for the three-dimensional permutohedron we
will truncate all the edges and all the vertices of our
tetrahedron, $K_5$ is located halfway.

To compare our picture of $K_5$ with the picture in \cite{DP10a},
and with what we had in Section~5, we will turn it so as to obtain
the following picture, with facets labelled as before in this
survey, and vertices labelled with s-constructions (the edges are
not labelled):
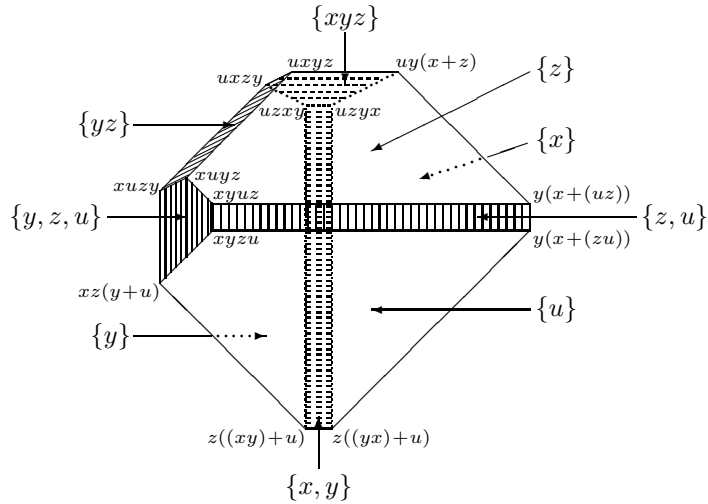
\begin{center}
\begin{picture}(170,170)(0,-10)
\put(80,5){\line(-1,1){55}} \put(90,5){\line(1,1){75}}
\put(115,140){\line(1,-1){50}}

\put(75,140){\line(-1,-1){40}} \put(65,135){\line(-1,-1){40}}

\put(45,90){\line(-1,1){10}} \put(45,80){\line(-1,-1){20}}
\put(35,100){\line(-2,-1){10}} \put(25,60){\line(0,1){35}}

\put(45,90){\line(1,0){120}} \put(45,80){\line(1,0){120}}
\put(45,80){\line(0,1){10}} \put(165,80){\line(0,1){10}}

\multiput(80,5.7)(0,2){61}{\line(0,1){.6}}
\multiput(90,5.7)(0,2){61}{\line(0,1){.6}}
\multiput(80.5,127.5)(2,0){5}{\thicklines\line(1,0){1}}
\put(80,5){\line(1,0){10}}

\put(75,140){\line(1,0){40}}
\multiput(113,139)(-2,-1){12}{\circle*{.5}}
\put(75,140){\line(-2,-1){10}}
\multiput(67,134)(2,-1){7}{\circle*{.5}}

\multiput(75,140)(-3,-3){14}{\line(-2,-1){10}}

\multiput(45,80)(3,0){41}{\line(0,1){10}}

\multiput(80,7)(0,2){60}{\line(1,0){2}}
\multiput(84,7)(0,2){60}{\line(1,0){2}}
\multiput(88,7)(0,2){60}{\line(1,0){2}}

\multiput(70,137.5)(3,0){13}{\line(1,0){1.5}}
\multiput(65,135)(3,0){13}{\line(1,0){1.5}}
\multiput(70,132.5)(3,0){10}{\line(1,0){1.5}}
\multiput(75,130)(3,0){7}{\line(1,0){1.5}}

\put(27,62){\line(0,1){34}} \put(29,64){\line(0,1){33}}
\put(31,66){\line(0,1){32}} \put(33,68){\line(0,1){31}}
\put(35,70){\line(0,1){30}} \put(37,72){\line(0,1){26}}
\put(39,74){\line(0,1){22}} \put(41,76){\line(0,1){18}}
\put(43,78){\line(0,1){14}} \put(44.5,79.5){\line(0,1){11}}

\put(73,142){\makebox(0,0)[bl]{\scriptsize $uxyz$}}
\put(114,142){\makebox(0,0)[bl]{\scriptsize $uy(x\pl z)$}}
\put(64,136){\makebox(0,0)[br]{\scriptsize $uxzy$}}
\put(71,126.5){\makebox(0,0)[t]{\scriptsize $uzxy$}}
\put(100,126.5){\makebox(0,0)[t]{\scriptsize $uzyx$}}
\put(38,99){\makebox(0,0)[bl]{\scriptsize $xuyz$}}
\put(25,95){\makebox(0,0)[br]{\scriptsize $xuzy$}}
\put(45,92){\makebox(0,0)[bl]{\scriptsize $xyuz$}}
\put(166,92){\makebox(0,0)[bl]{\scriptsize $y(x\pl(uz))$}}
\put(45,78){\makebox(0,0)[tl]{\scriptsize $xyzu$}}
\put(166,80){\makebox(0,0)[tl]{\scriptsize $y(x\pl(zu))$}}
\put(25,60){\makebox(0,0)[tr]{\scriptsize $xz(y\pl u)$}}
\put(80,5){\makebox(0,0)[tr]{\scriptsize $z((xy)\pl u)$}}
\put(90,5){\makebox(0,0)[tl]{\scriptsize $z((yx)\pl u)$}}

\put(95,137.5){\vector(0,-1){2}} \put(95,155){\line(0,-1){14.5}}
\put(95,158){\makebox(0,0)[b]{$\{xyz\}$}}

\put(15,120){\vector(1,0){38}}
\put(14,120){\makebox(0,0)[r]{$\{yz\}$}}

\put(165,140){\vector(-2,-1){60}}
\put(167,141){\makebox(0,0)[l]{$\{z\}$}}

\put(125,100){\vector(-3,-1){2}}
\multiput(125,100)(3,1){8}{\circle*{.5}}
\put(149,108){\line(3,1){15}}
\put(167,114){\makebox(0,0)[l]{$\{x\}$}}

\put(5,85){\vector(1,0){30}}
\put(3,85){\makebox(0,0)[r]{$\{y,z,u\}$}}

\put(205,85){\vector(-1,0){60}}
\put(207,85){\makebox(0,0)[l]{$\{z,u\}$}}

\put(165,50){\vector(-1,0){60}}
\put(167,50){\makebox(0,0)[l]{$\{u\}$}}

\put(85,7.5){\vector(0,1){2}} \put(85,-10){\line(0,1){14.5}}
\put(85,-12){\makebox(0,0)[t]{$\{x,y\}$}}

\put(63,40){\vector(1,0){2}}
\multiput(47,40)(3,0){5}{\circle*{.5}}
\put(15,40){\line(1,0){29.5}}
\put(14,40){\makebox(0,0)[r]{$\{y\}$}}

\end{picture}
\end{center}

\vspace{1ex}

\hspace{-2.5em} $(H'_{4331},H'_{4341},H^*_{4331})$\quad Next we
have three more cases with three edges truncated:
\begin{tabbing}
\hspace{1.5em}\=$H'_{4331}\:$\=$=H'_{4321}\cup\{\{x,z,u\}\}$,\\
\>$H'_{4341}$\>$=H'_{4331}\cup\{\{x,y,u\}\}$,\\
\>$H^*_{4331}$\>$=H_{4031}\cup\{\{x,z\},\{y,z\},\{u,z\}\}$.
\end{tabbing}

The first of these cases, for which we will draw no picture, is
interesting because $H'_{4331}$ is the intersection of the
hypergraphs of the hemiassociahedron (see $(H_{4431})$ below) and
of the three-dimensional cyclohedron (see $(H^\circ_{4441})$
below). We draw no picture for the second case either. Both of
these pictures that we do not draw are obtained easily from our
picture of $K_5$; we just extend with vertices the initial path of
truncated edges and vertices of the tetrahedron. (These two cases
are essentially hypergraphical.)

In the third case, we have that $H^*_{4331}$ is the saturated
closure~of
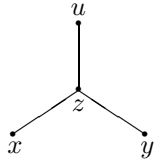
\begin{center}
\begin{picture}(50,55)(0,5)
\put(0,10){\line(3,2){25}} \put(50,10){\line(-3,2){25}}
\put(25,26.5){\line(0,1){25}}

\put(1,10){\makebox(0,0){\circle*{2}}}
\put(51,10){\makebox(0,0){\circle*{2}}}
\put(26,27){\makebox(0,0){\circle*{2}}}
\put(26,52){\makebox(0,0){\circle*{2}}}

\put(1,7){\makebox(0,0)[t]{$x$}} \put(51,7){\makebox(0,0)[t]{$y$}}
\put(25,23){\makebox(0,0)[t]{$z$}}
\put(25,55){\makebox(0,0)[b]{$u$}}
\end{picture}
\end{center}
(see $A^*$ in Section~3). The picture of a realization of
$\A(H^*_{4331})$, turned in the manner of the preceding picture,
is
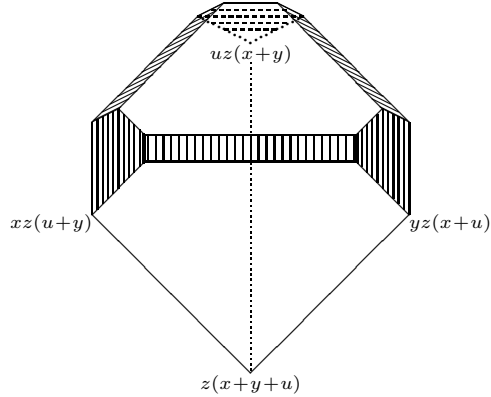
\begin{center}
\begin{picture}(170,150)(0,-5)
\put(85,0){\line(-1,1){60}} \put(85,0){\line(1,1){60}}

\put(75,140){\line(-1,-1){40}} \put(95,140){\line(1,-1){40}}
\put(65,135){\line(-1,-1){40}} \put(105,135){\line(1,-1){40}}

\put(45,90){\line(-1,1){10}} \put(45,80){\line(-1,-1){20}}
\put(35,100){\line(-2,-1){10}} \put(25,60){\line(0,1){35}}

\put(125,90){\line(1,1){10}} \put(125,80){\line(1,-1){20}}
\put(135,100){\line(2,-1){10}} \put(145,60){\line(0,1){35}}

\put(45,90){\line(1,0){80}} \put(45,80){\line(1,0){80}}
\put(45,80){\line(0,1){10}} \put(125,80){\line(0,1){10}}

\multiput(85,1)(0,2){59}{\line(0,1){.6}}

\put(75,140){\line(1,0){20}} \put(75,140){\line(-2,-1){10}}
\put(95,140){\line(2,-1){10}}
\multiput(67,134)(2,-1){10}{\circle*{.5}}
\multiput(103,134)(-2,-1){9}{\circle*{.5}}

\multiput(75,140)(-3,-3){14}{\line(-2,-1){10}}
\multiput(95,140)(3,-3){14}{\line(2,-1){10}}

\multiput(46,80)(3,0){27}{\line(0,1){10}}

\multiput(70,137.5)(3,0){10}{\line(1,0){1.5}}
\multiput(65,135)(2.95,0){14}{\line(1,0){1.5}}
\multiput(70,132.5)(3,0){10}{\line(1,0){1.5}}
\multiput(75,130)(3,0){7}{\line(1,0){1.5}}

\put(27,62){\line(0,1){34}} \put(29,64){\line(0,1){33}}
\put(31,66){\line(0,1){32}} \put(33,68){\line(0,1){31}}
\put(35,70){\line(0,1){30}} \put(37,72){\line(0,1){26}}
\put(39,74){\line(0,1){22}} \put(41,76){\line(0,1){18}}
\put(43,78){\line(0,1){14}} \put(44.5,79.5){\line(0,1){11}}

\put(143,62){\line(0,1){34}} \put(141,64){\line(0,1){33}}
\put(139,66){\line(0,1){32}} \put(137,68){\line(0,1){31}}
\put(135,70){\line(0,1){30}} \put(133,72){\line(0,1){26}}
\put(131,74){\line(0,1){22}} \put(129,76){\line(0,1){18}}
\put(127,78){\line(0,1){14}} \put(125.5,79.5){\line(0,1){11}}

\put(85,124){\makebox(0,0)[t]{\scriptsize $uz(x\pl y)$}}
\put(25,60){\makebox(0,0)[tr]{\scriptsize $xz(u\pl y)$}}
\put(145,60){\makebox(0,0)[tl]{\scriptsize $yz(x\pl u)$}}
\put(85,-1){\makebox(0,0)[t]{\scriptsize $z(x\pl y\pl u)$}}

\end{picture}
\end{center}
Here we label with s-constructions just those vertices whose
s-constructions are not permutations of $x$, $y$, $z$ and $u$.
This picture permits a comparison with the picture of Example 5.16
of \cite{DP10a}. The pictured polytope is the three-dimensional
stellohedron of \cite{PRW08} (Section 10.4; see also \cite{P09},
Section 8.4). This polytope is called $D4$ in \cite{ACDELM09}
(Figure 17), and in \cite{DP10a} the entirely Greek name
\emph{astrohedron} is suggested.

\vspace{1ex}

\hspace{-2.5em} $(H^*_{4341})$\quad Next we have
$H^*_{4341}=H^*_{4331}\cup\{\{x,y,u\}\}$, with the picture for a
realization of $\A(H^*_{4341})$ obtained from the preceding
picture, given for $\A(H^*_{4331})$, by truncating the vertex
labelled $z(x\pl y\pl u)$. (This is an essentially hypergraphical
case.)

\vspace{1ex}

\hspace{-2.5em} $(H_{4431})$\quad Next we have
$H_{4431}=H_{4331}\cup\{\{z,u\}\}=
H'_{4331}\cup\{\{x,z\}\}=H^*_{4331}\cup\{\{x,y\}\}$, which is the
saturated closure of
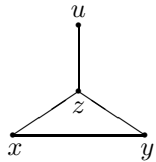
\begin{center}
\begin{picture}(50,61)
\put(0,10){\line(1,0){50}} \put(0,10){\line(3,2){25}}
\put(50,10){\line(-3,2){25}} \put(25,26.5){\line(0,1){25}}

\put(1,10){\makebox(0,0){\circle*{2}}}
\put(51,10){\makebox(0,0){\circle*{2}}}
\put(26,27){\makebox(0,0){\circle*{2}}}
\put(26,52){\makebox(0,0){\circle*{2}}}

\put(1,7){\makebox(0,0)[t]{$x$}} \put(51,7){\makebox(0,0)[t]{$y$}}
\put(25,23){\makebox(0,0)[t]{$z$}}
\put(25,55){\makebox(0,0)[b]{$u$}}
\end{picture}
\end{center}
The picture for a realization of $\A(H_{4431})$, turned in the
manner of the preceding two pictures, and labelled in the manner
of the preceding one, is
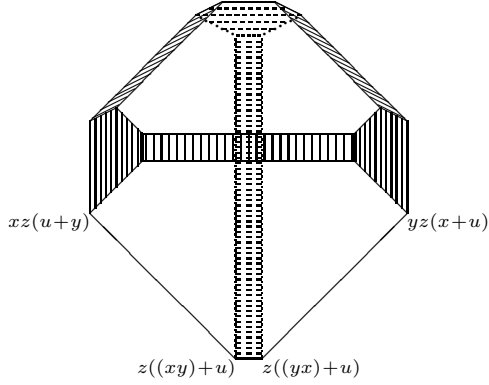
\begin{center}
\begin{picture}(170,145)
\put(80,5){\line(-1,1){55}} \put(90,5){\line(1,1){55}}

\put(75,140){\line(-1,-1){40}} \put(95,140){\line(1,-1){40}}
\put(65,135){\line(-1,-1){40}} \put(105,135){\line(1,-1){40}}

\put(45,90){\line(-1,1){10}} \put(45,80){\line(-1,-1){20}}
\put(35,100){\line(-2,-1){10}} \put(25,60){\line(0,1){35}}

\put(125,90){\line(1,1){10}} \put(125,80){\line(1,-1){20}}
\put(135,100){\line(2,-1){10}} \put(145,60){\line(0,1){35}}

\put(45,90){\line(1,0){80}} \put(45,80){\line(1,0){80}}
\put(45,80){\line(0,1){10}} \put(125,80){\line(0,1){10}}

\multiput(80,5.7)(0,2){61}{\line(0,1){.6}}
\multiput(90,5.7)(0,2){61}{\line(0,1){.6}}
\multiput(80.5,127.5)(2,0){5}{\thicklines\line(1,0){1}}
\put(80,5){\line(1,0){10}}

\put(75,140){\line(1,0){20}} \put(75,140){\line(-2,-1){10}}
\put(95,140){\line(2,-1){10}}
\multiput(67,134)(2,-1){7}{\circle*{.5}}
\multiput(103,134)(-2,-1){7}{\circle*{.5}}

\multiput(75,140)(-3,-3){14}{\line(-2,-1){10}}
\multiput(95,140)(3,-3){14}{\line(2,-1){10}}

\multiput(46,80)(3,0){27}{\line(0,1){10}}

\multiput(80,7)(0,2){60}{\line(1,0){2}}
\multiput(84,7)(0,2){60}{\line(1,0){2}}
\multiput(88,7)(0,2){60}{\line(1,0){2}}

\multiput(70,137.5)(3,0){10}{\line(1,0){1.5}}
\multiput(65,135)(2.95,0){14}{\line(1,0){1.5}}
\multiput(70,132.5)(3,0){10}{\line(1,0){1.5}}
\multiput(75,130)(3,0){7}{\line(1,0){1.5}}

\put(27,62){\line(0,1){34}} \put(29,64){\line(0,1){33}}
\put(31,66){\line(0,1){32}} \put(33,68){\line(0,1){31}}
\put(35,70){\line(0,1){30}} \put(37,72){\line(0,1){26}}
\put(39,74){\line(0,1){22}} \put(41,76){\line(0,1){18}}
\put(43,78){\line(0,1){14}} \put(44.5,79.5){\line(0,1){11}}

\put(143,62){\line(0,1){34}} \put(141,64){\line(0,1){33}}
\put(139,66){\line(0,1){32}} \put(137,68){\line(0,1){31}}
\put(135,70){\line(0,1){30}} \put(133,72){\line(0,1){26}}
\put(131,74){\line(0,1){22}} \put(129,76){\line(0,1){18}}
\put(127,78){\line(0,1){14}} \put(125.5,79.5){\line(0,1){11}}

\put(25,60){\makebox(0,0)[tr]{\scriptsize $xz(u\pl y)$}}
\put(145,60){\makebox(0,0)[tl]{\scriptsize $yz(x\pl u)$}}
\put(80,5){\makebox(0,0)[tr]{\scriptsize $z((xy)\pl u)$}}
\put(90,5){\makebox(0,0)[tl]{\scriptsize $z((yx)\pl u)$}}

\end{picture}
\end{center}
This picture permits a comparison with the picture of Example 5.14
of \cite{DP10a}, where the pictured polytope was called
\emph{hemiassociahedron} (it is called $X^a_4$ in \cite{ACDELM09},
Figure 17, and $P_{1,2}$ in \cite{B09}, Figure 6).

\vspace{1ex}

\hspace{-2.5em} $(H_{4441})$\quad Next we have
$H_{4441}=H_{4431}\cup\{\{x,y,u\}\}$, with the picture for a
realization of $\A(H_{4441})$ obtained from the preceding picture
of the hemiassociahedron, given for $\A(H_{4431})$, by truncating
the bottom edge, so as to obtain a quadrilateral. This edge
originates from a vertex, and in fact we truncate this vertex.
(This is an essentially hypergraphical case.)

\vspace{1ex}

\hspace{-2.5em} $(H^\circ_{4441})$\quad Next we have
$H^\circ_{4441}=H'_{4341}\cup\{\{x,u\}\}$, which is the saturated
closure~of
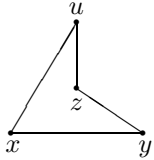
\begin{center}
\begin{picture}(50,60)(100,0)
\put(100,10){\line(1,0){50}} \put(100,10){\line(3,5){25}}
\put(150,10){\line(-3,2){25}} \put(125,26.5){\line(0,1){25}}

\put(101,10){\makebox(0,0){\circle*{2}}}
\put(151,10){\makebox(0,0){\circle*{2}}}
\put(126,27){\makebox(0,0){\circle*{2}}}
\put(126,52){\makebox(0,0){\circle*{2}}}

\put(101,7){\makebox(0,0)[t]{$x$}}
\put(151,7){\makebox(0,0)[t]{$y$}}
\put(125,23){\makebox(0,0)[t]{$z$}}
\put(125,55){\makebox(0,0)[b]{$u$}}
\end{picture}
\end{center}
(see $A^\circ$ in Section~3). The picture for a realization of
$\A(H^\circ_{4441})$, turned in the manner of the preceding three
pictures, and with all labels omitted, is
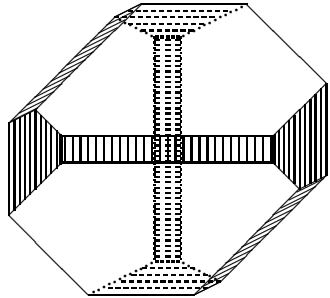
\begin{center}
\begin{picture}(170,110)(0,35)
\put(55,30){\line(-1,1){30}} \put(115,140){\line(1,-1){30}}

\put(145,75){\line(-1,-1){40}} \put(135,70){\line(-1,-1){40}}

\put(135,70){\line(-1,1){10}} \put(145,110){\line(-1,-1){20}}
\put(145,75){\line(-2,-1){10}} \put(145,75){\line(0,1){35}}

\put(75,140){\line(-1,-1){40}} \put(65,135){\line(-1,-1){40}}

\put(45,90){\line(-1,1){10}} \put(45,80){\line(-1,-1){20}}
\put(35,100){\line(-2,-1){10}} \put(25,60){\line(0,1){35}}

\put(45,90){\line(1,0){80}} \put(45,80){\line(1,0){80}}
\put(45,80){\line(0,1){10}} \put(125,80){\line(0,1){10}}

\multiput(80,43.2)(0,2){43}{\line(0,1){.6}}
\multiput(90,43.2)(0,2){43}{\line(0,1){.6}}
\multiput(80.5,127.5)(2,0){5}{\thicklines\line(1,0){1}}
\multiput(80.5,42.5)(2,0){5}{\thicklines\line(1,0){1}}

\put(55,30){\line(1,0){40}}
\multiput(57,31)(2,1){12}{\circle*{.5}}
\put(95,30){\line(2,1){10}}
\multiput(103,36)(-2,1){7}{\circle*{.5}}

\put(75,140){\line(1,0){40}}
\multiput(113,139)(-2,-1){12}{\circle*{.5}}
\put(75,140){\line(-2,-1){10}}
\multiput(67,134)(2,-1){7}{\circle*{.5}}

\multiput(75,140)(-3,-3){14}{\line(-2,-1){10}}
\multiput(95,30)(3,3){14}{\line(2,1){10}}

\multiput(46,80)(3,0){27}{\line(0,1){10}}

\multiput(80,44.5)(0,2){42}{\line(1,0){2}}
\multiput(84,44.5)(0,2){42}{\line(1,0){2}}
\multiput(88,44.5)(0,2){42}{\line(1,0){2}}

\multiput(100,32.5)(-3,0){13}{\line(-1,0){1.5}}
\multiput(105,35)(-3,0){13}{\line(-1,0){1.5}}
\multiput(100,37.5)(-3,0){10}{\line(-1,0){1.5}}
\multiput(95,40)(-3,0){7}{\line(-1,0){1.5}}

\multiput(70,137.5)(3,0){13}{\line(1,0){1.5}}
\multiput(65,135)(3,0){13}{\line(1,0){1.5}}
\multiput(70,132.5)(3,0){10}{\line(1,0){1.5}}
\multiput(75,130)(3,0){7}{\line(1,0){1.5}}

\put(27,62){\line(0,1){34}} \put(29,64){\line(0,1){33}}
\put(31,66){\line(0,1){32}} \put(33,68){\line(0,1){31}}
\put(35,70){\line(0,1){30}} \put(37,72){\line(0,1){26}}
\put(39,74){\line(0,1){22}} \put(41,76){\line(0,1){18}}
\put(43,78){\line(0,1){14}} \put(44.5,79.5){\line(0,1){11}}

\put(143,108){\line(0,-1){34}} \put(141,106){\line(0,-1){33}}
\put(139,104){\line(0,-1){32}} \put(137,102){\line(0,-1){31}}
\put(135,100){\line(0,-1){30}} \put(133,98){\line(0,-1){26}}
\put(131,96){\line(0,-1){22}} \put(129,94){\line(0,-1){18}}
\put(127,92){\line(0,-1){14}} \put(125.5,90.5){\line(0,-1){11}}

\end{picture}
\end{center}
This picture permits a comparison with the picture of Example 5.13
of \cite{DP10a}. The pictured polytope is the three-dimensional
cyclohedron (see \cite{BT94} and \cite{S97b}, Section~4).

\vspace{1ex}

\hspace{-2.5em} $(H_{4541})$\quad Next we have
$H_{4541}=H^\circ_{4441}\cup\{\{x,z\}\}$, which is the saturated
closure~of
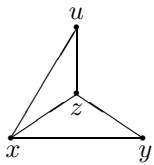
\begin{center}
\begin{picture}(50,65)
\put(0,10){\line(1,0){50}} \put(0,10){\line(3,5){25}}
\put(0,10){\line(3,2){25}} \put(50,10){\line(-3,2){25}}
\put(25,26.5){\line(0,1){25}}

\put(1,10){\makebox(0,0){\circle*{2}}}
\put(51,10){\makebox(0,0){\circle*{2}}}
\put(26,27){\makebox(0,0){\circle*{2}}}
\put(26,52){\makebox(0,0){\circle*{2}}}

\put(1,7){\makebox(0,0)[t]{$x$}} \put(51,7){\makebox(0,0)[t]{$y$}}
\put(25,23){\makebox(0,0)[t]{$z$}}
\put(25,55){\makebox(0,0)[b]{$u$}}
\end{picture}
\end{center}
The picture for a realization of $\A(H_{4541})$ is obtained from
the preceding picture of the three-dimensional cyclohedron, given
for $\A(H_{4441})$, by truncating one more edge---in this case,
the edge $\{\{x\},\{z\},\{x,y,z,u\}\}$, which is now in the
north-east. This picture permits a comparison with the picture of
Example 5.12 of \cite{DP10a}, where the pictured polytope is
called \emph{hemicyclohedron} (see also \cite{FS09}, Figure 10).

\vspace{1ex}

\hspace{-2.5em} $(H_{4641})$\quad Finally, we have
$H_{4641}=H_{4541}\cup\{\{y,z\}\}$, which is the saturated
closure~of
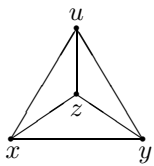
\begin{center}
\begin{picture}(50,65)
\put(0,10){\line(1,0){50}} \put(0,10){\line(3,5){25}}
\put(0,10){\line(3,2){25}} \put(50,10){\line(-3,5){25}}
\put(50,10){\line(-3,2){25}} \put(25,26.5){\line(0,1){25}}

\put(1,10){\makebox(0,0){\circle*{2}}}
\put(51,10){\makebox(0,0){\circle*{2}}}
\put(26,27){\makebox(0,0){\circle*{2}}}
\put(26,52){\makebox(0,0){\circle*{2}}}

\put(1,7){\makebox(0,0)[t]{$x$}} \put(51,7){\makebox(0,0)[t]{$y$}}
\put(25,23){\makebox(0,0)[t]{$z$}}
\put(25,55){\makebox(0,0)[b]{$u$}}
\end{picture}
\end{center}
Here we have selected for truncation all the vertices and all the
edges of the tetrahedron. The picture for a realization of
$\A(H_{4641})$, turned in the manner of the preceding four
pictures, is
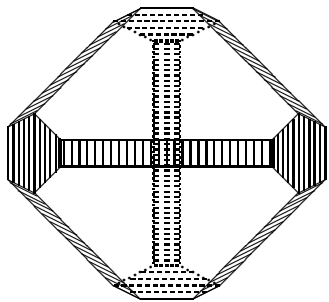
\begin{center}
\begin{picture}(170,110)(0,30)

\put(25,75){\line(1,-1){40}} \put(35,70){\line(1,-1){40}}
\put(145,75){\line(-1,-1){40}} \put(135,70){\line(-1,-1){40}}

\put(135,70){\line(-1,1){10}} \put(135,100){\line(-1,-1){10}}
\put(145,75){\line(-2,-1){10}} \put(145,95){\line(-2,1){10}}
\put(145,75){\line(0,1){20}}

\put(75,140){\line(-1,-1){40}} \put(95,140){\line(1,-1){40}}
\put(65,135){\line(-1,-1){40}} \put(105,135){\line(1,-1){40}}

\put(45,90){\line(-1,1){10}} \put(45,80){\line(-1,-1){10}}
\put(35,100){\line(-2,-1){10}} \put(35,70){\line(-2,1){10}}
\put(25,75){\line(0,1){20}}

\put(45,90){\line(1,0){80}} \put(45,80){\line(1,0){80}}
\put(45,80){\line(0,1){10}} \put(125,80){\line(0,1){10}}

\multiput(80,43.2)(0,2){43}{\line(0,1){.6}}
\multiput(90,43.2)(0,2){43}{\line(0,1){.6}}
\multiput(80.5,127.5)(2,0){5}{\thicklines\line(1,0){1}}
\multiput(80.5,42.5)(2,0){5}{\thicklines\line(1,0){1}}

\put(75,30){\line(1,0){20}} \multiput(67,36)(2,1){7}{\circle*{.5}}
\put(75,30){\line(-2,1){10}} \put(95,30){\line(2,1){10}}
\multiput(103,36)(-2,1){7}{\circle*{.5}}

\put(75,140){\line(1,0){20}} \put(75,140){\line(-2,-1){10}}
\put(95,140){\line(2,-1){10}}
\multiput(67,134)(2,-1){7}{\circle*{.5}}
\multiput(103,134)(-2,-1){7}{\circle*{.5}}

\multiput(75,140)(-3,-3){14}{\line(-2,-1){10}}
\multiput(95,140)(3,-3){14}{\line(2,-1){10}}
\multiput(75,30)(-3,3){14}{\line(-2,1){10}}
\multiput(95,30)(3,3){14}{\line(2,1){10}}

\multiput(46,80)(3,0){27}{\line(0,1){10}}

\multiput(80,44.5)(0,2){42}{\line(1,0){2}}
\multiput(84,44.5)(0,2){42}{\line(1,0){2}}
\multiput(88,44.5)(0,2){42}{\line(1,0){2}}

\multiput(100,32.5)(-3,0){10}{\line(-1,0){1.5}}
\multiput(105,35)(-2.95,0){14}{\line(-1,0){1.5}}
\multiput(100,37.5)(-3,0){10}{\line(-1,0){1.5}}
\multiput(95,40)(-3,0){7}{\line(-1,0){1.5}}

\multiput(70,137.5)(3,0){10}{\line(1,0){1.5}}
\multiput(65,135)(2.95,0){14}{\line(1,0){1.5}}
\multiput(70,132.5)(3,0){10}{\line(1,0){1.5}}
\multiput(75,130)(3,0){7}{\line(1,0){1.5}}

\put(27,74){\line(0,1){22}} \put(29,73){\line(0,1){24}}
\put(31,72){\line(0,1){26}} \put(33,71){\line(0,1){28}}
\put(35,70){\line(0,1){30}} \put(37,72){\line(0,1){26}}
\put(39,74){\line(0,1){22}} \put(41,76){\line(0,1){18}}
\put(43,78){\line(0,1){14}} \put(44.5,79.5){\line(0,1){11}}

\put(143,96){\line(0,-1){22}} \put(141,97){\line(0,-1){24}}
\put(139,98){\line(0,-1){26}} \put(137,99){\line(0,-1){28}}
\put(135,100){\line(0,-1){30}} \put(133,98){\line(0,-1){26}}
\put(131,96){\line(0,-1){22}} \put(129,94){\line(0,-1){18}}
\put(127,92){\line(0,-1){14}} \put(125.5,90.5){\line(0,-1){11}}

\end{picture}
\end{center}
This picture permits a comparison with the picture of Example 5.11
of \cite{DP10a}. The pictured polytope is the three-dimensional
permutohedron (see \cite{Z95}, Lecture~0, Example 0.10, and
\cite{GR63}). All its vertices are constructions of the type
$\{\{x\},\{x,y\},\{x,y,z\},\{x,y,z,u\}\}$; the corresponding
s-construction is $uzyx$, i.e.\ a permutation of $x$, $y$, $z$ and
$u$.

As for every $k\geq 3$ we have that $\A(H_{k0\ldots01})$ may be
realized as the ${(k\mn 1)}$-dimensional simplex, so, at the other
end, with $\A(H_{kn_1\ldots n_m1})$, where all the subscripts
$n_1\ldots n_m$ are maximal (i.e.\ where $n_i={k\choose {i+1}}$),
we obtain the ${(k\mn 1)}$-dimensional permutohedron. In between,
with lesser values of $n_1,\ldots, n_m$, but greater than the
minimal values $0,\ldots,0$, we obtain the ${(k\mn
1)}$-dimensional associahedron (the corresponding graph is a path
of ${k\mn 1}$ edges), the ${(k\mn 1)}$-dimensional cyclohedron
(the corresponding graph is a cycle of $k$ edges and $k$
vertices), and the ${(k\mn 1)}$-dimensional astrohedron (the
corresponding graph is a star-like graph with one vertex in the
middle and ${k\mn 1}$ vertices around joined by ${k\mn 1}$ edges).

With $k=3$, and dimension 2, the associahedron coincides with the
astro\-hedron---both are the pentagon---and the cyclohedron with
the permutohedron---both are the hexagon (see the cases of
$H_{321}$ and $H_{331}$ above). In the degenerate case when $k$ is
$2$, the simplex, the associahedron, the astrohedron and the
permutohedron of dimension 1 all coincide; they are all a single
edge with two incident vertices, which is the only polytope of
dimension 1 (see the case of $H_{21}$ above). If we take as in
\cite{H69} (Chapter 2) that a graph which is a cycle must have at
least 3 vertices, then there is no one-dimensional cyclohedron;
but we may stipulate by convention, as in \cite{S97b} (Section~4),
that the one-dimensional cyclohedron is also a single edge with
two incident vertices, and hence it coincides with the others. In
the degenerate case when $k$ is $1$, we have again just one
polytope of dimension 0; namely, a single vertex.

At the end, we give a chart of the types of hypergraphs
corresponding to some of the polytopes encountered in this
section, including those that are more interesting. (A chart with
all the types would be too intricate.) A line in this chart is
drawn when a hypergraph of one type is included in a hypergraph of
another type. The labels of types in boxes are those of cases
covered previously in \cite{CD06}, \cite{D09} and \cite{DP10a}
(these cases are not essentially hypergraphical). We have made
complete the upper part of the chart above $H^*_{4331}$,
$H_{4311}$ and $H'_{4321}$, which involves the truncation of at
least three edges. If $H_{4311}$ and the seven points with labels
not in boxes are omitted from this part of the chart, then we
obtain (upside down) the chart of \cite{DP10a} (after Example
5.16).

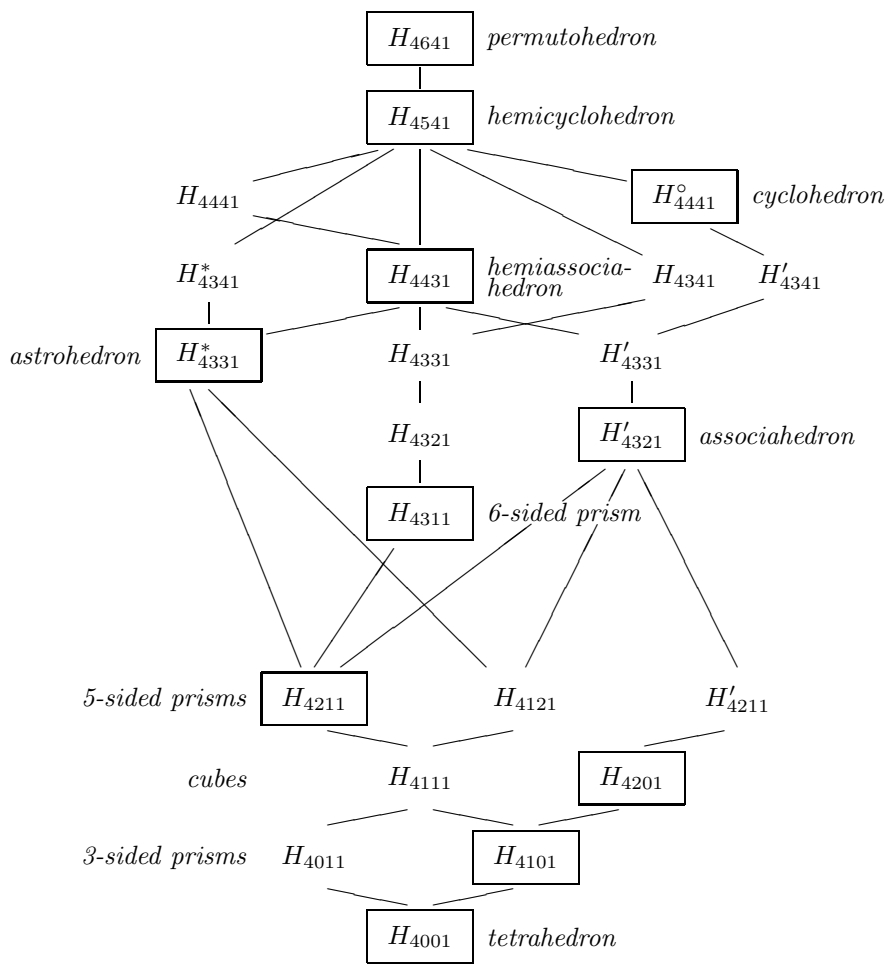
\begin{center}
\begin{picture}(320,360)(10,0)
\put(140,0){\line(1,0){40}} \put(140,20){\line(1,0){40}}
\put(140,0){\line(0,1){20}} \put(180,0){\line(0,1){20}}
\put(160,10){\makebox(0,0){$H_{4001}$}}
\put(185,10){\makebox(0,0)[l]{\emph{tetrahedron}}}

\put(120,40){\makebox(0,0){$H_{4011}$}}
\put(95,40){\makebox(0,0)[r]{\emph{3-sided prisms}}}

\put(180,30){\line(1,0){40}} \put(180,50){\line(1,0){40}}
\put(180,30){\line(0,1){20}} \put(220,30){\line(0,1){20}}
\put(200,40){\makebox(0,0){$H_{4101}$}}

\put(160,70){\makebox(0,0){$H_{4111}$}}
\put(95,70){\makebox(0,0)[r]{\emph{cubes}}}

\put(220,60){\line(1,0){40}} \put(220,80){\line(1,0){40}}
\put(220,60){\line(0,1){20}} \put(260,60){\line(0,1){20}}
\put(240,70){\makebox(0,0){$H_{4201}$}}

\put(100,90){\line(1,0){40}} \put(100,110){\line(1,0){40}}
\put(100,90){\line(0,1){20}} \put(140,90){\line(0,1){20}}
\put(120,100){\makebox(0,0){$H_{4211}$}}
\put(95,100){\makebox(0,0)[r]{\emph{5-sided prisms}}}

\put(200,100){\makebox(0,0){$H_{4121}$}}

\put(280,100){\makebox(0,0){$H'_{4211}$}}

\put(140,160){\line(1,0){40}} \put(140,180){\line(1,0){40}}
\put(140,160){\line(0,1){20}} \put(180,160){\line(0,1){20}}
\put(160,170){\makebox(0,0){$H_{4311}$}}
\put(185,170){\makebox(0,0)[l]{\emph{6-sided prism}}}

\put(160,200){\makebox(0,0){$H_{4321}$}}

\put(220,190){\line(1,0){40}} \put(220,210){\line(1,0){40}}
\put(220,190){\line(0,1){20}} \put(260,190){\line(0,1){20}}
\put(240,200){\makebox(0,0){$H'_{4321}$}}
\put(265,200){\makebox(0,0)[l]{\emph{associahedron}}}

\put(60,220){\line(1,0){40}} \put(60,240){\line(1,0){40}}
\put(60,220){\line(0,1){20}} \put(100,220){\line(0,1){20}}
\put(80,230){\makebox(0,0){$H^*_{4331}$}}
\put(55,230){\makebox(0,0)[r]{\emph{astrohedron}}}

\put(160,230){\makebox(0,0){$H_{4331}$}}

\put(240,230){\makebox(0,0){$H'_{4331}$}}

\put(80,260){\makebox(0,0){$H^*_{4341}$}}

\put(140,250){\line(1,0){40}} \put(140,270){\line(1,0){40}}
\put(140,250){\line(0,1){20}} \put(180,250){\line(0,1){20}}
\put(160,260){\makebox(0,0){$H_{4431}$}}
\put(185,264){\makebox(0,0)[l]{\emph{hemiassocia-}}}
\put(185,256){\makebox(0,0)[l]{\emph{hedron}}}

\put(260,260){\makebox(0,0){$H_{4341}$}}

\put(300,260){\makebox(0,0){$H'_{4341}$}}

\put(80,290){\makebox(0,0){$H_{4441}$}}

\put(240,280){\line(1,0){40}} \put(240,300){\line(1,0){40}}
\put(240,280){\line(0,1){20}} \put(280,280){\line(0,1){20}}
\put(260,290){\makebox(0,0){$H^\circ_{4441}$}}
\put(285,290){\makebox(0,0)[l]{\emph{cyclohedron}}}

\put(140,310){\line(1,0){40}} \put(140,330){\line(1,0){40}}
\put(140,310){\line(0,1){20}} \put(180,310){\line(0,1){20}}
\put(160,320){\makebox(0,0){$H_{4541}$}}
\put(185,320){\makebox(0,0)[l]{\emph{hemicyclohedron}}}

\put(140,340){\line(1,0){40}} \put(140,360){\line(1,0){40}}
\put(140,340){\line(0,1){20}} \put(180,340){\line(0,1){20}}
\put(160,350){\makebox(0,0){$H_{4641}$}}
\put(185,350){\makebox(0,0)[l]{\emph{permutohedron}}}


\put(155,22){\line(-5,1){30}} \put(165,22){\line(5,1){30}}
\put(125,52){\line(5,1){30}} \put(195,52){\line(-5,1){30}}
\put(205,52){\line(5,1){30}} \put(245,82){\line(5,1){30}}
\put(155,82){\line(-5,1){30}} \put(165,82){\line(5,1){30}}
\put(115,112){\line(-2,5){42}} \put(120,112){\line(2,3){30}}
\put(130,112){\line(4,3){72}} \put(214,175){\line(4,3){16}}
\put(185,112){\line(-1,1){105}} \put(200,112){\line(1,2){27}}
\put(231.5,175){\line(1,2){6}} \put(280,112){\line(-1,2){37.5}}
\put(160,182){\line(0,1){8}} \put(160,212){\line(0,1){8}}
\put(160,240){\line(0,1){8}} \put(80,242){\line(0,1){8}}
\put(240,212){\line(0,1){8}} \put(220,238){\line(-5,1){50}}
\put(250,238){\line(3,1){40}} \put(180,238){\line(5,1){65}}
\put(102,238){\line(5,1){50}} \put(152,272){\line(-5,1){55}}
\put(90,272){\line(5,3){60}} \put(160,272){\line(0,1){36}}
\put(244,268){\line(-2,1){80}} \put(160,331){\line(0,1){8}}
\put(238,296){\line(-5,1){60}} \put(290,268){\line(-2,1){20}}
\put(97,296){\line(4,1){47}}

\end{picture}
\end{center}

\vspace{3ex}

\noindent {\small \emph{Acknowledgement.} We would like to thank
Sonja \v Cuki\' c for leading us to the references mentioned in
the last paragraph of the Introduction. This work was supported by
the Ministry of Science of Serbia (Grant ON174026).}

\end{document}